# Topics In Primitive Roots

N. A. Carella

***Abstract***: This monograph considers a few topics in the theory of primitive roots modulo a prime $p \geq 2$. A few estimates of the least primitive roots $g(p)$ and the least prime primitive roots $g^*(p)$ modulo $p$, a large prime, are determined. One of the estimate here seems to sharpen the Burgess estimate $g(p) \ll p^{1/4+\epsilon}$ for arbitrarily small number $\epsilon > 0$, to the smaller estimate $g(p) \ll p^{5/\log\log p}$ uniformly for all large primes $p \geq 2$. The expected order of magnitude is $g(p) \ll (\log p)^c$, $c > 1$ constant. The corresponding estimates for least prime primitive roots $g^*(p)$ are slightly higher. The last topics deal with Artin conjecture and its generalization to the set of integers. Some effective lower bounds such as $\#\{\, p \leq x : \mathrm{ord}(g) = p - 1\,\} \gg x\,(\log x)^{-1}$ for the number of primes $p \leq x$ with a fixed primitive root $g \neq \pm 1$, $b^2$ for all large number $x \geq 1$ will be provided. The current results in the literature have the lower bound $\#\{\, p \leq x : \mathrm{ord}(g) = p - 1\,\} \gg x\,(\log x)^{-2}$, and have restrictions on the minimal number of fixed integers to three or more.

# Table Of Contents

# 26. Class Numbers And Primitive Roots

# 27. Appendix

# 28. References

# Introduction

A few topics in the theory of primitive roots modulo primes $p \geq 2$, and primitive roots modulo integers $n \geq 2$, are studied in this monograph. The topics investigated are listed below.

## 1.1. Least Primitive Roots

Chapter 9 deals with estimates of the least primitive roots $g(p)$ modulo $p$, a large prime. One of the estimate here seems to sharpen the Burgess estimate

$$g(p) \ll p^{1/4+\epsilon} \tag{1.1}$$

for arbitrarily small number $\epsilon > 0$, to the smaller estimate

$$g(p) \ll p^{5/\log\log p} \tag{1.2}$$

uniformly for all large primes $p \geq 2$. The expected order of magnitude is $g(p) \ll (\log p)^c$, $c > 1$ constant. The result is stated in Theorem 9.1.

## 1.2. Least Prime Primitive Roots

Chapter 10 provides the details for the analysis of some estimates for the least prime primitive root $g^*(p)$ in the cyclic group $\mathbb{Z}/(p-1)\mathbb{Z}$, $p \geq 2$ prime. The current literature has several estimates of the least prime primitive root $g^*(p)$ modulo a prime $p \geq 2$ such as

$$g^*(p) \ll p^c, \ c > 2.8. \tag{1.3}$$

The actual constant $c > 2.8$ depends on various conditions such as the factorization of $p - 1$, et cetera. These results are based on sieve methods and the least primes in arithmetic progressions, see [MA98], [MB97], [HA13]. Moreover, there are a few other conditional estimates such as $g(p) \leq g^*(p) \ll (\log p)^6$, see [SP92], and the conjectured upper bound $g^*(p) \ll (\log p)(\log\log p)^2$, see [BH97]. On the other direction, there is the Turan lower bound $g^*(p) \geq g(p) = \Omega(\log p \log\log p)$, refer to [BE68], [RN96, p. 24], and [MT91] for discussions. The result stated in Theorem 10.1 improves the current estimate to the smaller estimate

$$g^*(p) \ll p^{5/\log\log p} \tag{1.4}$$

uniformly for all large primes $p \geq 2$.

## 1.3. Subsets of Primes with a Fixed Primitive Roots

The main topic in Chapter 12 deals with an effective lower bound

$$\#\{\, p \leq x : \operatorname{ord}(g) = p - 1 \,\} \gg x\,(\log x)^{-1} \tag{1.5}$$

for the number of primes $p \leq x$ with a fixed primitive root $g \neq \pm 1$, $b^2$ for all large number $x \geq 1$. The current results in the literature have the lower bound

$$\#\{\, p \leq x : \operatorname{ord}(g) = p - 1 \,\} \gg x\,(\log x)^{-2}, \tag{1.6}$$

and have restrictions on the minimal number of fixed integers to three or more, refer to Chapter 11, [HB86], and [MR88] for other details. The result is stated in Theorem 12.1.

## 1.4. Subsets of Composites with a Fixed Primitive Roots

The main topic in Chapter 13 deals with an effective lower bound

$$\#\{\, N \leq x : \operatorname{ord}(u) = \lambda(N) \,\} \gg x\,(\log\log\log x)^{-1} \tag{1.7}$$

for the number of composite $N \leq x$ with a fixed primitive root $u \neq \pm 1$, $v^2$, and $\gcd(u,\ N) = 1$, for all large number $x \geq 1$. The result is stated in Theorem 13.1. The current results in the literature have the heuristic asymptotic formula $\#\{\, N \leq x : \operatorname{ord}(u) = \lambda(N) \,\} \sim B\,x$, where $B = B(u) > 0$ is a constant. A result similar to the heuristic expectation is stated in Theorem 13.2.

## 1.5. Multiplicative Subsets of Composites with Fixed Primitive Roots

The topics in Chapter 23 uses Wirsing theorem to develop new exact asymptotic formulas such as

$$\#\{\, n \leq x : \operatorname{ord}(u) = \lambda(n) \,\} = c\,x\,(\log x)^{-\alpha} + o(x\,(\log x)^{-\alpha}), \tag{1.8}$$

where $0 < \alpha < 1$, and $c > 0$ is a constant, for the number of composite integers $n \leq x$ with a fixed primitive root $u \neq \pm 1$, $v^2$, and $\gcd(u,\ n) = 1$, for all large number $x \geq 1$. The results are stated in Theorems 23.1 and 23.2.

**Chapter 2**

# Characters And Exponential Sums

Let $G$ be a finite group of order $q = \# G$. The order $\mathrm{ord}(u)$ of an element $u \in G$ is the smallest integer $d \mid q$ such that $u^d = 1$. An element $u \in G$ is called a *primitive element* if it has order $\mathrm{ord}(u) = q$. A cyclic group $G$ is a group generated by a primitive element $\tau \in G$. Given a primitive root $\tau \in G$, every element $0 \neq u \in G$ in a cyclic group has a representation as $u = \tau^v$, $0 \leq v < q$. The integer $v = \log u$ is called the *discrete logarithm* of $u$ with respect to $\tau$.

## 2.1 Simple Characters Sums

A character $\chi$ modulo $q \geq 2$, is a complex-valued periodic function $\chi : \mathbb{N} \longrightarrow \mathbb{C}$, and it has order $\mathrm{ord}(\chi) = d \geq 1$ if and only if $\chi(n)^d = 1$ for all integers $n \in \mathbb{N}$, $\gcd(n, q) = 1$. For $q \neq 2^r$, $r \geq 2$, a multiplicative character $\chi$ of order $\mathrm{ord}(\chi) = d \mid q$ has a representation as

$$\chi(u) = e^{i\,2\,\pi\,k\,\log(u)/d}, \tag{2.1}$$

where $v = \log u$ is the discrete logarithm of $u \neq 0$ with respect to some primitive root, and for some integer $k \in \mathbb{Z}$, see [LN97, p. 187], [MV07, p. 118], [SG06, p. 15], and [IK04, p. 271]. The principal character $\chi_0 = 1 \bmod q$ has order $d = 1$, and it is defined by the relation

$$\chi_0(n) = \begin{cases} 1 & \text{if } \gcd(n, q) = 1, \\ 0 & \text{if } \gcd(n, q) \neq 1, \end{cases} \tag{2.2}$$

and the nonprincipal character $\chi \neq 1 \bmod q$ of order $\mathrm{ord}(\chi) = d \mid q$ is defined by the relation

$$\chi(n) = \begin{cases} \omega & \text{if } \gcd(n, q) = 1, \\ 0 & \text{if } \gcd(n, q) \neq 1, \end{cases} \tag{2.3}$$

where $\omega \in \mathbb{C}$ is a $d$th root of unity.

**Lemma** 2.1.   For a fixed integer $u \neq 0$, and an integer $q \in \mathbb{N}$, let $\chi \neq 1$ be nonprincipal character mod $q$, then

(i) $\displaystyle\sum_{\mathrm{ord}(\chi) = \varphi(q)} \chi(u) = \begin{cases} \varphi(q) & \text{if } u \equiv 1 \bmod q, \\ -1 & \text{if } u \not\equiv 1 \bmod q. \end{cases}$

$$\tag{2.4}$$

(ii) $\displaystyle\sum_{1 \le a < \varphi(q)} \chi(a\,u) = \begin{cases} \varphi(q) & \text{if } u \equiv 1 \bmod q, \\ -1 & \text{if } u \not\equiv 1 \bmod q. \end{cases}$

An additive character $\psi$ of order $\mathrm{ord}(\psi) = p$ has a representation as

$$\psi(u) = e^{i\,2\,\pi\,k\,u/p}, \tag{2.5}$$

for some $k \ge 1$, $\gcd(k, \ p) = 1$, see [LN97, p. 187], [IK04, p. 271]. It has character sums similar to Lemma 2.1.

**Lemma 2.2.** For a fixed integer $u$, and an integer $q \in \mathbb{N}$, let $\psi$ be an additive character of order $\mathrm{ord}\,\psi = q$, then

(i) $\displaystyle\sum_{\psi} \psi(u) = \begin{cases} q & \text{if } u \equiv 0 \bmod q, \\ 0 & \text{if } u \not\equiv 0 \bmod q. \end{cases}$

(ii) $\displaystyle\sum_{0 \le a < q} \psi(a\,u) = \begin{cases} q & \text{if } u \equiv 0 \bmod q, \\ 0 & \text{if } u \not\equiv 0 \bmod q. \end{cases}$

(2.6)

## 2.2 Estimates of Exponential Sums

Exponential sums indexed by the powers of elements of nontrivial orders have applications in mathematics and cryptography. These applications have propelled the development of these exponential sums. There are many results on exponential sums indexed by the powers of elements of nontrivial orders, the interested reader should consult the literature, and references within the cited paper.

**Theorem 2.3.** ([BN04, BJ07, Theorem 2.1]) (i) Given $\delta > 0$, there is $\epsilon > 0$ such that if $\theta \in \mathbb{F}_p$ is of multiplicative order $t \ge t_1 > p^{-\delta}$, then

$$\max_{1 \le a \le p-1} \left| \sum_{1 \le m \le t_1} e^{i 2\pi a\,\theta^m/p} \right| < t_1\,p^{-\epsilon}. \tag{2.7}$$

(ii) If $H \subset \mathbb{Z}_N$ is a subset of cardinality $\# H \ge N^\delta$, $\delta > 0$, then

$$\max_{1 \le a \le \varphi(N)} \left| \sum_{x \in H} e^{i 2\pi a\,x/N} \right| < N^{1-\delta}. \tag{2.8}$$

Other estimates for exponential sum over arbitrary subsets $H \subset \mathbb{F}_p$ are also given in the [KS12]. For the finite rings $\mathbb{Z}/N\mathbb{Z}$ of the integer modulo $N \ge 1$, similar results have been proved [BJ07].

**Theorem 2.4.** ([FS02, Lemma 4]) For integers $a$, $k$, $N \in \mathbb{N}$, assume that $\gcd(a, N) = c$, and that $\gcd(k, t) = d$.

(i) If the element $\theta \in \mathbb{Z}_N$ is of multiplicative order $t \geq t_0$, then

$$\max_{1 \leq a \leq p-1} \left| \sum_{1 \leq x \leq t} e^{i2\pi \, a \, \theta^x / N} \right| < c \, d^{1/2} \, N^{1/2} . \tag{2.9}$$

(ii) If $H \subset \mathbb{Z} / N \mathbb{Z}$ is a subset of cardinality $\# H \geq N^\delta$, $\delta > 0$, then

$$\max_{\gcd(a, \varphi(N))=1} \left| \sum_{x \in H} e^{i2\pi \, a \, \theta^x / N} \right| < N^{1-\delta} . \tag{2.10}$$

A multiplicative version using double finite sum, which offers greater flexibility in some applications, satisfies the following bound.

**Lemma** 2.5. ([CC09]) Let $p \geq 2$ be a large prime, and let $x \geq 1$ be a large real number. Then the followings estimates hold:

(i) If the element $\tau \in \mathbb{F}_p$ is of multiplicative order $p - 1$, and $c \in \mathbb{Z}$, $\gcd(c, p) = 1$ constant, then

$$\sum_{1 \leq n \leq x} e^{i2\pi \, c \, \tau^n / p} = O\left(p^{1/2} \log p\right) . \tag{2.11}$$

(i) If the elements $\tau_1, \tau_2, \ldots, \tau_k \in \mathbb{F}_p$ are of multiplicative orders $p - 1$, $k = \varphi(p-1)$, and $c_i \in \mathbb{Z}$, $\gcd(c_i, p) = 1$ constants, then for any small number $\epsilon > 0$,

$$\frac{1}{\varphi(p-1)} \sum_{1 \leq n \leq x} e^{i2\pi \, (c_0 n + c_1 \tau_1^n + \cdots + c_k \tau_k^n)/p} = O\left(p^{23/24+\epsilon}\right) . \tag{2.12}$$

**Lemma** 2.6. ([SK11, Lemma 4]) For any nonprincipal multiplicative character $\chi \bmod p$ and arbitrary sets $\mathcal{A}, \mathcal{B} \subseteq \{0, 1, 2, \ldots, p-1\}$, the following double finite sum estimate holds:

$$\left| \sum_{a \in \mathcal{A}} \sum_{b \in \mathcal{B}} \chi\left(a+b\right) \right| \ll (\#\mathcal{A})^{1-1/2v} (\#\mathcal{B}) \, p^{1/4v} + (\#\mathcal{A})^{1-1/2v} (\#\mathcal{B})^{1/2} \, p^{1/2v} . \tag{2.13}$$

A generalization of many of these estimates is compactly expressed in the well known estimate summarized below.

**Lemma** 2.7. (Weil) Let $p \geq 2$ be a large prime, let $q \mid p - 1$, and let $\chi \bmod q$ be a nonprincipal multiplicative character. Let $f(x) = (x - x_1)^{d_1} (x - x_2)^{d_2} \cdots (x - x_r)^{d_r}$ be a polynomial of degree $\deg(f) = d_1 + d_2 + \cdots + d_r = d$, $\gcd(d_i, q) = 1$, and $x_i \neq x_j$ for $i \neq j$. Then

$$\left| \sum_{1 \le x \le p} \chi(f(x)) \right| \le (d - 1) \, p^{1/2} \, . \tag{2.14}$$

## 2.3 Some Exact Evaluations

There are a handful of closed form evaluations of characters and exponentials sums. The best known closed form evaluation is the Gauss quadratic character sum. For an integers triple $0 \ne a, b, p \in \mathbb{Z}$, define the finite sum

$$S(a, b, p) = \sum_{0 \le x \le |p| - 1} e^{i \, 2 \, \pi \, (a \, x^2 + b \, x) / p} \, . \tag{2.15}$$

**Lemma 2.8.** (Gauss) If $n \ge 2$ is an integer, and $0 \ne a \in \mathbb{Z} \, / \, n \, \mathbb{Z}$, and $\left( \frac{a}{n} \right)$ denotes the quadratic symbol, then

$$\sum_{0 \le x \le n-1} e^{i \, 2 \, \pi \, a \, x^2 / n} = \begin{cases} \left( \frac{a}{n} \right) \sqrt{n} & \text{if } n \equiv 1 \bmod 4, \\ i \left( \frac{a}{n} \right) \sqrt{n} & \text{if } n \equiv 3 \bmod 4, \\ 0 & \text{if } n \equiv 2 \bmod 4, \\ (1 + i) \left( \frac{a}{n} \right) \sqrt{n} & \text{if } n \equiv 0 \bmod 4. \end{cases} \tag{2.16}$$

There is an amazing variety of proofs of the quadratic Gauss sum, these techniques vary from the elementary to very advanced theory, some are given in [BW98, p. 12].

**Lemma 2.9.** (Reciprocity theorem for Gaussian sum) Let $a, b, c \in \mathbb{Z}$, $a \, c \ne 0$, be integers. Then

$$S(a, b, c) = \left| \frac{c}{a} \right|^{1/2} e^{i \, \pi \, (\mathrm{sign}(a \, c) - b^2 / a \, c) / 4} \, S(-c, -b, a) \, , \tag{2.17}$$

where $\mathrm{sgn} : \mathbb{Z} \longrightarrow \{-1, \, 0, \, 1\}$ is the sign function.

The special case

$$\frac{1}{\sqrt{p}} \sum_{0 \le x \le p-1} e^{i \, 2 \, \pi \, q \, x^2 / p} = \frac{e^{i \, \pi / 4}}{\sqrt{2 \, q}} \sum_{0 \le x \le q-1} e^{i \, 2 \, \pi \, p \, x^2 / q} \tag{2.18}$$

is known as Schaar formula.

**Lemma 2.10.** Let $f(x) = a\,x^2 + b\,x + c \in \mathbb{Z}[x]$ be a quadratic polynomial, and let $p \geq 3$ be a prime then

$$\sum_{0 \leq x \leq p-1} \left( \frac{f(x)}{p} \right) = \begin{cases} -\left( \frac{a}{p} \right) & \text{if } b^2 - 4\,a\,c \not\equiv 0 \bmod p, \\[2mm] (p-1)\left( \frac{a}{p} \right) & \text{if } b^2 - 4\,a\,c \equiv 0 \bmod p. \end{cases} \tag{2.19}$$

Two distinct proofs are realized in [BW98, p. 58] and [LN97, p. 218].

**Chapter 3**

# Representations of the Characteristic Function

The characteristic function $\Psi : G \longrightarrow \{0, 1\}$ of primitive element is one of the standard tools employed to investigate the various properties of primitive roots in cyclic groups $G$. Many equivalent representations of characteristic function $\Psi$ of primitive elements are possible. Several types of representations of the characteristic function of primitive elements in finite rings are investigated here.

## 3.1 Divisors Dependent Characteristic Functions

Three different forms of the characteristic function based on the structure of the divisors of the order $q = \# G$ of the finite group $G$ are given below. These representations are sensitive to the prime decompositions of the orders of the cyclic groups.

Let $n = p_1^{v_1} p_2^{v_2} \cdots p_t^{v_t}$, $p_i \geq 2$ prime, and $v_i \geq 1$, be the prime decomposition off the integer $n \geq 2$. The Mobius function occurs in various formulae. It is defined by

$$\mu(n) = \begin{cases} (-1)^t & \text{if } n = p_1 p_2 \cdots p_t, \ v_i = 1 \text{ all } i \geq 1, \\ 0 & \text{if } n \neq p_1 p_2 \cdots p_t, \ v_i \neq 1 \text{ some } i \geq 1. \end{cases} \tag{3.1}$$

**Lemma 3.1.** Let $G$ be a finite group of order $q = \# G$, and let $0 \neq u \in G$ be an invertible element of the group. Assume that $v = \log u$, and $e = \gcd(d, v)$. Then

(i) $\Psi(u) = \dfrac{\varphi(q)}{q} \displaystyle\sum_{d \,|\, q} \dfrac{\mu(d)}{\varphi(d)} \sum_{\text{ord}(\chi) = d} \chi(u) = \begin{cases} 1 & \text{if ord } u = q, \\ 0 & \text{if ord } u \neq q. \end{cases}$

(ii) $\Psi(u) = \dfrac{\varphi(q)}{q} \displaystyle\sum_{d \,|\, q} \mu(d) \dfrac{\mu(d \,/\, e)}{\varphi(d \,/\, e)} = \begin{cases} 1 & \text{if ord } u = q, \\ 0 & \text{if ord } u \neq q. \end{cases}$

$\qquad\qquad\qquad\qquad\qquad\qquad\qquad\qquad\qquad\qquad\qquad\qquad\qquad\qquad\qquad\qquad\qquad\qquad\qquad (3.2)$

(iii) $\Psi(u) = \displaystyle\prod_{p \,|\, q} \left( \dfrac{(p-1)^2}{(p-1)^2 - 1} \right) \sum_{d \,|\, q} \dfrac{\mu(d)}{\varphi(d)^2} \left| \sum_{\text{ord}(\chi) = d} \chi(u) \right|^2 = \begin{cases} 1 & \text{if ord } u = q, \\ 0 & \text{if ord } u \neq q. \end{cases}$

**Proof**: (i) A multiplicative character $\chi$ of order $\text{ord}(\chi) = d \,|\, q$ has the form $\chi(u) = \chi_k(u) = e^{i \, 2\pi k \log(u)/d}$, where $v = \log u$ is the discrete logarithm of $u \neq 0$, and $\gcd(k, d) = 1$, see (2.1). This immediately gives a closed form evaluation of the characters sum

$$\sum_{\text{ord}(\chi) = d} \chi(u) = \begin{cases} \varphi(d), & \text{if } v \equiv 0 \bmod d, \\ -1 & \text{if } v \not\equiv 0 \bmod d, \end{cases} \tag{3.3}$$

where $\chi \neq 1$. Replacing this information into the product

$$\frac{\varphi(q)}{q} \sum_{d \mid q} \frac{\mu(d)}{\varphi(d)} \sum_{\text{ord}(\chi) = d} \chi(u) = \frac{\varphi(q)}{q} \prod_{p \mid q} \left( 1 - \frac{\sum_{\text{ord}(\chi) = p} \chi(u)}{p - 1} \right) \tag{3.4}$$

shows that $\Psi(u) = 0$, it vanishes, if and only if the element $u \in G$ has order $\text{ord}(u) = d \mid q$, and $d < q$. Exempli gratia, if this occurs, then the order of the element $u \neq 0$ is divisible by prime $p \mid d < q$. This in turns means that the characters sum $\sum_{\text{ord}(\chi) = p} \chi(u) = \varphi(p)$, so the product vanishes.

(ii) Write the characters sum as

$$\sum_{\text{ord}(\chi) = d} \chi(u) = \sum_{\gcd(k,d) = 1} e^{j\, 2\,\pi\, k\, \log u/d} = \mu(d\,/\,e)\, \frac{\varphi(d)}{\varphi(d\,/\,e)}, \tag{3.5}$$

where $\gcd(d, k) = 1$, the integer $v = \log u$ is the discrete logarithm of $u$ with respect to $g$, and $e = \gcd(d, v)$. For the transformation from exponential to arithmetic functions, see [AP86, p. 160], [MV07, p. 110]. Replacing this information into the product

$$\frac{\varphi(q)}{q} \sum_{d \mid q} \frac{\mu(d)}{\varphi(d)} \sum_{\text{ord}(\chi) = d} \chi(u) = \frac{\varphi(q)}{q} \sum_{d \mid q} \frac{\mu(d)}{\varphi(d)}\, \mu(d\,/\,e)\, \frac{\varphi(d)}{\varphi(d\,/\,e)}. \tag{3.6}$$

Thus, if $\text{ord}\, u = q$, then $e = \gcd(d, v) = 1$ since a primitive root has the representation $u = g^v$ with $\gcd(v, q) = 1$. This immediately leads to

$$\frac{\varphi(q)}{q} \sum_{d \mid q} \frac{\mu^2(d)}{\varphi(d)} = 1. \tag{3.7}$$

This completes the proof. For (iii) use the same method as in (i) mutatis mutandis. ∎

Finer details on the characteristic function are given in [ES57, p. 863], [JN00], [LN97, p. 258], [MO04, p. 18], et alii. The characteristic function for multiple primitive roots is used in [CB98, p. 146] to study consecutive primitive roots. In [KS12] it is used to study the gap between primitive roots with respect to the Hamming metric. In [SK11] it is used to study Fermat quotients as primitive roots. And in [WR01] it is used to prove the existence of primitive roots in certain small subsets $A \subset \mathbb{F}_{p^n}$, $n \geq 1$. Other related decompositions of the characteristic function are given in [MA98], and [AM14].

**Extension to Arbitrary Finite Fields**

The trace $\mathrm{Tr} : \mathbb{F}_{p^m} \longrightarrow \mathbb{F}_p$ is defined by the relation

$$\mathrm{Tr}(\alpha) = \alpha + \alpha^p + \alpha^{p^2} + \cdots + \alpha^{p^{m-1}}. \tag{3.8}$$

This is a $p^{m-1}-$ to $-1$ map. Similarly, the norm $N : \mathbb{F}_{p^m} \longrightarrow \mathbb{F}_p$ is defined by the relation

$$N(\alpha) = \alpha \cdot \alpha^p \cdot \alpha^{p^2} \cdots \alpha^{p^{m-1}} = \alpha^{(p^m-1)/(p-1)}. \tag{3.9}$$

This is a $p-1-$ to $-1$ map, [LN97, p .54].

Let $q = p^m$, $m \geq 1$. The extension of the characteristic function to an arbitrary finite field $\mathbb{F}_q$ is simply the composition

$$\Psi(u) = \frac{\varphi(q-1)}{q-1} \sum_{d \,|\, q-1} \frac{\mu(d)}{\varphi(d)} \sum_{\mathrm{ord}(\chi)=d} \chi(\mathrm{Tr}(u)) = \begin{cases} 1 & \text{if ord } u = q-1, \\ 0 & \text{if ord } u \neq q-1. \end{cases} \tag{3.10}$$

## 3.2 Divisors Free Characteristic Functions

It is often difficult to derive any meaningful result using the usual divisors dependent characteristic function of primitive elements, compare Lemma 3.1. This difficulty is due to the large number of terms that can be generated by the divisors of the order of the cyclic group involved in the calculations, see [ES57], [KS12] for typical applications and [MP04, p. 19] for a discussion. For example, the integer $n = 2 \cdot 3 \cdot 5 \cdots p$, and similar highly composite integers have unmanageable expansions

$$\sum_{d \,|\, n} \frac{\mu(d)}{\varphi(d)} \sum_{\mathrm{ord}(\chi)=d} \chi(u) = \frac{\mu(1)}{\varphi(1)} \sum_{\mathrm{ord}(\chi)=1} \chi(u) + \frac{\mu(2)}{\varphi(2)} \sum_{\mathrm{ord}(\chi)=2} \chi(u) + \frac{\mu(3)}{\varphi(3)} \sum_{\mathrm{ord}(\chi)=3} \chi(u) +$$

$$\cdots + \frac{\mu(n)}{\varphi(n)} \sum_{\mathrm{ord}(\chi)=d} \chi(u) \,. \tag{3.11}$$

A new *divisors-free* representation of the characteristic function of primitive element is developed here. This representation can overcomes some of the limitations of its counterpart in certain applications. The *divisors representation* of the characteristic function of primitive roots, Lemma 3.1, detects the order ord $u \,|\, p-1$ of the element $u \in \mathbb{F}_p$ by means of the divisors of the Totient $p-1$. In contrast, the *divisors-free representation* of the characteristic function, Lemma 3.2, detects the order ord $u \,|\, p-1$ of the element $u \in \mathbb{F}_p$ by means of the solutions of the equation $\tau^n - u = 0$ in $\mathbb{F}_p$, where $u$, $\tau$ are constants, and $1 \leq n < p-1$, $\gcd(n, p-1) = 1$, is a variable. Two versions are given: a multiplicative version, and an additive version.

**Lemma 3.2.** Let $p \geq 2$ be a prime, and let $\tau$ be a primitive root mod $p$. For a nonzero element $u \in \mathbb{F}_p$, the followings hold:

(i) If $\chi \neq 1$ be a nonprincipal multiplicative character of order ord $\chi = p$, then

$$\Psi(u) = \frac{1}{p} \sum_{\gcd(n, p-1) = 1,} \sum_{0 \leq k \leq p-1} \chi\left((\tau^n \, \overline{u})^k\right) = \begin{cases} 1 & \text{if ord } u = p-1, \\ 0 & \text{if ord } u \neq p-1, \end{cases} \tag{3.12}$$

where $\overline{u}$ is the inverse of $u$ mod $p$.

(ii) If $\psi \neq 1$ be a nonprincipal additive character of order ord $\psi = p$, then

$$\Psi(u) = \frac{1}{p} \sum_{\gcd(n, p-1) = 1,} \sum_{0 \leq k \leq p-1} \psi((\tau^n - u)\, k) = \begin{cases} 1 & \text{if ord } u = p-1, \\ 0 & \text{if ord } u \neq p-1. \end{cases} \tag{3.13}$$

**Proof**: (i) As the index $n \geq 1$ ranges over the integers relatively prime to $p - 1$, the element $\tau^n \in \mathbb{F}_p$ ranges over the primitive roots mod $p$. Ergo, the equation $\tau^n \, \overline{u} = 1$ , where $\overline{u}$ is the inverse of $u$, has a solution if and only if the fixed element $u \in \mathbb{F}_p$ is a primitive root. Next apply Lemma 2.1. For (ii), use the equation $\tau^n - u = 0$, and apply Lemma 2.4. ∎

The normalization used in this characteristic function was chosen to match the associated exponential sum, which has a known nontrivial upper bound. This normalization leads to the main term

$$\frac{\varphi(p-1)}{p}. \tag{3.14}$$

The different normalizations used in the construction of the characteristic functions have different main terms. But there is a simple partial summation conversion to transform it from one normalization to another normalization, see the last step of the proof of Theorem 12.1. An extended concept of the normalized parameter is considered in Section 3.3.

**Extension to Arbitrary Finite Fields**

Let $q = p^m$, $m \geq 1$. The extension of the divisors-free characteristic function to an arbitrary finite field $\mathbb{F}_q$ has the form

$$\Psi(u) = \frac{1}{q} \sum_{\gcd(n, q-1) = 1,} \sum_{0 \leq k \leq q-1} \psi(k \operatorname{Tr}(\tau^n - u)) = \begin{cases} 1 & \text{if ord } u = q-1, \\ 0 & \text{if ord } u \neq q-1. \end{cases} \tag{3.15}$$

The trace function is defined in equation (3.8).

## 3.3 Moduli Free Characteristic Function for Finite Rings

Another construction of the characteristic function assembled below, has a normalization $\frac{\varphi(\lambda(N))}{p}$ that is nearly independent of the cardinality of the finite group. This facilitates the selections of the simpler exponential sums estimates of primes moduli.

**Lemma 3.3.** Let $N \geq 1$ be an integer, and let $p = N + o(N)$ be a prime number. Let $G \subset \mathbb{Z}_N^*$ be a maximal cyclic subgroup of order $\lambda(N) = \# G$ generated by a primitive root $\theta \in G$. If $z \in G$, and $\gcd(z, N) = 1$, then the followings holds.

(i) If $\chi \neq 1$ be a nonprincipal multiplicative character of order $\operatorname{ord} \chi = \varphi(N)$, then

$$\Psi(u) = \sum_{\gcd(n,\lambda(N))=1} \left( \frac{1}{p-1} \sum_{0 \leq k < p-1} \chi\big((\theta^n\,\overline{u})^k\big) \right) = \begin{cases} 1 & \text{if } \operatorname{ord} u = \lambda(N), \\ 0 & \text{if } \operatorname{ord} u \neq \lambda(N), \end{cases} \tag{3.16}$$

where $\overline{u}$ is the inverse of $u \bmod N$.

(ii) If $\psi \neq 1$ be a nonprincipal additive character of order $\operatorname{ord} \psi = p$, then

$$\Psi(u) = \sum_{\gcd(n,\lambda(N))=1} \left( \frac{1}{p} \sum_{0 \leq k \leq p-1} \psi\big((\theta^n - u)\,k\big) \right) = \begin{cases} 1 & \text{if } \operatorname{ord} u = \lambda(N), \\ 0 & \text{if } \operatorname{ord} u \neq \lambda(N). \end{cases} \tag{3.17}$$

**Proof**: Replaced the prime $p \sim \lambda(N)$ in the appropriate places, and proceed in the same way as in Lemma 3.2. ∎

## 3.4 Characteristic Functions for Finite Rings

There are various possibility for constructing an extension of the divisors-free characteristic functions of an arbitrary finite ring $\mathbb{Z}_N$. The simplest formulas matching the parameters of the multiplicative group $\mathbb{Z}_N^*$, see Section 3.3, seems to be the following.

**Lemma 3.4.** Let $N \geq 1$ be an integer, let $G \in \mathbb{Z}_N$ be a maximal subgroup, and let $\theta$ be a primitive root $\bmod N$. If $z \in G$, and $\gcd(z, N) = 1$, then

(i) $\Psi(z) = \dfrac{1}{\lambda(N)} \displaystyle\sum_{\gcd(n,\lambda(N))=1,} \sum_{0 \leq k < \lambda(N)} \psi\big((\theta^n - z)\,k\big) = \begin{cases} 1 & \text{if } \operatorname{ord} z = \lambda(N), \\ 0 & \text{if } \operatorname{ord} z \neq \lambda(N). \end{cases}$

(ii) $\Psi(z) = \dfrac{1}{\varphi(N)} \displaystyle\sum_{\gcd(n,\varphi(N))=1,} \sum_{0 \leq k < \varphi(N)} \psi\big((\theta^n - z)\,k\big) = \begin{cases} 1 & \text{if } \operatorname{ord} z = \lambda(N), \\ 0 & \text{if } \operatorname{ord} z \neq \lambda(N). \end{cases}$

$$\tag{3.18}$$

The normalization of formula (i) leads to the main term

$$\frac{\varphi(\lambda(N))}{\lambda(N)}. \tag{3.19}$$

Likewise, the equivalent normalization in formula (ii) leads to the main term $\varphi(\varphi(N)) / \varphi(N)$. This later normalization is used in [MA98, Lemma 4]. A related problem uses $\lambda(\lambda(N))$, see [FR01]. The normalization used in the construction of the characteristic function can have an effect on the main term depending on the finite rings. But there is a simple partial summation conversion to transform it from one normalization to another normalization.

## 3.5 Quasi Characteristic Functions

The quasi characteristic functions of primitive elements have the same detection properties characteristic functions but have different, nearly equal, values. These functions are useful in the investigation of primitive roots.

**Lemma** 3.5.    Let $p \geq 2$ be a prime, let $\tau$ be a primitive root mod $p$, and let $u \in \mathbb{F}_p$ be a nonzero element. If $\psi(u) = e^{j 2 \pi k u / p}$ is a nonprincipal additive character of order ord $\psi = p$, then

$$\Psi_*(u) = \frac{1}{p} \sum_{\gcd(n, p-1)=1, \, 0 \leq k \leq p-1} \sum \psi\big((\tau^n - u)\, k^2\big) = \begin{cases} 1 - \frac{1}{p} + \frac{\beta}{\sqrt{p}} \sum_{\gcd(n,p-1)=1} \left( \frac{\tau^n - u}{p} \right) & \text{if ord } u = p-1, \\[2ex] \frac{\beta}{\sqrt{p}} \sum_{\gcd(n,p-1)=1} \left( \frac{\tau^n - u}{p} \right) & \text{if ord } u \neq p-1. \end{cases}$$
$$\tag{3.20}$$

where $\beta = 1$ if $p \equiv 1 \bmod 4$ or $\beta = i$ if $p \equiv 3 \bmod 4$, and

$$\left( \frac{z}{p} \right) = \begin{cases} 1 & \text{if } z \equiv a^2 \bmod p, \\ -1 & \text{if } z \not\equiv a^2 \bmod p. \end{cases} \tag{3.21}$$

**Proof**: Assume that the multiplicative order ord$(u) \neq p - 1$. Then, $\tau^n - u \neq 0$ for all $n \in \mathbb{N}$, $\gcd(n, p - 1) = 1$. Ergo, by Lemma 2.8, the evaluation reduces to

$$\frac{1}{p} \sum_{\gcd(n, p-1)=1, \, 0 \leq k \leq p-1} \sum \psi\big((\tau^n - u)\, k^2\big) = \frac{1}{p} \sum_{\gcd(n, p-1)=1} \left( \frac{\tau^n - u}{p} \right) \sum_{0 \leq x \leq p-1} \left( \frac{x}{p} \right) e^{j 2 \pi x / p}$$
$$= \frac{\beta}{\sqrt{p}} \sum_{\gcd(n, p-1)=1} \left( \frac{\tau^n - u}{p} \right), \tag{3.22}$$

where $\beta = 1$ if $p \equiv 1 \bmod 4$ or $\beta = i$ if $p \equiv 3 \bmod 4$, and $\left(\frac{a}{p}\right)$ is the quadratic symbol, this follows from Lemma 2.8. The simplification

$$\sum_{0 \le x \le p-1} e^{i\,2\pi\,a\,x^2/p} = \sum_{0 \le x \le p-1} \left(1 + \left(\frac{x}{p}\right)\right) e^{i\,2\pi\,a\,x/p} = \left(\frac{a}{p}\right) \sum_{0 \le x \le p-1} \left(\frac{x}{p}\right) e^{i\,2\pi\,x/p} \tag{3.23}$$

was used to reduce the finite sum. The second part is similar, but uses $\tau^n - u = 0$ for all $n \in \mathbb{N}$, $\gcd(n,\,p-1) = 1$. $\blacksquare$

### 3.6 Basic Statistics of the Finite Rings $\mathbb{Z}_N$.

Let $N \in \mathbb{N}$ be an integer. The multiplicative group $\mathbb{Z}_N^* = \{\, z \ne 0 : \gcd(z,N) = 1 \,\}$ of the finite ring $\mathbb{Z}_N$ has $\varphi(N)$ elements. The exponent $\exp(\mathbb{Z}_N) \ge 1$ is the cardinality $\# G$ of a maximal cyclic subgroup $G \subset \mathbb{Z}_N^*$. The Carmichael function $\lambda(N) = \exp(\mathbb{Z}_N)$ gives the precise value. An element $\theta \in \mathbb{Z}_N^*$ of maximal order $\operatorname{ord}\theta = \lambda(N)$ is called a primitive element. Since $\theta^k$ has maximal order $\operatorname{ord}\theta^k = \lambda(N)$ if and only if $\gcd(k, \lambda(N)) = 1$, a maximal subgroup $G$ has precisely

$$\varphi(\lambda(N)) = \# \{\, k : \gcd(k, \lambda(N) = 1 \,\} \tag{3.24}$$

primitive elements, see Chapters 6 and 7 for further information on the arithmetic functions $\varphi$ and $\lambda$.

**Lemma 3.6.** Let $N \ge 1$ be an integer, and let $\mathbb{Z}_N$ be a finite ring. Let $G_1, G_2, \ldots, G_m$ be distinct maximal subgroups $G_1, G_2, \ldots, G_m$ of order $\lambda(N) = \# G_i$, where the index $m = \varphi(N)/\lambda(N)$. If $G_1 \bigcap G_2 \bigcap \cdots \bigcap G_m = \{\, 1 \,\}$, then, the number of primitive roots in the multiplicative group $\mathbb{Z}_N^*$ is given by

$$\varphi(\varphi(N)) = \varphi(N) \prod_{p\,|\,\varphi(N)} \left(1 - \frac{1}{p}\right). \tag{3.25}$$

**Proof**: The the multiplicative group $\mathbb{Z}_N^*$ has $\varphi(N)$ units, and a cyclic subgroup $G \subset \mathbb{Z}_N$ has a maximal order of $\lambda(N) = \# G$. Thus, the index $m = \varphi(N)/\lambda(N)$ is the number of distinct maximal cyclic groups of order $\lambda(N) = \# G_i$. The maximal cyclic subgroups constitute a partition of the multiplicative group $\mathbb{Z}_N^* = G_1 \bigcup G_2 \bigcup \cdots \bigcup G_m$, and each subgroup has precisely $\varphi(\lambda(N))$ primitive roots. Therefore, there is a total of

$$\varphi(\lambda(N)) \, \frac{\varphi(N)}{\lambda(N)} = \varphi(N) \, \frac{\varphi(\lambda(N))}{\lambda(N)} = \varphi(N) \prod_{p\,|\,\lambda(N)} \left(1 - \frac{1}{p}\right) \tag{3.26}$$

primitive roots. Lastly, observe that $\lambda(N)\,|\,\varphi(N)$. $\blacksquare$

For the prime factorization $N = p_1^{v_1} p_2^{v_2} \cdots p_t^{v_t}$, with $p_i$ prime and $v_i \ge 1$, the irreducible decomposition of the multiplicative group has the form

$$\mathbb{Z}_N^* \cong U\left(p_1^{v_1}\right) \times U\left(p_2^{v_2}\right) \times \cdots \times U\left(p_t^{v_t}\right). \tag{3.27}$$

Each factor $U\left(p_i^{v_i}\right)$ is a cyclic group if $p_i^{v_i}$ is an odd prime power, confer Lemma 22.1 and see [ML04], and [CP14].

**Lemma 3.7.** Let $N \geq 1$ be an integer, and let $\mathbb{Z}_N$ be a finite ring. Then, the number of primitive roots in the multiplicative group $\mathbb{Z}_N^*$ is given by

$$\varphi(N) \prod_{p \,|\, \varphi(N)} \left(1 - \frac{1}{p^{m(p)}}\right), \tag{3.28}$$

where $m(p)$ is the number of elementary divisors of the irreducible decomposition of the multiplicative group whose $p$-part is maximal.

## 3.7 Summatory Functions

**Lemma 3.8.** Let $p \geq 2$ be a prime, let $G$ be the finite cyclic group of order $p - 1 = \# G$, and let $x \geq 1$ be a large number. Then

(i) $\displaystyle\sum_{n \leq x} \Psi(n) = \frac{\varphi(p-1)}{p-1}\, x + o(x),$

(ii) $\displaystyle\sum_{n \leq k\,(p-1)} \Psi(n) = k\, \varphi(p-1), \ \ k \geq 1.$

$$\tag{3.29}$$

**Lemma 3.9.** Let $N \geq 2$ be an integer, let $\mathbb{Z}_N$ be a finite ring, and let $x \geq 1$ be a large number. Then

(i) $\displaystyle\sum_{n \leq x} \Psi(n) = \frac{\varphi(\varphi(N))}{\varphi(N)}\, x + o(x),$

(ii) $\displaystyle\sum_{n \leq k\,\varphi(N)} \Psi(n) = k\, \varphi(\varphi(N)), \ \ k \geq 1.$

$$\tag{3.30}$$

More specific information on the arithmetic functions $\varphi$ and $\lambda$ are given in Chapters 5 to 7.

# Prime Divisors Counting Functions

For $n \in \mathbb{N}$, the prime counting function is defined by $\omega(n) = \# \{\, p \mid n \,\}$. Somewhat similar proofs of the various properties of the arithmetic function $\omega(n)$ are given in [HW08, p. 473], [MV07, p. 55], [CR06, p. 34], [TM95, p. 83], and other references.

## 4.1 Average Orders

The individual values of the prime divisors counting function have an irregular pattern, but the average value is stable.

**Lemma** 4.1.　　Let $x \geq 1$ be a large number. Then

$$\sum_{n \leq x} \omega(n) = x \log\log x + B_1 x + (\gamma - 1)\, x / \log x + O\!\left(x\, e^{-c \sqrt{\log x}}\right),$$

(4.1)

where $\gamma$ and $B_1$ are Euler and Mertens constants, and $c > 0$ is an absolute constant.

**Proof**:　Let $\{ x \} = x - [\, x \,]$ be the fractional part function, and expand the finite sum:

$$
\begin{aligned}
\sum_{n \leq x} \omega(n) &= \sum_{p \leq x} \left[ \frac{x}{p} \right] \\
&= \sum_{p \leq x} \left( \frac{x}{p} - \left\{ \frac{x}{p} \right\} \right) \\
&= x \left( \log\log x + B_1 + O\!\left( e^{-c \sqrt{\log x}} \right) \right) - \sum_{p \leq x} \left\{ \frac{x}{p} \right\}.
\end{aligned}
$$

(4.2)

Lastly, apply Lemmas 17.1 and 17.9 to complete the calculations of the estimate. ■

The Euler constant is defined by the asymptotic formula $\gamma = \lim_{x \to \infty} (\sum_{n \leq x} 1 / n - \log x)$, and the Mertens constant is defined by the asymptotic formula $B_1 = \lim_{x \to \infty} \left( \sum_{p \leq x} 1 / p - \log\log x \right)$. There are additional details in Chapter 18.

The best conditional error term, assuming the RH, has the form

$$\sum_{n \leq x} \omega(n) = x \log\log x + B_1 x + (\gamma - 1) x / \log x + O\left(x^{1/2} \log x\right).$$

(4.3)

And the oscillations properties of the sequence of prime numbers induces the omega result

$$\sum_{n \leq x} \omega(n) = x \log\log x + B_1 x + (\gamma - 1) x / \log x + \Omega_{\pm}\left(x^{1/2}\right).$$

(4.4)

**Lemma 4.2.** Let $n \geq 1$ be a large integer. Then

(i) Almost every integer $n \geq 1$ satisfies the inequality

$$\omega(n) \leq \log\log n + B_1 + (\gamma - 1) / \log n + O\left(e^{-c\sqrt{\log n}}\right).$$

(4.5)

(ii) Every integer $n \geq 1$ satisfies the inequality

$$\omega(n) \leq \log n / \log\log n + O\left(\log n \, e^{-c\sqrt{\log n}}\right).$$

(4.6)

**Proof**: (i) This follows from Lemma 4.1. For (ii), it is sufficient to let $n = \prod_{p \leq x} p$ be the product of the primes up $p \leq x$. Taking logarithm produces

$$\log n = \sum_{p \leq x} \log(p) = x + O\left(x \, e^{-c\sqrt{\log x}}\right).$$

(4.7)

The right side follows from the prime number theorem, [MV07, p. 179], [TM95, p. 166]. Ergo, $p \leq x \leq \log n$. Next, observe that there are $\pi(x) \leq \pi(\log n) \leq \log n / \log\log n + O\left(\log n \, e^{-c\sqrt{\log n}}\right)$ primes in the very short interval $[1, x] \subset [1, \log n]$, and $\omega(n) \leq \pi(x) \leq \pi(\log n)$. ∎

**Lemma 4.3.** Let $n \geq 1$ be a large integer. Then

(i) Almost every integer $n \geq 1$ satisfies the average number of squarefree divisors inequality

$$2^{\omega(n)} \leq 2^{\log\log n + B_1 + (\gamma - 1)/\log n + o(1)}.$$

(4.8)

(ii) For every integer $n \geq 1$, the number of squarefree divisors satisfies the inequality

$$2^{\omega(n)} \leq 2^{\log n/\log\log n \, + o(\log n/\log\log n)} \, . \tag{4.9}$$

The proof of Lemma 4.3 follows from Lemma 4.2. Alternatively, this can be proved utilizing the expressions for squarefree integers

$$2^{\omega(n)} = \sum_{d \,|\, n} \mu^2(d) \leq \sum_{d \,\leq\, x} \mu^2(d) = \frac{6}{\pi^2}\, x + O\!\left(x^{1/2}\right), \tag{4.10}$$

for some $x \geq 1$. The actual value of $x$ is implied by the maximal order $\log n \, / \log\log n$ of the arithmetic function $\omega(n) \geq 1$, see [TM, p. 83]. Specifically, $x = n^{\log 2/\log\log n \, -\left(6\pi^{-2}\right)/\log n \, + o(1/\log\log n)}$.

## 4.2 Distribution of Prime Divisors

The average order, the variance $\sum_{n \,\leq\, x} (\omega(n) - \log\log n)^2 = x \log\log x + o(x \log\log x)$, and other information lead to a refinement of these results, known as the Erdos-Kac theorem. This refinement shows that the average number of prime factors of a large integer $n \geq 1$ has a normal distribution of mean $\mu = \log\log n$ and standard deviation $\sigma = \sqrt{\log\log n}$, confer [CR06, p. 34] for other details.

## 4.3 Problems And Exercises

**Problem 4.1.** Let $x \in \mathbb{R}$ be a large real number, and $\omega : \mathbb{N} \longrightarrow \mathbb{N}$ be the prime divisors counting function. Prove

$$\sum_{n \,\leq\, x} \frac{1}{\omega(n)} = \frac{x}{\log\log x} + O\!\left(\frac{x}{(\log\log x)^2}\right). \tag{4.11}$$

**Problem 4.2.** Let $x \in \mathbb{R}$ be a large real number, and $\omega, \, \Omega : \mathbb{N} \longrightarrow \mathbb{N}$ be the prime divisors counting functions. Prove

$$\sum_{n \,\leq\, x} \frac{\Omega(n)}{\omega(n)} = x + O\!\left(\frac{x}{\log\log x}\right). \tag{4.12}$$

**Problem 4.3.** Let $x \in \mathbb{R}$ be a large real number, and $\omega : \mathbb{N} \longrightarrow \mathbb{N}$ be the prime divisors counting function. Prove

$$\sum_{n \,\leq\, x} \Omega(n) = x \log\log x + C_1 \, x + O\!\left(\frac{x}{\log x}\right), \text{ where the constant is } C_1 = B_1 + \sum_{p \,\geq\, 2} \frac{1}{p(p-1)} = 1.034653 \ldots . \tag{4.13}$$

**Problem 4.4.** Let $x \in \mathbb{R}$ be a large real number, and $\omega : \mathbb{N} \longrightarrow \mathbb{N}$ be the prime divisors counting function. Prove

$$\sum_{p \,\leq\, x} \omega(p-1) = x \log\log x + C_2 \, x + O\!\left(\frac{x}{\log x}\right), \text{ where } C_2 > 0 \text{ is } a \text{ constant.} \tag{4.14}$$

**Problem 4.5.** Let $x \in \mathbb{R}$ be a large real number, and $\omega : \mathbb{N} \longrightarrow \mathbb{N}$ be the prime divisors counting function. For $a < q < x$, prove

$$\sum_{n \leq x, \, n \equiv a \bmod q} \omega(n) = x \log\log x + C_3 \, x + O\left(\frac{x}{\log x}\right), \quad \text{where } C_3 > 0 \text{ is } a \text{ constant.}$$

(4.15)

**Problem 4.6.** Let $x \in \mathbb{R}$ be a large real number, and $\omega : \mathbb{N} \longrightarrow \mathbb{N}$ be the prime divisors counting function. For $a = O(\log x)$, approximate the finite sum

$$\sum_{n \leq x} \omega(n) \, \omega(n + a) = ? \; .$$

(4.16)

**Chapter 5**

# Euler Totient Function

The Euler totient function counts the number of relatively prime integers $\varphi(n) = \# \{ k : \gcd(k, n) = 1 \}$. This is compactly expressed by the analytic formulae given below.

**Lemma 5.1.** If $n \in \mathbb{N}$ is an integer, then

(i) $\varphi(n) = n \prod_{p \mid n} (1 - 1 / p)$.

(ii) $\varphi(n) = n \sum_{d \mid n} \dfrac{\mu(d)}{d}$. $\hspace{3cm}$ (5.1)

**Lemma 5.2.** (Fermat-Euler) If $a \in \mathbb{Z}$ is an integer such that $\gcd(a, n) = 1$, then $a^{\varphi(n)} \equiv 1 \bmod n$.

## 5.1 Estimates and Average Orders of the Euler Totient Function

The Euler totient function has a trivial upper bound $\varphi(n) \leq n / 2$. The derivation of some lower bounds required some works. Various inequalities for the totient function are computed in [RS62].

**Lemma 5.3.** For any integer $n \in \mathbb{N}$ the Euler totient function satisfies the followings estimates

(i) $\varphi(n) \gg n / \log\log\log n$, $\hspace{2cm}$ for almost every integer $n \in \mathbb{N}$.

(ii) $\varphi(n) \gg n / \log\log n$, $\hspace{2cm}$ for every integer $n \in \mathbb{N}$. $\hspace{2cm}$ (5.2)

(iii) $\varphi(n) \geq n / (e^{\gamma} \log\log n + 2.5 / \log\log n)$, $\hspace{1cm}$ for every integer $n \in \mathbb{N}$.

**Lemma 5.4.** ([CL12]) For large number $x \geq 1$, the average order of the Euler totient function are as follows:

(i) $\sum_{n \leq x} \varphi(n) = \dfrac{3}{\pi^2} x^2 + O(x)$,

(ii) $\sum_{n \leq x} \dfrac{\varphi(n)}{n} = \dfrac{6}{\pi^2} x + O(1)$, $\hspace{3cm}$ (5.3)

(iii) $\displaystyle\sum_{n \le x} \varphi(n)\,\{x\,/\,n\} = \left(\frac{6}{\pi^2} - \frac{1}{2}\right) x^2 + O(x).$

**Proof**: (iii) Let $\{z\} = z - [z]$ be the fractional part function. Expanding the Sylvester formula yields

$$\frac{x(x+1)}{2} = \sum_{n \le x} \varphi(n)\left[\frac{x}{n}\right] = x \sum_{n \le x} \frac{\varphi(n)}{n} - \sum_{n \le x} \varphi(n)\,\{x\,/\,n\}. \tag{5.4}$$

Substituting (ii) completes the verification. $\blacksquare$

**Lemma 5.5.** For large number $x \ge 1$, the average order of the Euler totient function are as follows:

(i) $\displaystyle\sum_{n \le x} \frac{1}{\varphi(n)} = c_0 \log x + c_1 + O\!\left(\frac{\log x}{x}\right),$

(ii) $\displaystyle\sum_{n \le x} \frac{1}{n^a\,\varphi(n)^b} = c_{ab}\, x^{1-a-b} + o\!\left(x^{1-a-b}\right),$ $\qquad(5.5)$

(iii) $\displaystyle\sum_{n \le x} \frac{n}{\varphi(n)} = b_0\, x + b_1 \log x + o(\log x),$

where $b_0,\ b_1,\ c_0,\ c_1,$ and $c_{ab}$ are constants.

## 5.2 Average Gap

The nth gap $d_n = \varphi(n+1) - \varphi(n)$ between the values of the Euler totient function $\varphi(n)$ has extreme values of unknown frequencies. For example, at the primes $n + 1 = 1 + \prod_{p \le z} p$, it has the extreme values

$$\varphi(n+1) - \varphi(n) = \left(1 + O\!\left(\frac{1}{\log z}\right)\right)\prod_{p \le z} p \sim \prod_{p \le z} p, \tag{5.6}$$

where $z \ge 2$ is a number, refer to [CW95] for some information on the sequence of primes $1 + \prod_{p \le z} p$. However, on average, the gaps between the values of the Euler totient function are small:

$$\frac{1}{x}\sum_{n \le x}(\varphi(n+1) - \varphi(n)) = (\log n)\,\Big/\,e^{c_0(\log_3(x) - \log_4(x))^2 + c_1 \log_3(x) + c_2 \log_4(x) + c_3}, \tag{5.7}$$

see Theorem 5.5. These information imply that the infinite sum $\sum_{n \ge 1} 1\,/\,\varphi(n) = \infty$ diverges. These ideas are the crux of the next results.

## 5.3 Density of the Euler Totient Function

By definition, a subset $R \subset \mathbb{Q}$ of rational numbers is called dense in $\mathbb{R}$ if for each fixed number $\alpha \in \mathbb{R}$, and any small number $\epsilon > 0$, the inequality $| \alpha - r | < \epsilon$ has infinitely many rational solutions $r \in R$.

**Theorem** 5.6. (i) The set of points $\{ \varphi(n) / n : n \in \mathbb{N} \}$ is dense in the unit interval $(0, 1)$.

(ii) The set of points $\{ \varphi(n) : n \in \mathbb{N} \}$ is dense in the interval $(0, \infty)$.

**Proof**: (i) Given $\epsilon > 0$, and a real number $x_0 \in (0, 1)$, let $n_0 \geq \liminf \{ n \in \mathbb{N} : \varphi(n) / n \geq x_0 \}$. For example, choose a prime such that $n_0 = p > 1 / (1 - x_0)$. Then, it is clear that there is an infinite sequence of integers, i.e., $\{ n = n_0 \, m : \text{some } m \in \mathbb{N} \}$, such that

$$\left| x_0 - \frac{\varphi(n)}{n} \right| < \epsilon \, . \tag{5.8}$$

This implies the dense property. For (ii), given any small number $\epsilon > 0$, and any real number $w_0 \in (0, \infty)$, let $n_0 \geq \liminf \left\{ n \in \mathbb{N} : w_0 \leq -\log \prod_{p \,|\, n} (1 - 1 / p) \right\}$. Observe that

$$f(n) = -\log \frac{\varphi(n)}{n} = -\log \prod_{p \,|\, n} (1 - 1 / p) = -\sum_{p \,|\, n} \log(1 - 1 / p) \tag{5.9}$$

is a finite sum of continuous functions, so it is a continuous function $f : \mathbb{N} \longrightarrow (0, \infty)$ with respect to some discrete metric or valuation $\nu : \mathbb{N} \longrightarrow \mathbb{R}$. The range of $f$ follows from $0 < \varphi(n) / n < 1$, and the fact that the series

$$\sum_{p \,\leq\, x} \log(1 - 1 / p) = C_1 + \sum_{p \,\leq\, x} \frac{1}{p} \geq \log\log x, \tag{5.10}$$

where $C_1 > 0$ is a constant, diverges to $\infty$ as $z \longrightarrow \infty$.

Now, given any small number $\epsilon > 0$, and any real number $w_0 \in (0, \infty)$, there exists a small number $\delta > 0$ such that $| f(n) - w_0 | < \epsilon$ for all integers $n \in D(n_0) = \{ n \in \mathbb{N} : \nu(n - n_0) < \delta \}$, which is a small open disk in some discrete topology with respect to the valuation $\nu$. This proves (ii). ∎

**Example** 5.1. The real number $1 / \sqrt{2} = 0.707106781186547524400844 \ldots$ has the sequence of totient numbers approximations, but not unique: Take a small prime $n_0 = p \geq 5 > 1 / \left(1 - 1 / \sqrt{2}\right)$. The sequence of rational appoximations is

$$\frac{\varphi(5)}{5} = \frac{4}{5}, \quad \frac{\varphi(5 \cdot 11)}{5 \cdot 11} = \frac{8}{11}, \quad \frac{\varphi(5 \cdot 11 \cdot 37)}{5 \cdot 11 \cdot 37} = \frac{288}{407}, \quad \frac{\varphi(5 \cdot 11 \cdot 37 \cdot 1409)}{5 \cdot 11 \cdot 37 \cdot 1409} = \frac{36\,864}{52\,133}, \quad \ldots . \tag{5.11}$$

Thus, given any small real number $\epsilon > 0$, the inequality $\left| 1 / \sqrt{2} - \frac{\varphi(n)}{n} \right| < \epsilon$ is satisfied by infinitely many integers $n \in \mathbb{N}$. In fact, much more can be stated, for example, $\left| 1 / \sqrt{2} - \frac{\varphi(n)}{n} \right| < 2 / n^2$.

## 5.4 Range of Values of the Totient Function

The subset of integers $\mathcal{V} = \{ \varphi(n) : n \in \mathbb{N} \} \subset \mathbb{N}$ represents the image of the Euler totient function. The natural density defined by

$$\delta(\mathcal{V}) = \lim_{x \longrightarrow \infty} \frac{\# \{ \varphi(n) \leq x : n \in \mathbb{N} \}}{x} \tag{5.12}$$

is a well defined quantity. The cardinality of the range of the Euler totient function over the interval $[1, x]$ is defined by $V(x) = \# \{ \varphi(n) \leq x : n \in \mathbb{N} \}$. The totient numbers theorem, which is analogous to the prime numbers theorem, takes the form given here.

**Theorem 5.7.** ([FD11, Theorem 1]) For all large numbers $x \geq 1$, the cardinality of the range of the Euler totient function has the asymptotic formula

$$V(x) = \frac{x}{\log x} e^{c_0(\log_3(x) - \log_4(x))^2 + c_1 \log_3(x) + c_2 \log_4(x) + c_3}, \tag{5.13}$$

where $c_0 = 0.81781464\ldots$, $c_1 = 2.17696874\ldots$, $c_2 < 0$, and $c_3 \in \mathbb{R}$ are constants.

For $k \geq 1$, the notation $\log_k x = \log(\cdots \log(x) \cdots)$ denotes the $k$-fold iterated logarithm.

## 5.5 Average Order of the Number of Primitive Roots in the Multiplicative Group

As computed in Lemma 3.5, the number of primitive roots modulo some odd integers $n \geq 1$ is exactly $\varphi(\varphi(n))$. In these cases, the average order of the number of primitive roots in the multiplicative group $\mathbb{Z}_n^*$ is as follows.

**Lemma 5.8.** For large number $x \geq 1$. Then, there is an absolute constant $v_0 > 0$ such that

$$\sum_{\varphi(n) \leq x} \frac{\varphi(\varphi(n))}{\varphi(n)} = v_0 \frac{x}{\log x} e^{C(\log_3(x) - \log_4(x))^2 + O(\log_3(x))} + O(\log x). \tag{5.14}$$

**Proof**: Use the identity $\frac{\varphi(n)}{n} = \sum_{d \mid n} \frac{\mu(d)}{d}$ to expand the finite sum:

$$\sum_{\varphi(n) \leq x} \frac{\varphi(\varphi(n))}{\varphi(n)} = \sum_{n \in [1, x] \cap \mathcal{V}} \frac{\varphi(n)}{n} = \sum_{n \in [1, x] \cap \mathcal{V}} \sum_{d \mid n} \frac{\mu(d)}{d} = \sum_{d \in [1, x] \cap \mathcal{V}} \frac{\mu(d)}{d} \sum_{n \in [1, x] \cap \mathcal{V}_d} 1, \tag{5.15}$$

where $\mathcal{V}_d = \{\, \varphi(n) \equiv 0 \bmod d : n \in \mathbb{N} \,\}$. By Theorem 5.7, this reduces to

$$
\begin{aligned}
\sum_{\varphi(n) \leq x} \frac{\varphi(\varphi(n))}{\varphi(n)} &= \sum_{d \in [1,\,x] \cap \mathcal{V}} \frac{\mu(d)}{d} \left( \frac{V(x)}{d} - \{V(x)/d\} \right) \\
&= V(x) \sum_{d \in [1,\,x] \cap \mathcal{V}} \frac{\mu(d)}{d^2} - \sum_{d \in [1,\,x] \cap \mathcal{V}} \frac{\mu(d)}{d} \{V(x)/d\},
\end{aligned}
\tag{5.16}
$$

where $V(x) = \#\{\, \varphi(n) \leq x : n \in \mathbb{N} \,\}$, and $\{z\} = z - [z]$ is the fractional part function. Since the infinite sum $v_0 = \sum_{n \in \mathcal{V}} \mu(n)\, n^{-2}$ converges, the claim is immediate. $\blacksquare$

Often, the maximal subgroup and the multiplicative group coincides: $G \subset \mathbb{Z}_N^*$. This occurs whenever the index $\varphi(n)/\lambda(n) = 1$. For example, $N = 2,\ 4,\ p^m,\ 2\,p^m$, with $p > 2$ an odd prime and $m \geq 1$, confer Lemma 22.1. Some numerical data is available in [CP14]. In light of this idea, define the subset of integers $\mathcal{R} = \{\, n \in \mathbb{N} : \varphi(n)/\lambda(n) = 1 \,\}$, and its cardinality $R(x) = \#\{\, n \leq x : \varphi(n)/\lambda(n) = 1 \,\}$.

**Theorem 5.9.** ([ML04, Theorem 2]) The ratios of arithmetic functions

$$
\frac{\varphi(\lambda(n))}{\lambda(n)} = \frac{\varphi(\varphi(n))}{\varphi(n)}
\tag{5.17}
$$

is satisfied on a subset of integers of zero density. In particular, $R(x) \gg x / \log^B x$, where $B > 0$ is a constant.

## 5.6 Multiplicity of the Totient Function

Fix an integer $m \in \mathbb{N}$, and consider the multiplicity $\#\{\, k \in \mathbb{N} : \varphi(k) = n \,\}$ of the totient function $\varphi : \mathbb{N} \longrightarrow \mathbb{N}$. As the totient function increases $\varphi(n) \longrightarrow \infty$ with increasing argument $n \longrightarrow \infty$, it is trivial that the multiplicity of the equation $\varphi(n) = m$ is finite.

**Theorem 5.10.** ([BT72]) Let $x \geq 1$ be a large number, then

$$
\#\{\, n \in \mathbb{N} : \varphi(n) \leq x \,\} = \frac{\zeta(2)\,\zeta(3)}{\zeta(6)} x + O\!\left(e^{-c(\log x)^{1/2}}\right),
\tag{5.18}
$$

where $c > 0$ is a constant.

**Proof**: Let $f(n) = \#\{\, k \in \mathbb{N} : \varphi(k) = n \,\}$ be the multiplicity of $n \geq 1$. As done in [BT72, p. 331], proceed to apply the Wiener-Ikehara theorem to the analytic function

$$F(s) = \sum_{k \geq 1} \frac{1}{\varphi(k)^s} = \sum_{n \geq 1} \frac{f(n)}{n^s} = \prod_{p \geq 2} \left( 1 + \frac{1}{\varphi(p)^s} + \frac{1}{\varphi(p^2)^s} + \frac{1}{\varphi(p^3)^s} + \cdots \right) = \zeta(s) \prod_{p \geq 2} \left( 1 - \frac{1}{p^s} + \frac{1}{(p-1)^s} \right) \tag{5.19}$$

to verify the claim. ∎

The amalgamation of the last two results produces a very precise estimate of the mean multiplicity.

***Corollary* 5.11.** The average multiplicity of the totient function on the interval $[1, x]$ is

$$\overline{V_k}(x) = \frac{\zeta(2)\,\zeta(3)}{\zeta(6)} \left( (\log x) \,\Big/\, e^{-c_0(\log_3(x) - \log_4(x))^2 - c_1 \log_3(x) - c_2 \log_4(x) - c_3} \right) + o\left( \frac{\log x}{e^{-c(\log_3(x) - \log_4(x))^2}} \right), \tag{5.20}$$

where $c > 0$ is an absolute constant.

## 5.7 Squarefree Values of the Totient Function

The Euler totient function has a nonsquarefree value, that is, $\mu(\varphi(n)) = 0$, for almost every integer $n \in \mathbb{N}$. Further, the subset of integers $\{\, n \in \mathbb{N} : \mu(\varphi(n)) \neq 0 \,\}$ for which it is squarefree has a very precise description.

***Lemma* 5.12.** Let $x \geq 1$ be a large number, and let $B > 1$ be a constant. Then

$$\sum_{n \leq x} \mu^2(\varphi(n)) = \frac{3}{2} \, \text{li}(x) \prod_{p \geq 2} \left( 1 - \frac{1}{p(p-1)} \right) + O\left( \frac{x}{\log^B x} \right). \tag{5.21}$$

***Proof***: Let $p^\alpha \,\|\, n$, $\alpha \geq 0$, be the maximal prime power divisor of the integer $n \in \mathbb{N}$. The product formula

$$\varphi(n) = \prod_{p^\alpha \| n} \varphi(p^\alpha) = \prod_{p^\alpha \| n} p^{\alpha - 1}(p - 1) \tag{5.22}$$

shows that $\varphi(n)$ is squarefree if and only if $\varphi(p) = p - 1$ is squarefree for each $p^\alpha \,\|\, n$, and $\alpha \leq 2$. For example, at the integers $n = 2,\ 4,\ p,\ 2\,p,\ p^2,$ and $2\,p^2$, with $p > 2$ an odd prime. Therefore, there is a total of

$$\sum_{n \leq x} \mu^2(\varphi(n)) = \#\,\{\, p \leq x : p - 1 \text{ is squarefree} \,\} + \#\,\{\, p \leq x/2 : p - 1 \text{ is squarefree} \,\} + \cdots$$
$$= \pi_{\text{sf}}(x) + \pi_{\text{sf}}(x/2) + \pi_{\text{sf}}\big(x^{1/2}\big) + \cdots, \tag{5.23}$$

where $\pi_{\text{sf}}(x) = \#\,\{\, p \leq x : \mu(p - 1) \neq 0 \,\}$. The exact formula follows from Lemma 6.3. ∎

The set of integers $\mathbb{N} = \{0,\ 1,\ 2,\ 3,\ \dots\}$ contains a subset of integers relatively primes to the image of the totient function. The proof of the corresponding asymptotic formula is surprisingly involved. This is covered in [MV07, p. 397] and a few other places in the literature.

***Theorem*** **5.13.** (Erdos) Let $x \geq 1$ be a large number, and let $N(x) = \#\ \{\ n \leq x : \gcd(n,\ \varphi(n)) = 1$. Then

$$N(x) = \frac{e^{-\gamma}\ x}{\log\log\log x}\left(1 + O\!\left(\frac{\log\log\log\log x}{\log\log\log x}\right)\right). \tag{5.24}$$

## 5.6 Average Order Over Arithmetic Progressions

***Theorem*** **5.14.** Let $x \geq 1$ be a large number, and let $a < q = O(\log x)$. Then

$$\sum_{n \leq x,\ n \equiv a \bmod q} \varphi(n) = \frac{3}{\pi^2}\ \frac{q}{\varphi_2(q)}\prod_{p\,|\,c}\left(1 - 1\,/\,p^2\right)x^2 + O(x\log x), \tag{5.25}$$

where $\varphi_s(n) = q^s \prod_{p\,|\,q}\left(1 - 1\,/\,p^s\right)$, and $c = \gcd(a,\ q)$.

The complete proof of this result appears in [PV86, Chapter 4].

## 5.9 Problems And Exercises

**Problem 5.1.** Let $k \geq 1$ be a fixed integer, and let $x \in \mathbb{R}$ be a large real number. Find the kth moment of the totient function:

$$\sum_{n \leq x} n^k\ \varphi(n) = c_k\ x^{k+2} + O\!\left(x^{k+1}\right),\ \text{ where } c_k > 0 \text{ is } a \text{ constant.} \tag{5.26}$$

**Problem 5.2.** Show that the nth gap of the totient function has the asymptotic

$$\varphi(n + 1) - \varphi(n) \sim n. \tag{5.27}$$

**Problem 5.3.** Let $x \in \mathbb{R}$ be a large real number, and $k = O(\log x)$. Approximatwe the finite sum

$$\sum_{n \leq x} \varphi(n)\ \varphi(n + k) = ?\ . \tag{5.28}$$

**Problem 5.4.** Prove that the average of the integer primitive root gap modulo $n \geq 2$ is approximated by

$$\sum_{n \leq x} \frac{n}{\varphi(n)} = \frac{\zeta(2)\ \zeta(3)}{\zeta(6)}\ x + O(\log x). \ \text{ Use the identity }\ \frac{n}{\varphi(n)} = \sum_{d\,|\,n} \frac{\mu^2(d)}{\varphi(d)}. \tag{5.29}$$

**Problem 5.5.** Let $\mathcal{V} = \{\, \varphi(n) : n \in \mathbb{N} \,\}$ be a subset of integers. Compute the constant

$$v_0 = \sum_{n \in \mathcal{V}} \mu(n)\, n^{-2}. \tag{5.30}$$

**Problem 5.6.** Prove an asymptotic formula for the average order of squarefree totients:

$$\sum_{n \leq x} \mu(n)^2\, \varphi(n) = ? \tag{5.31}$$

**Problem 5.7.** Prove that, use the identity $\frac{n}{\varphi(n)} = \sum_{d \mid n} \frac{\mu^2(d)}{\varphi(d)}$,

$$\sum_{n \geq 1} \frac{1}{\varphi(n)^s} = \zeta(s) \prod_{p \geq 2} \left( 1 - \frac{1}{p^s} + \frac{1}{(p-1)^s} \right) = \zeta(s)\, \frac{\zeta(2s)\, \zeta(3s)}{\zeta(6s)}. \tag{5.32}$$

**Problem 5.8.** Compute an asymptotic formula for the average order of the iterated totient:

$$\sum_{n \leq x} \varphi(\varphi(n)) = c\, x^2\,?, \ \ c > 0 \text{ constant}. \tag{5.33}$$

**Chapter 6**

# Finite Sums And Products Of Totients

Several results for finite sums of totients and products over the primes numbers are investigated here. Some of these estimates can be used to obtain exact solutions or nearly exact solutions of the asymptotic formula for the number of primes with a fixed primitive root, other are useful for obtaining estimates.

## 6.1 Finite Products over the Primes

Various estimates for the product of the consecutive primes are investigated in [RS62]. A useful case is stated here.

***Lemma* 6.1.**    Let $x \geq 1$ be a large number, then the prime product satisfies

$$\prod_{p \leq x} \left( 1 - \frac{1}{p} \right) = \frac{e^{-\gamma}}{\log x} \left( 1 + O\left( \frac{1}{\log x} \right) \right). \tag{6.1}$$

***Lemma* 6.2.**    Let $p \geq 3$ be a large prime, then the totient ratio satisfies

(i)    $\dfrac{\varphi(p-1)}{p-1} \leq \dfrac{1}{2},$                                    for any prime $p \geq 3$.

(ii)    $\dfrac{\varphi(p-1)}{p-1} \gg \dfrac{1}{\log\log p},$                        for any prime $p \geq 3$.

(iii)    $\dfrac{\varphi(p-1)}{p-1} \asymp \dfrac{1}{\log\log\log p},$                for almost any prime $p \geq 3$.

(iv)    $\dfrac{\varphi(p-1)}{p-1} \asymp \dfrac{1}{\log\log p},$                    on a subset of primes of zero density.

$$\tag{6.2}$$

***Proof***: (iii) Let $\epsilon, \delta > 0$ be small numbers, and let $\omega(n)$ be the number of prime divisors of an integer $n \in \mathbb{N}$. By Lemma 4.2, or Erdos-Kac theorem, the number of prime divisors function satisfies $\omega(n) \geq \log\log n - (\log\log n)^{1/2+\delta}$ for almost all integers $n \in \mathbb{N}$. Next observe that the interval $[1, \ a \log\log n]$ contains more than $a(\log\log n)^{1-\epsilon} < \omega(n)$ primes, see [HW08, Theorem 431], [MV07, Theorem 2.12]. Replace this in Lemma 5.1 to reach the inequalities

$$\frac{\varphi(p-1)}{p-1} = \prod_{r\,|\,p-1}\left(1-\frac{1}{r}\right) \geq \prod_{p\,\leq\, a\,\log\log p}\left(1-\frac{1}{p}\right) \gg \frac{1}{\log\log\log p}\,, \tag{6.3}$$

where $a \geq 1$ is a constant. Similarly, $\omega(n) \leq \log\log n + (\log\log n)^{1/2+\delta}$ for almost all integers $n \in \mathbb{N}$. And the interval $[1,\ b\log\log n]$ contains more than $b(\log\log n)^{1-\epsilon} < \omega(n)$ primes. Thus, for $n = p - 1$,

$$\prod_{p\,\leq\, x}\left(1-\frac{1}{p}\right) \ll \prod_{p\,\leq\, b\,\log\log p}\left(1-\frac{1}{p}\right)$$
$$\ll \frac{1}{\log\log\log p}\,. \tag{6.4}$$

The proof of (iv) is similar, and the proofs of (i) and (ii) use simpler techniques, ipso facto. ∎

**Lemma 6.3.** Let $x \geq 1$ be a large number, and let $\varphi(n)$ be the Euler totient function. Then

(i) $\displaystyle\sum_{p\,\leq\, x}\frac{\varphi(p-1)}{p-1} \leq \frac{x}{2}$,      for any large number $x \geq 3$.

(ii) $\displaystyle\sum_{p\,\leq\, x}\frac{\varphi(p-1)}{p-1} \gg \frac{x}{(\log x)^2}$,      for any large number $x \geq 3$.

(iii) $\displaystyle\sum_{p\,\leq\, x}\frac{\varphi(p-1)}{p-1} \gg \frac{x}{(\log x)\log\log x}$,      for almost any large number $x \geq 3$.

(iv) $\displaystyle\sum_{p\,\leq\, x}\frac{\varphi(p-1)}{p-1} \asymp \frac{x}{(\log x)\log\log\log x}$,      on a subset of numbers $x \geq 1$ of zero density.

$$\tag{6.5}$$

**Proof**: (iii) By Lemma 6.2, the inequality

$$\sum_{p\,\leq\, x}\frac{\varphi(p-1)}{p-1} \geq c_1\sum_{p\,\leq\, x}\frac{1}{\log\log p} = \int_2^x \frac{1}{\log\log t}\,d\,\pi(t)\,, \tag{6.6}$$

where $\pi(x) = x/\log x + O\!\left(x/\log^2 x\right)$ is the prime counting function, and $c_1 > 0$ is a constant. Evaluate the integral to complete the estimate. The verifications of the other statements are similar. ∎

**Lemma 6.4.** Let $x \geq 1$ be a large number, and let $\varphi(n)$ be the Euler totient function, and let $\omega(n)$ be the number of prime divisors of an integer $n \in \mathbb{N}$. Then

$$\sum_{p \leq x} \frac{\varphi(p-1)}{p-1} 2^{\omega(p-1)} \ll x \log x \,.$$

(6.7)

**Proof**: It is sufficient to consider a weak upper bound

$$\sum_{p \leq x} \frac{\varphi(p-1)}{p-1} 2^{\omega(p-1)} \leq \sum_{n \leq x} 2^{\omega(n)} \ll x \log x \,.$$

(6.8)

These complete the verifications. ∎

For information on the summatory function $\sum_{n \leq x} 2^{\omega(n)} = \left( 6 \big/ \pi^2 \right) x \log x + O(x)$, see [RM08, p. 328], [TM95, p. 53], et alii.

## 6.2 The Squarefree Totients

The subset $\{ p \in \mathbb{P} : p - 1 \text{ is squarefree} \}$ of primes $p \geq 2$ with squarefree totients $p - 1$ is the backbone of the theory of primitive root modulo $p$. The most basic statistic on this subset of primes is the density/counting function stated below. Let $\pi_{\mathrm{sf}}(x) = \# \{ p \leq x : p - 1 \text{ is squarefree} \}$ be the corresponding counting function.

**Lemma** 6.5.    For a large number $x \geq 1$, the number of primes $p \geq 2$ for which $p - 1$ is a squarefree integer has the asymptotic formula

$$\pi_{\mathrm{sf}}(x) = \mathrm{li}(x) \prod_{p \geq 2} \left( 1 - \frac{1}{p(p-1)} \right) + O\left( \frac{x}{\log^B x} \right),$$

(6.9)

where $\mathrm{li}(x) = \int_2^x (\log t)^{-1} \, d\,t$ is the logarithm integral, and $B > 1$ is a constant.

**Proof**: This require an evaluation of the finite sum

$$\sum_{p \leq x} \mu(p-1) = \sum_{p \leq x} \sum_{d^2 \mid p-1} \mu(d) \,.$$

(6.10)

Proceed to use standard analytic techniques to complete the verification. ∎

The constant $C_0 = \prod_{p \geq 2} \left( 1 - \frac{1}{p(p-1)} \right) = \sum_{n \geq 1} \mu(n) \left( n \, \varphi(n) \right)^{-1}$, which is the density of squarefree totients $p - 1$, is ubiquitous in the theory of primitive roots.

## 6.3 Moments of the Density Function

**Lemma 6.6.** ([ST69], [VN73]) Let $x \geq 1$ be a large number, and let $\varphi(n)$ be the Euler totient function. For $k \geq 1$, the $k$th moment

$$\sum_{p \leq x} \left( \frac{\varphi(p-1)}{p-1} \right)^k = \mathrm{li}(x) \prod_{p \geq 2} \left( 1 - \frac{\left( 1 - (1 - 1/p)^k \right)}{(p-1)} \right) + O\left( \frac{x}{\log^B x} \right), \tag{6.11}$$

where $\mathrm{li}(x)$ is the logarithm integral, and $B > 1$ is a constant, as $x \to \infty$.

**Proof**: To simplify the notation, assume that $k = 1$. Use the identity $\frac{\varphi(n)}{n} = \sum_{d \mid n} \frac{\mu(d)}{d}$ to expand the finite sum:

$$\begin{aligned}
\sum_{p \leq x} \frac{\varphi(p-1)}{p-1} &= \sum_{p \leq x, \ d \mid p-1} \frac{\mu(d)}{d} \\
&= \sum_{d \leq x} \frac{\mu(d)}{d} \sum_{p \leq x, \ p \equiv 1 \bmod d} 1 \\
&= \sum_{d \leq x} \frac{\mu(d)}{d} \pi(x, 1, d) \\
&= \sum_{d \leq \log^B x} \frac{\mu(d)}{d} \pi(x, 1, d) + \sum_{\log^B x < d \leq x} \frac{\mu(d)}{d} \pi(x, 1, d).
\end{aligned} \tag{6.12}$$

The dyadic summation in the fourth line assumes that $B > 2$ is a constant. Applying the Siegel-Walfisz theorem, see Theorem 16.3, to the first sum yields

$$\begin{aligned}
\sum_{d \leq \log^B x} \frac{\mu(d)}{d} \pi(x, 1, d) &= \mathrm{li}(x) \sum_{d \leq \log^B x} \frac{\mu(d)}{d \, \varphi(d)} + O\left( x \, e^{-c(\log x)^{1/2}} \sum_{\log^B x < d \leq x} \frac{1}{d \, \varphi(d)} \right) \\
&= \mathrm{li}(x) \sum_{d \geq 1} \frac{\mu(d)}{d \, \varphi(d)} + \mathrm{li}(x) \sum_{d > \log^B x} \frac{\mu(d)}{d \, \varphi(d)} + O\left( x \, e^{-c(\log x)^{1/2}} \right) \\
&= \mathrm{li}(x) \sum_{d \geq 1} \frac{\mu(d)}{d \, \varphi(d)} + O\left( x \, e^{-c(\log x)^{1/2}} \right).
\end{aligned} \tag{6.13}$$

An application of the Brun–Titchmarsh theorem, see Theorem 16.4, to the second sum yields

$$\sum_{\log^B x < d \le x} \frac{\mu(d)}{d} \pi(x, 1, d) \le \frac{x}{\log^B x} \sum_{d \ge 1} \frac{1}{d} = O\left(\frac{x}{\log^{B-1} x}\right).$$

(6.14)

Combining the two sums yields the result. ∎

## 6.4. Evaluation Of The Relatively Prime Counting Measure

An estimate of a finite sum over the numbers $n \le x$ relatively prime to $\varphi(q)$ is calculated here.

**Lemma 6.7.** Let $x \ge 1$ be a large number, and let $p \ge 3$ be a prime. Then

(i) $\displaystyle\sum_{\gcd(n, p-1)=1, \, n \le x} 1 \le \frac{x}{2}$ .

(ii) $\displaystyle\sum_{\gcd(n, p-1)=1, \, n \le x} 1 = \frac{\varphi(p-1)}{p-1} x + O(x^\epsilon)$ ,

(6.15)

where $\epsilon > 0$ is a small number.

**Proof**: The first statement (i) follows from Lemma 6.2. To verify the second statement (ii), use an equivalent finite sum, then evaluate it:

$$\sum_{\gcd(n, p-1)=1, \, n \le x} 1 = \sum_{n \le x, \, d \mid n, \, d \mid p-1} \mu(d) = \sum_{d \mid p-1} \mu(d) \sum_{n \le x, \, d \mid n,} 1 = \sum_{d \mid p-1} \mu(d) \left[\frac{x}{d}\right].$$

(6.16)

Substitute the largest integer function $[z] = z - \{z\}$ in the last equation to find

$$\sum_{\gcd(n, p-1)=1, \, n \le x} 1 = x \sum_{d \mid p-1} \frac{\mu(d)}{d} - \sum_{d \mid p-1} \mu(d) \{x / d\} = \frac{\varphi(p-1)}{p-1} x + O(x^\epsilon) ,$$

(6.17)

where $\sum_{d \mid n} 1 = O(n^\epsilon)$, $\epsilon > 0$ is a small number, see [HW08, Theorem 319]. ∎

Finer details on the error term are discussed in [MV07, p. 80].

## 6.4 Average Order Over Primes

**Theorem 6.8.** Let $x \ge 1$ be a large number. Then

$$\sum_{p \leq x} \varphi(p-1) = \frac{\alpha_2}{2} \frac{x^2}{\log^2 x} + o\left(\frac{x^2}{\log^2 x}\right), \tag{6.18}$$

where $\alpha_2 = \prod_{p \geq 2} (1 - 1 / p(p-1))$ is Artin constant.

A detailed and straight Forward proof appears in [GD09].

## 6.5 Problems And Exercises

**Problem 6.1.** Let $k \geq 1$ be a fixed integer, and let $x \in \mathbb{R}$ be a large real number. Find the $k$th moment of the totient function:

$$\frac{n}{\varphi(n)} \sum_{d \mid n} \frac{\mu(d)^2}{\varphi(d)}. \tag{6.19}$$

**Problem 6.2.** Compute the average order

$$\sum_{p \leq x} \varphi(p-1) = \frac{\alpha_2}{2} \frac{x^2}{\log x} + O\left(\frac{x^2}{\log^2 x}\right), \text{ where } \alpha_2 > 0 \text{ is Artin constant.} \tag{6.20}$$

**Problem 6.3.** Show that the average gap between primitive roots modulo $p \leq x$ has the asymptotic

$$\sum_{p \leq x} \frac{p-1}{\varphi(p-1)} = \prod_{p \geq 2}\left(1 - \frac{1}{(p-1)^2}\right) + O\left(\frac{\log\log x}{\log x}\right). \tag{6.21}$$

**Problem 6.4.** Let $k \geq 1$. Compute the kth moment

$$\sum_{p \leq x} \left(\frac{p-1}{\varphi(p-1)}\right)^k = ? \tag{6.22}$$

**Problem 6.5.** Let $a \geq 1$., see Theorem 6.8. Compute the shifted prime sum

$$\sum_{p \leq x} \varphi(p-a) = ? \tag{6.23}$$

**Problem 6.6.** Let $\alpha \in \mathbb{R} - \mathbb{Z}$ be an irrational number. Show that $\left| \alpha - \frac{\varphi(n)}{n} \right| < 2 / n^2$.

**Problem 6.7.** Compute the $k$th moment

$$\sum_{p \leq x} \left(\frac{\varphi(p-1)}{p-1}\right)^k = \text{li}(x) \prod_{p \geq 2} \left(1 - \left(1 - (1-1/p)^k\right) / (p-1)\right)^{-1} + O\left(\frac{x}{\log^B x}\right), \tag{6.24}$$

where $k \geq 1$ is an integer, and $B > 1$ is a constant, see [ST69] and [VN73].

**Problem 6.8.** Let $\epsilon > 0$ be a small number. Show $\varphi(p-1)/(p-1) > \epsilon$ and $\varphi(p-1)/(p-1) < \epsilon$ infinitely often as $p \longrightarrow \infty$.

**Problem 6.9.** Verify the totient function identity

$$\sum_{p \leq x} p^k \, \varphi(p-1) \overset{?}{=} b_k \, \frac{x^{k+2}}{\log x} + O\left(\frac{x^{k+1}}{\log x}\right), \quad \text{where } b_k > 0 \text{ is } a \text{ constant.} \tag{6.25}$$

　

<div style="background:gray">**Chapter 7**</div>

# Carmichael Function

The Carmichael function is basically an refinement of the Euler totient function to the finite ring $\mathbb{Z}_N$. Given an integer $N = p_1^{v_1} \, p_2^{v_2} \cdots p_t^{v_t}$, the Carmichael function is defined by

$$\lambda(N) = \text{lcm}\big(\lambda\big(p_1^{v_1}\big), \, \lambda\big(p_2^{v_2}\big) \cdots \lambda(p_t^{v_t})\big) = \prod_{p^v \,\|\, \lambda(N)} p^v, \tag{7.1}$$

where the symbol $p^v \,\|\, n$, $v \geq 0$, denotes the maximal prime power divisor of $n \geq 1$, and

$$\lambda(p^v) = \begin{cases} \varphi(p^v) & \text{if } p \geq 3 \text{ or } v \leq 2, \\ 2^{v-2} & \text{if } p = 2 \text{ and } v \geq 3. \end{cases} \tag{7.2}$$

The two functions coincide, that is, $\varphi(n) = \lambda(n)$ if $n = 2, 4, p^m$, or $2\,p^m$, $m \geq 1$. And $\varphi(2^m) = 2\,\lambda(2^m)$. In other cases, there are some simple relationships between $\varphi(n)$ and $\lambda(n)$. In fact, it seamlessly improves the Fermat-Euler Theorem: The improvement provides the least exponent $\lambda(n) \mid \varphi(n)$ such that $a^{\lambda(n)} \equiv 1 \bmod n$.

**Lemma** 7.1. (CR07, p. 233)   If $a, n \in \mathbb{N}$ are integers, such that $\gcd(a, n) = 1$, then $a^{\lambda(n)} \equiv 1 \bmod n$.

## 7.1 Estimates and Average Orders of the Carmichel Function

The Carmichael function has a trivial upper bound $\lambda(n) \leq n\,/\,2$. The derivation of some lower bounds required some works.

The value of the Carmichael function $\lambda(N) \geq 1$ measures the largest order of the elements of the cyclic group of units of a finite ring $\mathbb{Z}_N$. For each integer $n \in \mathbb{N}$, it has the lower bound

$$\lambda(n) > (\log n)^{c\, \log\log\log n}, \ \ c > 0 \text{ constant.} \tag{7.3}$$

**Theorem** 7.2.   ([EP91, Theorem 3])   For all $x \geq 16$, the average order of the Carmichael function has the asymptotic formula

$$\sum_{n \le x} \lambda(n) = \frac{x^2}{\log x} e^{B \frac{\log\log x}{\log\log\log x} (1+o(1))},$$ (7.4)

where the constant is defined by the product $B = e^{-\gamma} \prod_{p \ge 2} \left(1 - (p-1)^{-2}(p+1)^{-1}\right) = .34537 \ldots$.

***Corollary* 7.3.** For all $x \ge 16$, the average order of the normalized Carmichael function is

$$\sum_{n \le x} \frac{\lambda(n)}{n} = (1 + o(1)) \frac{x}{\log x} e^{B \frac{\log\log x}{\log\log\log x} (1+o(1))}.$$ (7.5)

***Proof***: Use partial summation, Theorem A-1 in Appendix A, to derive the following:

$$\sum_{n \le x} \frac{\lambda(n)}{n} = \sum_{n \le x} \frac{1}{n} \lambda(n) = \int_1^x \frac{1}{t} \, d \, R(t),$$ (7.6)

where $R(t) = \sum_{n \le x} \lambda(n)$. Continue to evaluate the integral as

$$\int_1^x \frac{1}{t} \, d \, R(t) = \frac{x}{\log x} e^{B \frac{\log\log x}{\log\log\log x} (1+o(1))} + \int_1^x \frac{1}{\log t} e^{B \frac{\log\log x}{\log\log\log x} (1+o(1))} \, d \, t.$$ (7.7)

***Lemma* 7.4.** Let $n \in \mathbb{N}$ be an integer, and let $x \ge 1$ be a large number. Then

(i) $\varphi(\lambda(n)) = \lambda(n) \prod_{p \mid \lambda(n)} \left(1 - p^{-1}\right) = \lambda(n) \prod_{p \mid \varphi(n)} \left(1 - p^{-1}\right),$

(ii) $\varphi(\lambda(n)) \ge \dfrac{1}{10} \dfrac{\lambda(n)}{\log\log\log n}$        for almost every integer $n \in \mathbb{N}$, (7.8)

(iii) $\varphi(\lambda(n)) \ge \dfrac{1}{e^{\gamma}} \dfrac{\lambda(n)}{\log\log n + 2.5 / \log\log n}$        for every integer $n \in \mathbb{N}$.

***Proof***: For (ii), apply Lemmas 4.2 and 5.2. The first statement (i) follows from $\lambda(N) \mid \varphi(N)$. The last inequality (iii) is similar to 5.3-iii. ∎

***Lemma* 7.5.** Let $x \ge 1$ be a large number. Then

(i) $\displaystyle\sum_{n \le x} \frac{\varphi(\lambda(n))}{\lambda(n)} \gg \frac{x}{\log\log x}$ .

(ii) $\displaystyle\sum_{n \le x} \frac{\varphi(\lambda(n))}{\lambda(n)} \le \frac{x}{2}$ .

**Proof**: (i) By Lemma 7.4, the average order is rewritten as

$$\sum_{n \le x} \frac{\varphi(\lambda(n))}{\lambda(n)} \gg \sum_{n \le x} \frac{1}{\log\log n} \gg \int_{20}^{x} \frac{1}{\log\log t}\, d\,t \gg \frac{x}{\log\log x} \ . \tag{7.10}$$

The first finite sum dominates the second finite sum, so it is absorbed by the first finite sum. Accordingly, a single integral yields an effective approximation. Lastly, for any integer $n \ge 2$, the product is bounded:

$$\frac{\varphi(\lambda(n))}{\lambda(n)} = \prod_{r \mid \lambda(n)} \left(1 - \frac{1}{p}\right) \le \frac{1}{2} \ . \tag{7.11}$$

This proves (ii). ∎

## 7.2 Range of Values of the Carmichael Function

The subset of integers $\mathcal{U} = \{\, \lambda(n) : n \in \mathbb{N} \,\}$ represents the image of the Carmichael function. And its counting function is defined by $U(x) = \#\{\, \lambda(n) \le x : n \in \mathbb{N} \,\}$.

**Theorem** 7.6. ([FE14, Theorem 1]) For all large numbers $x \ge 1$, the cardinality of the range of the Carmichael function has the asymptotic formula

$$\#\{\, \lambda(n) \le x : n \in \mathbb{N} \,\} = \frac{x}{(\log x)^{\eta + o(1)}}, \tag{7.12}$$

where $\eta = 1 - (1 + \log\log 2) / \log 2 = .08607\ldots$ is a constant.

The $\lambda(n)$-number theorem is the cumulative efforts of many years of works by scores of authors. The complexity of this analysis is comparable to or exceeds the complexity of the proof of the prime numbers theorem.

**Theorem** 7.7. For all large numbers Let $x \ge 1$ be a large number, and let $a$, $q \in \mathbb{N}$ be fixed integers, such that $\gcd(a, q) = 1$, and $q = O(\log^B x)$, $B \ge 0$. Then, the Carmichael function in an arithmetic progression has the asymptotic formula

$$U(x, q, a) = \frac{\varphi(q)}{q} \frac{x}{(\log x)^{\eta + o(1)}},$$ (7.13)

where $\eta = 1 - (1 + \log\log 2) / \log 2 = .08607\ldots$ is a constant.

**Lemma 7.8.** Let $x \geq 1$ be a large number, and let $n \in \mathbb{N}$ be an integer. Then, the followings holds.

(i) If $\varphi(n) \leq x$, then $n \leq e^\gamma x \log\log x$.

(ii) If $\lambda(n) \leq x$, then $n \leq e^\gamma x \log\log x$. (7.14)

The upper bound (i) is discussed in the short paper, [SA14], the second one should be similar.

## 7.3 Average Order of the Number of Primitive Roots in Subgroups of the Multiplicative Group

The number of primitive roots modulo $N \geq 1$ in a maximal subgroup $G \subset \mathbb{Z}_N^*$ is exactly $\varphi(\lambda(N))$. The average order of the number of primitive roots in a maximal multiplicative subgroup $\mathbb{Z}_N^*$ is as follows.

**Theorem 7.9.** For large number $x \geq 1$, there is an absolute constant $u_0 > 0$ such that

$$\sum_{\lambda(n) \leq x} \frac{\varphi(\lambda(n))}{\lambda(n)} = u_0 \frac{x}{(\log x)^{\eta + o(1)}} + O(\log x),$$ (7.15)

where $\eta = 1 - (1 + \log\log 2) / \log 2 = .08607\ldots$ is a constant.

**Proof**: Use the identity $\varphi(n) / n = \sum_{d \mid n} \mu(d) / d$ to expand the finite sum:

$$\sum_{\lambda(n) \leq x} \frac{\varphi(\lambda(n))}{\lambda(n)} = \sum_{n \in [1, x] \cap \mathcal{U}} \frac{\varphi(n)}{n} = \sum_{n \in [1, x] \cap \mathcal{U}, } \sum_{d \mid n} \frac{\mu(d)}{d} = \sum_{d \in [1, x] \cap \mathcal{U}} \frac{\mu(d)}{d} \sum_{n \in [1, x] \cap \mathcal{U}_d} 1,$$ (7.16)

where $\mathcal{U}_d = \{ \lambda(n) \equiv 0 \bmod d : n \in \mathbb{N} \}$. This reduces to

$$\sum_{\lambda(n) \leq x} \frac{\varphi(\lambda(n))}{\lambda(n)} = \sum_{d \in [1, x] \cap \mathcal{U}} \frac{\mu(d)}{d} \left( \frac{U(x)}{d} - \{U(x) / d\} \right)$$

$$= U(x) \sum_{d \in [1, x] \cap \mathcal{U}} \frac{\mu(d)}{d^2} - \sum_{d \in [1, x] \cap \mathcal{U}} \frac{\mu(d)}{d} \{U(x) / d\},$$ (7.17)

where $U(x) = \#\{\lambda(n) \leqslant x : n \in \mathbb{N}\}$, and $\{z\} = z - [z]$ is the fractional part function. Since the infinite sum $u_0 = \sum_{n \in \mathcal{U}} \mu(n)\, n^{-2}$ converges, the claim is immediate from Theorem 7.6. ∎

Often, the maximal subgroup and the multiplicative group coincides: $G \subset \mathbb{Z}_N^*$. This occurs whenever the index $\varphi(n)/\lambda(n) = 1$. For example, $N = 2,\ 4,\ p^m,\ 2\,p^m$, with $p > 2$ an odd prime and $m \geqslant 1$, see Lemma 22.1. Some numerical data is available in [CP14]. Define the subset of integers $\mathcal{R} = \{n \in \mathbb{N} : \varphi(n)/\lambda(n) = 1\}$, and its cardinality $R(x) = \#\{n \leqslant x : \varphi(n)/\lambda(n) = 1\}$. In Theorem 6.9, it is shown that $R(x) \gg x / \log^B x$, where $B > 0$ is a constant.

## 7.4 Average Gap

For $n \geqslant 1$, the gaps $e_n = \lambda(n+1) - \lambda(n)$ between the values of the Carmichael function $\lambda(n)$ should have extreme values of unknown frequencies as the totient function, see Section 5.5 for details.

**Corollary 7.12.** The average gaps of the Carmichael numbers is

$$\overline{\Delta\lambda(n)} = \frac{1}{x} \sum_{n \leqslant x} (\lambda(n+1) - \lambda(n)) = (\log x)^{\eta + o(1)}. \tag{7.18}$$

## 7.5 Squarefree Values

The Carmichael function has a squarefree integer value whenever $n = p_1^{v_1} \cdot p_2^{v_2} \cdots p_k^{v_k}$, with $p_i - 1$ squarefree, and $v_i \leqslant 2$. The counting function $L(x) = \#\{n \leqslant x : \mu(\lambda(n)) \neq 0\}$ of the subset of squarefree values $\mathcal{L} = \{n \in \mathbb{N} : \mu(\lambda(n)) \neq 0\}$ has a nice derivation.

**Lemma 7.13.** ([PP03]) Let $x \geqslant 1$ be a large number. Then, there is a pair of constants $\alpha = \prod_{p \geqslant 2}\left(1 - p^{-1}(p-1)^{-1}\right)$, and $\kappa > 0$ for which

$$\sum_{n \leqslant x} \mu(\lambda(n)) = (\kappa + o(1))\, \frac{x}{\log^{1-\alpha}}. \tag{7.19}$$

## 7.6 Problems And Exercises

**Problem 7.1.** Use $\varphi(n)\frac{\varphi(\lambda(n))}{\lambda(n)}$ the number of primitive roots modulo $n$ to verify that the average gap of the primitive root $r \neq \pm 1,\ r^2$ modulo an integer $n \geqslant 3$ is

$$\overline{g} = \frac{\lambda(n)}{\varphi(\lambda(n))} = \prod_{p \,|\, \lambda(n)} (1 - 1/p)^{-1} \ll \prod_{p \,\leq\, \log n} (1 - 1/p)^{-1} \ll \log\log n. \qquad (7.20)$$

**Problem 7.2.** Prove an asymptotic formula for the average order

$$\sum_{n \leq x} \frac{\lambda(n)}{\varphi(\lambda(n))} = \sum_{n \leq x} \prod_{p \,|\, \lambda(n)} (1 - 1/p)^{-1} ? \qquad (7.21)$$

of integer primitive root gap modulo $n \geq 2$.

**Problem 7.3.** Let $\varphi$ be the Euler totient function, and let $\lambda$ be the Carmichael totient function. Compute the exact asymptotic formula for

$$\sum_{n \leq x} \frac{\varphi(\lambda(n))}{\lambda(n)} \gg \frac{x}{\log\log\log n}\,. \qquad (7.22)$$

**Problem 7.4.** Prove an asymptotic formula for the average order of the index

$$\sum_{n \leq x} \frac{\varphi(n)}{\lambda(n)} = ? \qquad (7.23)$$

of the finite group $\mathbb{Z}_n$. What is the distribution of the sequence of numbers $\varphi(n)/\lambda(n)$? Find an upper bound for the sequence of numbers $\varphi(n)/\lambda(n) \ll (\log n)^{\beta}$, $\beta > 0$ constant?

**Problem 7.5.** Let $\mathcal{U} = \{\, \lambda(n) : n \in \mathbb{N} \,\}$ be a subset of integers. Compute the constant and determine if it is irrational:

$$u_0 = \sum_{n \in \mathcal{U}} \mu(n)\, n^{-2}. \qquad (7.24)$$

**Chapter 8**

# Basic *L*-Functions Estimates

For a character $\chi$ modulo $q \geq 2$, and a complex number $s \in \mathbb{C}$, $\mathcal{R}e(s) > 1/2$, an *L*-function is defined by the infinite sum $L(s, \chi) = \sum_{n \geq 1} \chi(n)\, n^{-s}$. There is a vast literature on the theory of *L*-functions. However, only a few simple elementary estimates associated with the characteristic function of primitive roots, and a few other elementary estimates are calculated here.

## 8.1 An *L*-Function for the Least Primitive Roots

The analysis of the least primitive root mod $p$ is based on the Dirichlet series

$$L(s, \Psi) = \sum_{n \geq 1} \frac{\Psi(n)}{n^s}, \qquad (8.1)$$

where $s \in \mathbb{C}$ is a complex number, and the characteristic function of primitive roots is defined by

$$\Psi(n) = \begin{cases} 1 & \text{if } n \text{ is a primitive root,} \\ 0 & \text{if } n \text{ is not a primitive root,} \end{cases} \qquad (8.2)$$

see Lemmas 3.1 and 3.2 for the exact formulae. The function $L(s, \Psi)$ is zerofree, and analytic on the complex half plane $\{ s \in \mathbb{C} : \mathcal{R}e(s) = \sigma > 1 \}$. Furthermore, since $\sum_{n \leq x} \Psi(n) \gg x / \log x$, for sufficiently large $x \gg p$, see Lemma 3.4, and it has a pole at $s = 1$.

***Lemma* 8.1.** Let $p \geq 2$ be a prime number, and let $s \in \mathbb{C}$ be a complex number, $\mathcal{R}e(s) = \sigma > 1$. Then, there is a nonnegative constant $\kappa_1 > 0$, depending on $p$, such that

$$L(s, \Psi) = \sum_{n \geq 1} \frac{\Psi(n)}{n^s} = \kappa_1. \qquad (8.3)$$

***Proof***: Fix a prime $p \geq 2$, and let $r_1, r_2, \ldots, r_{\varphi(p-1)}$ be the primitive roots mod $p$. Then, the series

$$\sum_{n \geq 1} \frac{\Psi(n)}{n^s} = \frac{1}{r_1^s} + \frac{1}{r_2^s} + \cdots + \frac{1}{r_{\varphi(p-1)}^s} + \frac{1}{(p+r_1)^s} + \frac{1}{(p+r_2)^s} + \cdots \qquad (8.4)$$

converges to a constant $\kappa_1 = \kappa_1(p) > 0$ whenever $\mathcal{R}e(s) = \sigma > 1$ is a real number. $\blacksquare$

***Lemma* 8.2.** Let $x \geq 2$ be a large number. Let $p \geq 2$ be a prime, and let $\chi$ be the characters modulo $d \mid p - 1$ If $s \in \mathbb{C}$ be a complex number, $\mathcal{R}e(s) = \sigma > 1$, then

$$\sum_{n > x} \frac{1}{n^\sigma} \sum_{d \mid p-1} \frac{\mu(d)}{\varphi(d)} \sum_{\text{ord}(\chi) = d} \chi(n) = O\left(\frac{2^{\omega(p-1)}}{x^{\sigma-1}}\right). \qquad (8.5)$$

***Proof***: Assume $s \in \mathbb{C}$ is a complex number, $\mathcal{R}e(s) = \sigma \geq 1$. Rearrange it, to reach the following:

$$\left| \sum_{d \mid p-1} \frac{\mu(d)}{\varphi(d)} \sum_{\text{ord}(\chi) = d} \sum_{n > x} \frac{\chi(n)}{n^s} \right| \leq \sum_{d \mid p-1} \mu^2(d) \left| \sum_{n > x} \frac{\chi(n)}{n^s} \right|$$

$$= O\left( \frac{1}{x^{\sigma-1}} \sum_{d \mid p-1} \mu^2(d) \right). \qquad (8.6)$$

To estimate the infinite sum observe that, see Theorem A-1, in Appendix A,

$$\sum_{n > x} \frac{\chi(n)}{n^s} = \int_x^\infty \frac{1}{t^s} \, d\, R(t) \,, \qquad (8.7)$$

where

$$R(x) = \sum_{n \leq x} \chi(n) = O(x) \,. \qquad (8.8)$$

The claim quickly follows from this observation. $\blacksquare$

***Example* 8.1.** Take the prime $p = 13$, and its $\varphi(p-1) = 4$ primitive roots $r_1 = 2$, $r_2 = 6$, $r_3 = 7$, $r_4 = 11 \mod p$. The numerical value of the series evaluated at $s = 2$ is

$$\sum_{n \geq 1} \frac{\Psi(n)}{n^s} = \frac{1}{2^2} + \frac{1}{6^2} + \frac{1}{7^2} + \frac{1}{11^2} + \frac{1}{(13+2)^2} + \frac{1}{(13+6)^2} + \frac{1}{(13+7)^2} \cdots = 0.321802 \ldots . \qquad (8.9)$$

## 8.2 An L-Function for the Least Prime Primitive Roots

The analysis of the least prime primitive root mod $p$ is based on the Dirichlet series

$$L(s, \Psi\Lambda) = \sum_{n \geq 1} \frac{\Psi(n)\,\Lambda(n)}{n^s}, \tag{8.10}$$

where $s \in \mathbb{C}$ is a complex number, and $\Psi(n)\,\Lambda(n)\,/\,n^s$ is the weighted characteristic function of prime power primitive roots. This is constructed using the characteristic function of primitive roots, which is defined by

$$\Psi(n) = \begin{cases} 1 & \text{if } n \text{ is a primitive root,} \\ 0 & \text{if } n \text{ is not a primitive root,} \end{cases} \tag{8.11}$$

see Lemmas 3.1 and 3.2 for the exact formulae, and the vonMangoldt function

$$\Lambda(n) = \begin{cases} \log p & \text{if } n = p^k, k \geq 1, \\ 0 & \text{if } n \neq p^k, k \geq 1, \end{cases} \tag{8.12}$$

where $p \geq 2$ is a prime. The function $L(s, \Psi\Lambda)$ is zerofree, and analytic on the complex half plane $\{\, s \in \mathbb{C} : \mathcal{R}e(s) = \sigma \geq 1 \,\}$. Furthermore, it has a pole at $s = 1$. This scheme has a lot of flexibility and does not require delicate information on the zerofree regions $\{\, s \in \mathbb{C} : 0 < \mathcal{R}e(s) = \sigma < 1 \,\}$ of the associated $L$-functions. Moreover, the analysis is much simpler than the sieve methods used in the current literature —— lex parsimoniae.

**Lemma 8.3.** Let $p \geq 2$ be a prime number, and let $s \in \mathbb{C}$ be a complex number, $\mathcal{R}e(s) = \sigma > 1$. Then, there is a nonnegative constant $\kappa_2 > 0$, depending on $p$, such that

$$L(s, \Psi\Lambda) = \sum_{n \geq 1} \frac{\Psi(n)\,\Lambda(n)}{n^s} = \kappa_2. \tag{8.13}$$

**Proof:** Fix a prime $p \geq 2$, and let $r_1$, $r_2$, ..., $r_{\varphi(p-1)}$ be the primitive roots mod $p$. Then,

$$\sum_{n \geq 1} \frac{\Psi(n)\,\Lambda(n)}{n^s} = \frac{\Lambda(r_1)}{r_1^s} + \frac{\Lambda(r_2)}{r_2^s} + \cdots + \frac{\Lambda\big(r_{\varphi(p-1)}\big)}{r_{\varphi(p-1)}^s} + \frac{\Lambda(p+r_1)}{(p+r_1)^s} + \frac{\Lambda(p+r_2)}{(p+r_2)^s} + \cdots \tag{8.14}$$

converges to a constant $\kappa_2 = \kappa_2(p) > 0$ whenever $\mathcal{R}e(s) = \sigma > 1$ is a real number. $\blacksquare$

**Lemma 8.4.** Let $x \geq 2$ be a large number. Let $p \geq 2$ be a prime, and let $\chi$ be the characters modulo $d \mid p - 1$ If $s \in \mathbb{C}$ is a complex number, $\mathcal{R}e(s) = \sigma > 1$, then

$$\sum_{n > x} \frac{\Lambda(n)}{n^s} \sum_{d \,|\, p-1} \frac{\mu(d)}{\varphi(d)} \sum_{\text{ord}(\chi) \,=\, d} \chi(n) = O\left(\frac{2^{\omega(p-1)}}{x^{\sigma-1}}\right). \tag{8.15}$$

***Proof***:   Rearranging the absolutely convergent sum yields

$$\left| \sum_{d \,|\, p-1} \frac{\mu(d)}{\varphi(d)} \sum_{\text{ord}(\chi) \,=\, d,} \sum_{n > x} \frac{\chi(n)\,\Lambda(n)}{n^s} \right| \leq \sum_{d \,|\, p-1} \mu^2(d) \left| \sum_{n > x} \frac{\chi(n)\,\Lambda(n)}{n^s} \right|. \tag{8.16}$$

The infinite sum is estimated using, see Theorem A-1, in Appendix A,

$$\sum_{n > x} \frac{\chi(n)\,\Lambda(n)}{n^s} = \int_x^\infty \frac{1}{t^s}\, d\,\psi_\chi(t)\,, \tag{8.17}$$

where

$$\left| \psi_\chi(x) \right| = \left| \sum_{n \leq x} \chi(n)\,\Lambda(n) \right| \ll x\,. \tag{8.18}$$

Estimate the integral to complete the calculation. ∎

***Example*** **8.2.** Take the prime $p = 13$, and its $\varphi(p-1) = 4$ primitive roots $r_1 = 2$, $r_2 = 6$, $r_3 = 7$, $r_4 = 11 \mod p$. The numerical value of the series evaluated at $s = 2$ is

$$\sum_{n \geq 1} \frac{\Psi(n)\,\Lambda(n)}{n^s} = \frac{\log 2}{2^2} + \frac{\log 7}{7^2} + \frac{\log 11}{11^2} + \frac{\log 19}{19^2} + \frac{\log 37}{(2 \cdot 13 + 11)^2} + \frac{\log 41}{(3 \cdot 13 + 2)^2} + \cdots = 0.243611\ldots\,. \tag{8.19}$$

## 8.3 An L-Function for Fixed Primitive Roots

Information on the primes with a fixed primitive root can be encoded in an *L*-series. This concept leads to a Dirichlet series defined by

$$L(s, \psi) = \sum_{\gcd(n, p-1) = 1,\, n \geq 1} \frac{\psi((\tau^n - u)\,k)}{n^s}, \tag{8.20}$$

where $\psi$ is an additive character, and $s \in \mathbb{C}$ is a complex variable, see Lemma 3.2. This series is absolutely convergent on the complex half plane $\{\, s \in \mathbb{C} : \mathcal{R}e(s) = \sigma > 1 \,\}$. An estimate of the average of the partial sum

$$\sum_{\gcd(n,p-1)=1} \frac{\psi((\tau^n - u)\, k)}{n^s} = \sum_{n \le p-1} \frac{\psi((\tau^n - u)\, k)\, \chi_0(n)}{n^s}, \tag{8.21}$$

where $\chi_0 = 1$ is the principal character mod $p - 1$, is the subject of the next result.

**Lemma 8.5.** Let $p \ge 2$ be a large prime, and let $\tau$ be a primitive root mod $p$. For $k, u \in \mathbb{F}_p$, let $\psi(u) = e^{i\, 2\pi k u/p} \ne 1$ be an additive character, and let $s \in \mathbb{C}$, $\mathcal{R}e(s) = \sigma > 1$, be a complex number. If the element $u$ is not a primitive root, then

$$\sum_{n < p,\, \gcd(n,p-1)=1} \frac{1}{n^s} \sum_{0 < k \le p-1} \psi((\tau^n - u)\, k) = -\frac{s}{s-1}\, \frac{\varphi(p-1)}{p-1} + O\left(\frac{1}{p^{s-1}}\right). \tag{8.22}$$

**Proof:** By hypothesis $\tau^n - u \ne 0$ for any integer $n \ge 1$, $\gcd(n,\, p-1) = 1$. Use this information to write the finite sum as

$$\begin{aligned} \sum_{\gcd(n,p-1)=1} \frac{1}{n^s} \sum_{0 < k \le p-1} \psi((\tau^n - u)\, k) &= -\sum_{\gcd(n,p-1)=1} \frac{1}{n^s} \\ &= -\int_1^x \frac{1}{t^s}\, d\, W(t), \end{aligned} \tag{8.23}$$

where $x = p - 1$, and the discrete measure is defined by

$$W(t) = \sum_{\gcd(n,p-1)=1,\, n \le t} 1 = \frac{\varphi(p-1)}{p-1}\, t + O(t^\epsilon), \tag{8.24}$$

where $\epsilon > 0$ is a small number, see Lemma 6.7. Evaluating the integral proves the claim. ∎

Another method for estimating the double finite sum occurring in the nonexistence equations in Theorems 12.1, 13.1, and 13.2 is illustrated here. This technique requires weaker conditions, but the error term is weaker.

**Lemma 8.6.** Let $p \ge 2$ be a large prime, and let $\tau$ be a primitive root mod $p$. If $\psi(u) = e^{i\, 2\pi k u/p}$ is an additive character, with $0 \ne k \in \mathbb{F}_p$ constant, and $s \in \mathbb{C}$, $\mathcal{R}e(s) = \sigma > 1$, is complex number, then

$$\sum_{n < p,\, \gcd(n,p-1)=1} \frac{\psi(\tau^n k)}{n^s} = e^{i\, 2\pi k\, \tau/p} + O\left(\frac{1}{p^{s-1+\epsilon}}\right). \tag{8.25}$$

**Proof:** Let $W(t) = \sum_{\gcd(n,p-1)=1,\, n \le t} \psi(\tau^n k)$. Applying the Abel summation formula, see Theorem A-1, in Appendix A, produces the integral representation

$$\sum_{n < p,\ \gcd(n, p-1)=1} \frac{\psi(\tau^n k)}{n^s} = \frac{W(1)}{1^s} + \int_{n_0}^{p} \frac{1}{t^s}\, d\, W(t)$$

$$= \frac{W(1)}{1^s} + \frac{W(p)}{p^s} + s \int_{n_0}^{p} \frac{W(t)\, d\, t}{t^{s+1}},$$

(8.26)

where $n_0 = \min\{1 < n < p - 1 : \gcd(n,\ p-1) = 1\}$. By Theorem 2.3, the summatory function $W(t)$ has the upper bound

$$W(p) = \sum_{\gcd(n, p-1)=1,\ n \neq 1} \psi(\tau^n k)$$

$$\leq \max_{0 < k \leq p-1} \sum_{v \in H} \psi(k\, v)$$

(8.27)

$$= O\big(p^{1-\epsilon}\big),$$

where $H = \{\tau^n : 1 < n < p - 1,\ \text{and } \gcd(n,\ p - 1) = 1\}$, the element $\tau$ is a primitive root modulo $p$, and $\epsilon > 0$ is an arbitrarily small number. Continuing, this yields

$$\sum_{n < p,\ \gcd(n, p-1)=1} \frac{\psi(\tau^n k)}{n^s} = e^{i\, 2\pi k\, \tau/p} + O\left(\frac{1}{p^{s-1+\epsilon}}\right).$$

(8.28)

This completes the verification. ∎

**Lemma** 8.7. Let $\epsilon > 0$ be a small number. Let $p \geq 2$ be a large prime, and let $\tau$ be a primitive root mod $p$. If $\psi(u) = e^{i\, 2\pi k\, u/p}$ is an additive character, with $0 \neq k \in \mathbb{F}_p$ constant, and $s \in \mathbb{C}$, $\mathcal{R}e(s) = \sigma \geq 2$, is complex number, then

$$\sum_{n < p,\ \gcd(n, p-1)=1} \frac{1}{n^s} \sum_{0 < k \leq p-1} \psi((\tau^n - u)\, k) = -1 + O\left(\frac{1}{p^{s-2+\epsilon}}\right).$$

(8.29)

**Proof**: An approximation of this character sum is computed as follows:

$$\sum_{n < p,\ \gcd(n, p-1)=1} \frac{1}{n^s} \sum_{0 < k \leq p-1} \psi((\tau^n - u)\, k) = \sum_{0 < k \leq p-1} \psi(-u\, k) \sum_{n < p,\ \gcd(n, p-1)=1} \frac{\psi(\tau^n k)}{n^s}$$

$$= \sum_{0 < k \leq p-1} \psi(-u\, k) \left( e^{i\, 2\pi k\, \tau/p} + O\left(\frac{1}{p^{s-1+\epsilon}}\right) \right)$$

(8.30)

$$= -1 + O\left(\frac{1}{p^{s-1+\epsilon}} \left| \sum_{0 < k \leq p-1} \psi(-u\,k) \right| \right)$$

$$= -1 + O\left(\frac{1}{p^{s-2+\epsilon}}\right),$$

where $\epsilon > 0$ is a small constant. The second line on right side of the (8.30) follows from Lemma 8.6. And the fourth line utilizes the trivial estimate $\left| \sum_{0 < k \leq p-1} \psi(-u\,k) \right| \leq p$. ∎

**Chapter 9**

# Least Primitive Roots

Let $N \in \mathbb{N}$ be an integer, and let $G$ be a finite group of order $\#G = N$. A group is cyclic if it has a generator, called a primitive root, such that $G = \{ g^n : 0 \leq n < N \}$. A few centuries ago, continuing the works of his predecessors, Gauss proved that a group is cyclic if and only if the order $N = \varphi(2)$, $\varphi(4)$, $\varphi(p^m)$, $\varphi(2\,p^m)$, with $p > 2$ an odd prime and $m \geq 1$, [AP86, p. 209], [NW00, p. 15-23], [RE94].

The order of magnitude of the smallest primitive root is a topic of much interest in theory and practice. The current literature has the Burgess result

$$g(p) \ll p^{1/4+\epsilon}, \tag{9.1}$$

with $\epsilon > 0$ an arbitrarily small number, see [BR62], and a few other almost identical results as the least primitive root modulo an odd prime $p \geq 3$. The average value satisfies

$$\pi(x)^{-1} \sum_{p \leq x} g(p) \ll (\log x)^2 \, (\log\log x)^4 \,, \tag{9.2}$$

see [BE68], [BS71]. Moreover, there are a few conditional estimates such as

$$g(p) \ll \omega(p-1)^4 \log(\omega(p-1)+1)\,(\log p)^2 \tag{9.3}$$

obtained in [SP92]. And the conjectured upper bound $g(p) \ll (\log p)\,(\log\log p)^2$ is derived in [BH97]. The conjectured upper bound is close to the Turan lower bound $g(p) = \Omega(\log p \log\log p)$, refer to [BE68], [RN96, p. 24], [MT91] for discussions.

## 9.1 The Proof

The goal of this Section is to present some effective estimates of the least primitive root in a cyclic group $G$ of order $\#G = \varphi(p)$ using basic techniques based on an $L$-function associated with primitive root. A primitive root of a prime $p$ or $p^2$ lifts to a primitive root modulo $p^m$, and $2\,p^m$ for all $m \geq 1$, so it is sufficient to consider only the primes $p$ or $p^2$.

**Theorem 9.1.** Let $p \geq 3$ be a large prime. Then the followings holds.

(i) Almost every large prime $p$ has a primitive root $g(p) \ll (\log p)\,(\log\log p)$.

(ii) Every large prime $p$ has a primitive root $g(p) \ll p^{5/\log\log p}$.

**Proof**: Fix a large prime $p \geq 2$. Let $p - 1 = \varphi(p)$, and consider summing the weighted characteristic function $\Psi(n) / n^s$ over a small range of integers $n \leq x$, see Lemma 3.1. Then, the nonexistence equation

$$\sum_{n \leq x} \frac{\Psi(n)}{n^s} = \sum_{n \leq x} \frac{1}{n^s} \left( \frac{\varphi(p-1)}{p-1} \sum_{d \,|\, p-1} \frac{\mu(d)}{\varphi(d)} \sum_{\mathrm{ord}(\chi) = d} \chi(n) \right) = 0 \,, \tag{9.4}$$

where $s \in \mathbb{C}$ is a complex number, $\mathcal{R}e(s) \geq 1$, holds if and only if there are no primitive roots in the interval $[1, \ x]$.

Applying Lemma 8.1, and Lemma 8.2 yield, with $s = 2$, yield

$$\begin{aligned}
0 &= \sum_{n \leq x} \frac{1}{n^s} \left( \sum_{d \,|\, p-1} \frac{\mu(d)}{\varphi(d)} \sum_{\mathrm{ord}(\chi) = d} \chi(n) \right) \\
&= \sum_{n \geq 1} \frac{\Psi(n)}{n^s} \ - \sum_{n > x} \frac{\Psi(n)}{n^s} \\
&= \kappa_1 + O\!\left( \frac{2^{\omega(p-1)}}{x} \right),
\end{aligned} \tag{9.5}$$

where $\kappa_1 = \kappa_1(p) > 0$ is a constant, which depends on the prime $p \geq 2$.

**Case (i): $2^{\omega(p-1)} \ll \log p$. Restriction to the Average Integers $p - 1$, Refer to Lemmas 4.2 and 4.3.**

Let $x = (\log p)(\log\log p)$, and suppose that the short interval $[2, \ (\log p)(\log\log p)]$ does not contain primitive roots.

Then, replacing these information in (9.3) yields

$$\begin{aligned}
0 &= \sum_{n \leq x} \frac{1}{n^s} \left( \sum_{d \,|\, p-1} \frac{\mu(d)}{\varphi(d)} \sum_{\mathrm{ord}(\chi) = d} \chi(n) \right) \\
&= \kappa_1 + O\!\left( \frac{2^{\omega(p-1)}}{x} \right) \\
&= \kappa_1 + O\!\left( \frac{1}{\log\log p} \right) > 0,
\end{aligned} \tag{9.6}$$

where $\kappa_1 > 0$ is a constant. But this is a contradiction for all sufficiently large $p \geq 2$.

**Case (ii): $2^{\omega(p-1)} \ll p^{4/\log\log p}$. No Restrictions on the Integers $N$, Refer to Lemmas 4.2 and 4.3.**

Let $x \leq p^{5/\log\log p}$, and suppose that the short interval $\left[2, \ p^{5/\log\log p}\right]$ does not contain primitive roots.

Then, replacing these information in (9.3) yield

$$
\begin{aligned}
0 &= \sum_{n \leq x} \frac{1}{n^s} \left( \sum_{d \mid p-1} \frac{\mu(d)}{\varphi(d)} \sum_{\mathrm{ord}(\chi) = d} \chi(n) \right) \\
&= \kappa_1 + O\left( \frac{2^{\omega(p-1)}}{x} \right) \\
&= \kappa_1 + O\left( \frac{1}{p^{1/\log\log p}} \right) > 0,
\end{aligned}
\tag{9.7}
$$

where $\kappa_1 > 0$ is a constant. But this is a contradiction for all sufficiently large $p \geq 2$. $\blacksquare$

Case (i) above shows why primitive roots are very frequently very small integers. Equation (9.2) offers lot flexibility in the selection of the complex number $s \in \mathbb{C}$. The selection of a complex number in the zerofree region $\mathcal{R}e(s) = \sigma \geq 1$ is probably simpler since the analysis is well established.

The pair of graphs below displays the least primitve root $g \geq 2$ modulo the $n$th prime $p_n$ for two range of integers.

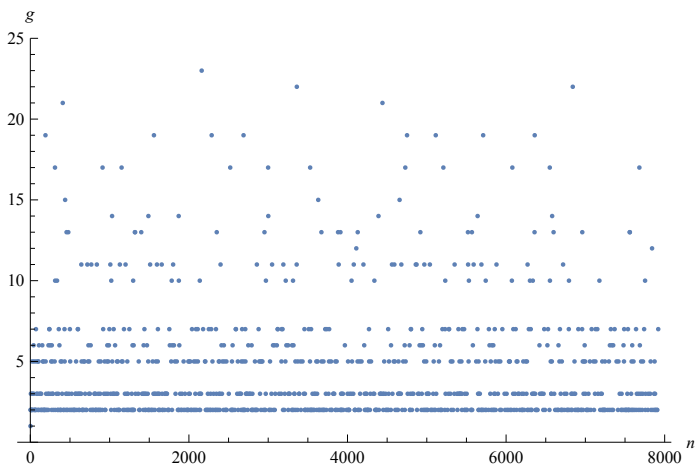

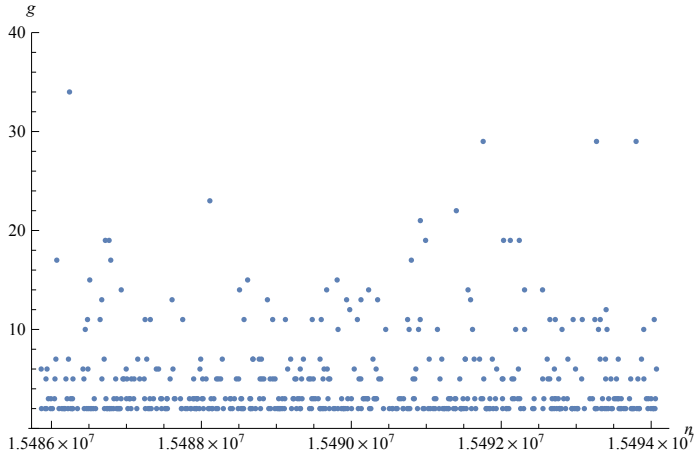

## 9.2 Problems And Exercises

**Problem 9.1.** Let $p \geq 2$ be a large prime number. Improve the known average least primitive root $\overline{g(p)} = \pi(x)^{-1} \sum_{p \leq x} g(p) \ll (\log x)^2 (\log\log x)^4$, see [BE68], [BS71].

**Problem 9.2.** Let $p \geq 2$ be a large prime number. Compute the variance of the least primitive root $V(p) = \pi(x)^{-1} \sum_{p \leq x} (g(p) - \overline{g}(p))^2$?

**Problem 9.3.** Let $p \geq 2$ be a large prime number. Determine the distribution of the least primitive root modulo $p$, Poisson?, Gaussian? or otherwise?

**Problem 9.4.** Give an easy proof of the average gap between primitive roots mod $p$, show that it is quite small $\frac{p-1}{\varphi(p-1)} \ll \log p$.

**Chapter 10**

# The Least Prime Primitive Roots

The current literature has several estimates of the least prime primitive root $g^*(p)$ modulo a prime $p \geq 2$ such as $g^*(p) \ll p^c$, $c > 2.8$. The actual constant $c > 2.8$ depends on various conditions such as the factorization of $p - 1$, et cetera. These results are based on sieve methods and the least primes in arithmetic progressions, see [MA98], [MB97], [HA13]. Moreover, there are a few other conditional estimates such as $g(p) \leq g^*(p) \ll (\log p)^6$, see [SP92], and the conjectured upper bound $g^*(p) \ll (\log p)(\log\log p)^2$. This is implied by the conjectured limit supremum

$$\limsup_{n \to \infty} \frac{g^*(p)}{\log p (\log\log p)^2} = e^\gamma \tag{10.1}$$

derived in [BH97]. An easy proof of the estimate $g^*(p) \ll p$, conditional on the GRH, is sketched in [RM08, p. 400]. On the other direction, there is the Turan lower bound $g^*(p) \geq g(p) = \Omega(\log p \log\log p)$, refer to [BE68], [RN96, p. 24], and [MT91, p. 48] for discussions.

## 10.1 The Proof

It is expected that there are estimates for the least prime primitive roots which are quite similar to the estimates for the least primitive roots. This is a very small quantity and nowhere near the currently proved results. In fact, the numerical tables confirm that the least prime primitive roots are very small, and nearly the same magnitude as the least primitive roots, but have more complex patterns, see [PA02], [PB02].

**Theorem 10.1.** Let $p \geq 3$ be a large prime. Then the followings holds.

(i) Almost every prime $p \geq 3$ has a prime primitive root $g^*(p) \ll (\log p)^c$, $c > 1$ constant.

(ii) Every prime $p \geq 3$ has a prime primitive root $g^*(p) \ll p^{5/\log\log p}$.

**Proof**: Fix a large prime $p \geq 2$, and consider summing the weighted characteristic function $\Psi(n) \Lambda(n) / n^s$ over a small range of prime powers $q^k \leq x$, $k \geq 1$, see Lemma 3.1. Then, the nonexistence equation

$$\sum_{n \leq x} \frac{\Psi(n) \Lambda(n)}{n^s} = \sum_{n \leq x} \frac{\Lambda(n)}{n^s} \left( \frac{\varphi(p-1)}{p-1} \sum_{d \mid p-1} \frac{\mu(d)}{\varphi(d)} \sum_{\mathrm{ord}(\chi) = d} \chi(n) \right) = 0, \tag{10.2}$$

where $s \in \mathbb{C}$ be a complex number, $\mathcal{R}e(s) = \sigma \geqslant 1$, holds if and only if there are no prime power primitive roots in the interval $[1, x]$.

Applying Lemma 8.3, and Lemma 8.4 with $q = p - 1$, and $s = 2$, yield

$$
\begin{aligned}
0 &= \sum_{n \leq x} \frac{\Lambda(n)}{n^s} \left( \frac{\varphi(p-1)}{p-1} \sum_{d \mid p-1} \frac{\mu(d)}{\varphi(d)} \sum_{\text{ord}(\chi) = d} \chi(n) \right) \\
&= \sum_{n \geq 1} \frac{\Psi(n)\,\Lambda(n)}{n^s} \;-\; \sum_{n > x} \frac{\Psi(n)\,\Lambda(n)}{n^s} \\
&= \kappa_2 + O\left( \frac{2^{\omega(p-1)}}{x} \right),
\end{aligned}
\tag{10.3}
$$

where $\kappa_2 = \kappa_2(p) > 0$ is a constant, which depends on the prime $p \geqslant 2$.

**Case (i): $2^{\omega(p-1)} \ll \log p$. Restriction to the Average Integers $p - 1$, Refer to Lemmas 4.2 and 4.3.**

Let $x = (\log p)^{1+\epsilon}$, $\epsilon > 0$, and suppose that the short interval $\left[2, (\log p)^{1+\epsilon}\right]$ does not contain prime primitive roots. Then, replacing these information in (10.3) yields

$$
\begin{aligned}
0 &= \sum_{n \leq x} \frac{\Lambda(n)}{n^2} \left( \sum_{d \mid p-1} \frac{\mu(d)}{\varphi(d)} \sum_{\text{ord}(\chi) = d} \chi(n) \right) \\
&= \kappa_2 + O\left( \frac{2^{\omega(p-1)}}{x} \right) \\
&= \kappa_2 + O\left( \frac{1}{(\log p)^{\epsilon}} \right) > 0.
\end{aligned}
\tag{10.4}
$$

Since $\kappa_2 > 0$ is a constant, this is a contradiction for all sufficiently large prime $p \geqslant 3$.

**Case (ii): $2^{\omega(p-1)} \ll p^{4/\log\log p}$. No Restrictions on the Integers $p - 1$, Refer to Lemmas 4.2 and 4.3.**

Let $x = p^{5/\log\log p}$, and suppose that the short interval $\left[2, p^{5/\log\log p}\right]$ does not contain prime primitive roots. Then, replacing these information in (10.3) yield

$$
\begin{aligned}
0 &= \sum_{n \leq x} \frac{\Lambda(n)}{n^2} \left( \sum_{d \mid p-1} \frac{\mu(d)}{\varphi(d)} \sum_{\text{ord}(\chi) = d} \chi(n) \right) \\
&= \kappa_2 + O\left( \frac{2^{\omega(p-1)}}{x} \right) \\
&= \kappa_2 + O\left( \frac{1}{p^{1/\log\log p}} \right) > 0.
\end{aligned}
\tag{10.5}
$$

Since $\kappa_2 > 0$ is a constant, this is a contradiction for all sufficiently large prime $p \geqslant 3$. $\blacksquare$

## 10.2 Problems And Exercises

**Problem 10.1.** Let $p \geq 2$ be a large prime number. Estimate the average least prime primitive root $\overline{g^*(p)} = \pi(x)^{-1} \sum_{p \leq x} g^*(p)$?, this appears to be unknown.

**Problem 10.2.** Let $p \geq 2$ be a large prime number. Compute the variance of the least prime primitive root $V^*(p) = \pi(x)^{-1} \sum_{p \leq x} \left(g^*(p) - \overline{g^*(p)}\right)^2$?, this appears to be unknown.

**Problem 10.3.** Let $p \geq 2$ be a large prime number. Determine the distribution of the least prime primitive root modulo $p$, Poisson?, Gaussian? or otherwise? this appears to be unknown.

**Chapter 11**

# Artin Conjecture

Let the symbol $\mathbb{Z} = \{0, \pm 1, \pm 2, \dots\}$ denotes the set of integers, and the symbol $\mathbb{P} = \{2, 3, 5, \dots\}$ denotes the set of prime numbers. The subset of primes $\mathcal{P}_g = \{p \in \mathbb{P} : \mathrm{ord}\,(g) = p - 1\} \subset \mathbb{P}$ with a fixed primitive root $g \in \mathbb{Z}$, is expected to have a nonzero density $\delta(g) > 0$.

**Conjecture 11.1.** Let $g \neq \pm 1$, $b^2$ be a fixed integer. Then, $g$ is a primitive root modulo $p$ for infinitely many primes $p \geq 2$ as $x \to \infty$. In particular,

$$\pi_g(x) = \#\{p \leq x : \mathrm{ord}(g) = p - 1\} = \delta(g)\,\mathrm{li}(x) + O\big(x(\log x)^{-2}\big) \tag{11.1}$$

where $\mathrm{li}(x)$ is the logarithm integral, as $x \to \infty$.

## 11.1 Conditional Result

A conditional proof of this result was achieved in [HY67], and a simplified sketch of the proof appears in [MP04, p. 8]. The determination of the constant $\delta(g) \geq 0$ for a fixed integer $g \in \mathbb{Z}$ is an interesting technical subject, [HY67, p. 218], [LA11], [LB77], et alii. An introduction to its historical development, and its calculations is covered in [MP04, p. 3-10], [SN03].

**Theorem 11.2.** ([HY67]) If it be assumed that the extended Riemann hypothesis hold for the Dedekind zeta function over Galois fields of the type $\mathbb{Q}\big(\sqrt[n]{a}, \sqrt[d]{1}\,\big)$, where $n$ is a squarefree integer, and $d \mid n$. For a given nonzero integer $a \neq \pm 1$, $b^2$, let $N_a(x)$ be the number of primes $p \leq x$ for which $a$ is a primitive root modulo $p$. Let $a = a_1^m \cdot a_2^{2m}$, where $a_1 > 1$ is square-free, and $m \geq 1$ is odd, and let

$$C(m) = \prod_{q \mid m}\left(1 - \frac{1}{q-1}\right)\prod_{q \nmid m}\left(1 - \frac{1}{q(q-1)}\right). \tag{11.2}$$

(i) If $a \neq \pm 1$, $b^2$ is a fixed integer, then there are infinitely many primes $p$ for which $a$ is a primitive root modulo $p$.

(ii) If $a_1 \not\equiv 1 \bmod 4$, then

$$N_a(x) = C(m) \frac{x}{\log x} + O\left(\frac{x \log\log x}{\log^2 x}\right) \quad \text{as } x \longrightarrow \infty \,. \tag{11.3}$$

While if $a_1 \equiv 1 \bmod 4$, then

$$N_a(x) = C(m)\left(1 - \mu(\,|\,a_1\,|) \prod_{q\,|\,m,\,q\,|\,a_1} \frac{1}{q-2} \prod_{q\,|\,m,\,q\,\nmid\,a_1} \frac{1}{q^2+q-1}\right) \frac{x}{\log x} + O\left(\frac{x \log\log x}{\log^2 x}\right) \quad \text{as } x \longrightarrow \infty \,. \tag{11.4}$$

There are a few equivalent formulae for the density $\delta(g) \geq 0$. The average value of the density is the Artin constant

$$A_1 = \sum_{n \geq 1} \frac{\mu(n)}{n\,\varphi(n)} = \prod_{p \geq 2}\left(1 - \frac{1}{p(p-1)}\right) = .3739558136\ldots\,. \tag{11.5}$$

corresponds to $A_1 = C(1)$.

## 11.3 Average Results

The Artin primitive root conjecture on average was proved in [GD68] unconditionally, and refined in [ST69]. These works had shown that almost all admissible integers $g \in \mathbb{Z} - \{-1,\ 1,\ b^2 : b \in \mathbb{Z}\}$ are primitive roots for infinitely many primes.

**Theorem** 11.2 ([GD68], [ST69]) (i) If $N > e^{4\sqrt{\log x \log\log x}}$, then the average number of primes $p \leq x$ with a fixed primitive root is

$$N^{-1} \sum_{g \leq N} \pi_g(x) = A \operatorname{li}(x) + O\big(x\,(\log x)^{-D}\big), \tag{11.6}$$

where $A = \prod_{p \geq 2} \left( 1 - p^{-1}(p-1)^{-1} \right)$ is Artin constant, and $D > 1$ is an arbitrary constant.

(ii) If $N > e^{6\sqrt{\log x \log\log x}}$, then the variance of the number of primes $p \leq x$ with a fixed primitive root is

$$N^{-1} \sum_{g \leq N} \left( \pi_g(x) - A\,\mathrm{li}(x) \right)^2 \ll x^2 \left( \log x \right)^{-E}, \tag{11.7}$$

where $E > 2$ is an arbitrary constant.

The number of exceptions is a subset of zero density. The individual quantity $\pi_g(x)$ in (11.2) can be slightly different from the average quantity in (11.4). The variations, discovered by the Lehmers using numerical experiments, depend on the primes decomposition of the fixed value $g$, see [HY67, p. 220], [MO, p. 3] and similar references for the exact formula for the density $\delta(g)$.

## 11.4 Current Results

The closer approximations to the Artin primitive root conjecture was established a few decades ago. The current results in the literature have reduced the number of possible exceptions to a small finite set. For example, in [GM84] it was proved that for a fixed primes triple $q,\ r,\ s$, the subset of integers

$$\mathcal{A}(q,r,s) = \left\{ q^a\, r^b\, s^c : 0 \leq a,\, b,\, c \leq 3 \right\} \tag{11.8}$$

contains a primitive root for infinitely many primes.

***Definition*** **11.3.** A subset of integers $\{\, q, r, s \,\}$ are multiplicative independent if $q^a\, r^b\, s^c = 1$, with $a,\, b,\, c \in \mathbb{Z}$, then $a = b = c = 0$.

For a triple of distinct integers $q, r, s \in \mathbb{Z}$, let $S = \left\{\, q\,s^2,\ q^2\,r^2,\ q^2\,r,\ r^3\,s^2,\ r^2\,s,\ q^2\,s^3,\ q\,r^3,\ q^3\,r\,s^2,\ r\,s^2,\ r\,s^3,\ q^2\,r^3\,s,\ q^3\,s,\ q\,r^2\,s^3,\ q\,r\,s \,\right\}$.

***Theorem*** **11.4.** ([GM84]) For some $a \in S$, there is a real number $\delta > 0$ such that for at least $x \log^{-2} x$ primes $p \leq x$, the integer $a$ is a primitive root mod $p$.

This result was later reduced to the smaller subset $\mathcal{B}(q,r,s) = \left\{\, q^a\, r^b\, s^c : 0 \leq a,\, b,\, c \leq 1 \,\right\}$, see [HB86].

***Theorem*** **11.5.** ([HB86]) Let $q, r, s$ be three nonzero integers which are multiplicatively independent. Suppose that none of $q, r, s, -3qr, -3qs, -3rs, qrs$ is a square. Then the number $N(x)$ of primes $p \leq x$ for which at least one of $q, r, s$ is a primitive root mod $p$ satisfies $N(x) \gg x(\log x)^{-2}$ as $x \to \infty$.

These results have been extended to quadratic numbers fields in [CJ07], [NW00], [RH02], [AM14], and cyclic groups of elliptic curves in [GM84]. Other related results are given in [AP14], [FX11], [PA02], et alii.

In both of these results, the lower bound for the number of primes $p \leq x$ with a fixed primitive root in either of the subset $\mathcal{A}(q, r, s)$ or $\mathcal{B}(q, r, s)$ is $\#\{ p \leq x : \mathrm{ord}(g) = p - 1 \} \gg x (\log x)^{-2}$ for all large number $x \geq 1$. The technique explored in Chapter 12 provides the expected lower bound for the number of primes $p \leq x$ with a fixed primitive root $g$ mod $p$ for all large number $x \geq 1$. Furthermore, the fixed primitive root $g$ is not restricted to a small finite subset such as $\mathcal{A}(q, r, s)$ or $\mathcal{B}(q, r, s)$.

**Chapter 12**

# Subsets of Primes with Fixed Primitive Roots

The representations of the characteristic function of primitive roots, Lemma 3.1, and Lemma 3.2, easily detect certain local and global properties of the elements $u \in \mathbb{F}_p$ in a finite field. Exempli gratia, it vanishes

$$\Psi(u) = 0 \qquad \text{if } u = \pm 1, v^2. \tag{12.1}$$

Ergo, the constraint $u \neq \pm 1, v^2$ is a necessary global condition to be a primitive element in $\mathbb{F}_p$, $p \geq 2$ an arbitrary prime. The requirement of being an $n$th power nonresidue mod $p$ or primitive root are local properties.

## 12.1 A Proof of the Result

A summation technique utilizing a new weighted divisors-free representation of the characteristic function appears to be effective in solving the fixed primitive root/variable primes counting problem.

**Theorem 12.1.** A fixed integer $u \neq \pm 1, v^2$ is a primitive root mod $p$ for infinitely many primes $p \geq 2$. In addition, the density of these primes satisfies

$$\#\{ p \leq x : \text{ord}(u) = p - 1 \} = \alpha_u \, \text{li}(x) + o(\text{li}(x)), \tag{12.2}$$

where $\text{li}(x)$ is the logarithm integral, and $\alpha_u > 0$ is a constant for all large numbers $x \geq 1$.

**Proof**: Suppose that $u \neq \pm 1, v^2$ is not a primitive root for all primes $p \geq x_0$, with $x_0 \geq 1$ constant. Let $x > x_0$ be a large number, and consider the nonexistence equation

$$0 = \sum_{x \leq p \leq x^2} \left( \frac{1}{p} \sum_{\gcd(n, p-1)=1} \frac{1}{n^s} \sum_{0 \leq k \leq p-1} \psi((\tau^n - u) k) \right), \tag{12.3}$$

where $\tau \in \mathbb{F}_p$ is a primitive root mod $p$, and $\psi \neq 1$ is an additive character of order ord $\psi = p$, see Lemma 3.2 for details on the characteristic function. The complex variable $s \in \mathbb{C}$ will be set to $s \geq 2$.

Regroup the triple finite sum as

$$0 = \sum_{x \le p \le x^2} \frac{1}{p} \sum_{\gcd(n,\,p-1)=1} \frac{1}{n^s} + \sum_{x \le p \le x^2} \frac{1}{p} \left( \sum_{\gcd(n,\,p-1)=1} \frac{1}{n^s} \sum_{0 < k \le p-1} \psi((\tau^n - u)\,k) \right). \qquad (12.4)$$

The main term is determined by a finite sum over the trivial additive character $\psi = 1$, and the error term is determined by a finite sum over the nontrivial additive characters $\psi = e^{i\,2\pi\,k/p} \ne 1$.

Put $s = 2$. Applying Lemma 8.5, the trivial upper bound $\varphi(p-1)\,/\,(p-1) \le 1\,/\,2$ in Lemma 6.2, and simplifying return

$$
\begin{aligned}
0 &= \sum_{x \le p \le x^2} \frac{1}{p} \sum_{\gcd(n,\,p-1)=1} \frac{1}{n^s} + \sum_{x \le p \le x^2} \frac{1}{p} \left( -\frac{s}{s-1} \frac{\varphi(p-1)}{p-1} + O\!\left( \frac{1}{p^{s-1}} \right) \right) \\
&\ge \sum_{x \le p \le x^2} \frac{1}{p} \sum_{\gcd(n,\,p-1)=1} \frac{1}{n^2} + \sum_{x \le p \le x^2} \left( -\frac{1}{p} + O\!\left( \frac{1}{p^2} \right) \right) \\
&\ge \sum_{x \le p \le x^2} \frac{1}{p} \sum_{\gcd(n,\,p-1)=1,\; n \ge 3} \frac{1}{n^2} + O\!\left( \sum_{x \le p \le x^2} \frac{1}{p^2} \right) \\
&> \frac{1}{3^2} \sum_{\substack{x \le p \le x^2, \\ p \equiv 5 \bmod 6}} \frac{1}{p} + O\!\left( \sum_{x \le p \le x^2} \frac{1}{p^2} \right) \\
&\ge \frac{1}{2} \frac{1}{3^2} \log(2) + O\!\left( \frac{1}{\log x} \right) \\[4pt]
&> 0.
\end{aligned}
\qquad (12.5)
$$

The last two inequalities follow from the fact that the primes in the arithmetic progression $\{\, p = 6\,m + 5 : m \ge 1 \,\}$ has density $1\,/\,2$ in the set of primes, and $\gcd(3,\,p-1) = 1$ for any prime $p = 6\,m + 5$. Confer Lemma 17.3 for information on the prime harmonic sum $\sum_{p \le x} 1\,/\,p \sim \log\log x$.

Clearly, this is a contradiction for all sufficiently large numbers $x \ge x_0$. Hence, the short interval $\left[x,\,x^2\right]$ contains primes $p \ge x$ such that $u \ne \pm 1$, $v^2$ is a primitive root mod $p$.

To compute the asymptotic primes counting formula, take

$$\sum_{p \le x} \frac{1}{p} \sum_{\gcd(n,\,p-1)=1,\; 0 \le k \le p-1} \psi((\tau^n - u)\,k) = \sum_{p \le x} \frac{\varphi(p-1)}{p} + o\!\left( \sum_{p \le x} \frac{\varphi(p-1)}{p} \right). \qquad (12.6)$$

The right side (12.6) follows from (12.5), which implies that the main term in (12.6) is larger than the error term. Then, by Lemma 6.6, and partial summation, the density of the corresponding subset of primes over the short interval $[1, x]$ is

$$\sum_{p \leq x} \frac{\varphi(p-1)}{p} + o\left(\sum_{p \leq x} \frac{\varphi(p-1)}{p}\right) = \sum_{p \leq x} \frac{p-1}{p} \frac{\varphi(p-1)}{p-1} + o\left(\sum_{p \leq x} \frac{p-1}{p} \frac{\varphi(p-1)}{p-1}\right)$$

$$= \int_1^x \frac{t-1}{t} \, d \, R(t) + o\left(\int_1^x \frac{t-1}{t} \, d \, R(t)\right) \tag{12.7}$$

$$= \alpha_u \operatorname{li}(x) + o(\operatorname{li}(x)),$$

where $\operatorname{li}(x) = \int_1^x (\log t)^{-1} \, d \, t$ is the logarithm integral, and

$$R(x) = \sum_{p \leq x} \frac{\varphi(p-1)}{p-1}, \tag{12.8}$$

see Lemma 5.6, and $\alpha_u > 0$ is a constant, for all large numbers $x \geq 1$.

The determination of the constant $\alpha_u = \delta(u) > 0$, which is the density of primes with a fixed primitive root $u$, is a very technical subject in algebraic number theory, see [HY67, p. 218], [LA11], [LB77], [MP04, p. 10], [SN92]. The formula, without proof, appears in (11.3-5) above.

Some calculations of the constants for some cases of primitive roots in quadratic fields are given in [RH02], [CJ07], et cetera. Other calculations of the constants for the related cases for elliptic primitive roots are given in [BC11], et alii.

## 12.2 The Densities Of Fixed Primitive Roots In Arithmetic Progressions

For a fixed pair of integers $a < q$, $\gcd(a, q) = 1$, consider the arithmetic progression of primes $\{\, p \in \mathbb{P} : p \equiv a \bmod q \,\}$. Little changes in the proof of Theorem 12.1 yields a result for the primes in the arithmetic progressions with a fixed primitive root.

**Theorem 12.2.** Let $a < q$, $\gcd(a, q) = 1$, be a small parameters. A fixed integer $u \neq \pm 1$, $v^2$ is a primitive root mod $p$ for infinitely many primes $p \equiv a \bmod q$. In addition, the density of these primes satisfies

$$\#\{\, p \leq x : \operatorname{ord}(u) = p-1 \text{ and } p \equiv a \bmod q \,\} = \alpha_u(q, a) \, \frac{\varphi(q)}{q} \operatorname{li}(x) + o(\operatorname{li}(x)), \tag{12.9}$$

where $\operatorname{li}(x)$ is the logarithm integral, and $\alpha_u(q, a) > 0$ is a constant for all large numbers $x \geq 1$.

The topic of densities of primes with fixed primitive roots in arithmetic progressions is a very interesting topic in algebraic number theory. Various results have been derived in [MP07].

## 12.3 The Densities of A Few Fixed Primitive Roots

The counting functions of the first few fixed primitive roots are computed in this Section. These calculations are based on the well established formulae for the constants. The average density is given by the Artin constant

$$\alpha_1 = \sum_{n \geq 1} \frac{\mu(n)}{n \, \varphi(n)} = \prod_{p \geq 2} \left( 1 - \frac{1}{p(p-1)} \right) = .3739558136 \ldots . \tag{12.10}$$

The average fixed primitive root has an average constant. The dependence of the constant on the fixed primitive root was discovered by the Lehmers in numerical experiments, see [SN03]. The corrected formula in (11.5) is the product of about a half a century of works by scores of mathematicians.

$$\tag{12.11}$$

(1) $\pi_2(x) = \# \{ \, p \leq x : \operatorname{ord}(2) = p - 1 \, \} = \alpha_1 \operatorname{li}(x) + O\!\left( x (\log x)^{-2} \right).$

(2) $\pi_3(x) = \# \{ \, p \leq x : \operatorname{ord}(3) = p - 1 \, \} = \alpha_1 \operatorname{li}(x) + O\!\left( x (\log x)^{-2} \right).$

(3) $\pi_4(x) = \# \{ \, p \leq x : \operatorname{ord}(4) = p - 1 \, \} = 0 \, .$

(4) $\pi_5(x) = \# \{ \, p \leq x : \operatorname{ord}(5) = p - 1 \, \} = \dfrac{20}{19} \alpha_1 \operatorname{li}(x) + O\!\left( x (\log x)^{-2} \right), \ \text{see formula (11.4).}$

(5) $\pi_6(x) = \# \{ \, p \leq x : \operatorname{ord}(6) = p - 1 \, \} = \alpha_1 \operatorname{li}(x) + O\!\left( x (\log x)^{-2} \right).$

(6) $\pi_7(x) = \# \{ \, p \leq x : \operatorname{ord}(7) = p - 1 \, \} = \alpha_1 \operatorname{li}(x) + O\!\left( x (\log x)^{-2} \right).$

(7) $\pi_8(x) = \# \{ \, p \leq x : \operatorname{ord}(8) = p - 1 \, \} = \dfrac{4}{5} \alpha_1 \operatorname{li}(x) + O\!\left( x (\log x)^{-2} \right), \ \text{see formula (11.4).}$

(8) $\pi_9(x) = \# \{ \, p \leq x : \operatorname{ord}(9) = p - 1 \, \} = 0.$

(9) $\pi_{10}(x) = \# \{ \, p \leq x : \operatorname{ord}(10) = p - 1 \, \} = \alpha_1 \operatorname{li}(x) + O\!\left( x (\log x)^{-2} \right).$

## 12.4 Nonzero Rational Numbers As Primitive Roots

The sufficient condition $u \neq \pm 1$, $s^2$ of Artin primitive root conjecture can be simplified to $u \neq 0$. Observe that $u + u^{-1} = \pm 1$ or $x^2$ has no integer solutions $u \in \mathbb{Z}$. To verify this, write the discriminant of the resulting Diophantine equation as $y^2 = x^4 - 4$, and show that it has no integral solutions in the set of integers $\mathbb{Z}$. This idea seems to simplify the condition of Artin primitive root conjecture. The precise idea is as follows.

**Theorem 12.3** For any fixed integer $u \neq 0$, the rational number $v = u + u^{-1}$ is a primitive root mod $p$ for infinitely many primes $p \geq 2$. In addition, the density of these primes satisfies

$$\#\{\, p \leq x : \operatorname{ord}(v) = p - 1 \,\} = C_v \operatorname{li}(x) + o(\operatorname{li}(x)) \,, \tag{12.12}$$

where $C_v > 0$ is a constant for all large numbers $x \geq 1$.

## 12.5 Problems And Exercises

**Problem 12.1.** Compute the folowing statistics:

$\pi_2(x) = \#\{\, p \leq x : \operatorname{ord}(2) = p - 1 \,\}$, for $x = 10^3$, $10^6$ $10^9$ $10^{12}$.

$\pi_3(x) = \#\{\, p \leq x : \operatorname{ord}(3) = p - 1 \,\}$, for $x = 10^3$, $10^6$ $10^9$ $10^{12}$.

$\pi_5(x) = \#\{\, p \leq x : \operatorname{ord}(2) = p - 1 \,\}$, for $x = 10^3$, $10^6$ $10^9$ $10^{12}$.

$\pi_6(x) = \#\{\, p \leq x : \operatorname{ord}(2) = p - 1 \,\}$, for $x = 10^3$, $10^6$ $10^9$ $10^{12}$.

**Problem 12.2.** Explain why the integers $2, 3, 5, 6, \ldots, u \neq v^2, u \equiv 1 \bmod 4$, have the same asymptotic order $\#\{\, p \leq x : \operatorname{ord}(u) = p - 1 \,\} = \alpha_2 \operatorname{li}(x) + o(\operatorname{li}(x))$, but the numerical data shows that 2 has a higher frequency of occurence.

**Problem 12.3.** Assume that a finite group $G$ of order $N = \# G$ can be generated by a subset $\{g_1, g_2, \ldots, g_M\} \subset G$ of at most $M = O(\log N)$ elements. For example, $< g_1, g_2, \ldots, g_M > = G$. Does this imply that a finite cyclic group has a primitive element $g \in G$ which is a product of a few elements, that is, $g = h_1 h_2 \cdots h_m$, where $h_i \in \{g_1, g_2, \ldots, g_M\}$, and $m \leq 6$ is small.

**Problem 12.4.** Determine the least primitive roots in arithmetic progressions $\{\, q\,n + a : n \in \mathbb{Z} \,\}$, $\gcd(a,\ q) = 1$, for small $q \geq 2$.

**Problem 12.5.** Determine the least primitive roots in quadratic arithmetic progressions $\{\, a\,n^2 + b\,n + c : n \in \mathbb{Z} \,\}$, $\gcd(a, b, c) = 1$.

**Problem 12.6.** Let $g_1 < g_2 < \cdots < g_{p-1}$ be the set of primitive roots modulo $p \geq 2$. Verify the average gap is $\overline{g} = \frac{p-1}{\varphi(p-1)} = \prod_{q \mid p-1} (1 - 1/q)^{-1} \ll \prod_{q \leq \log p} (1 - 1/q)^{-1} \ll \log\log p$. What is the distribution primitive roots gaps $d_n = g_{n+1} - g_n$?

**Problem 12.7.** Given a primitive root $g \neq \pm 1$, $r^2$ modulo a prime $p \geq 3$, the subset $\{\, g^n \bmod p : n \geq 1 \,\}$ of integers contains all the primes $q < p$. What is the least exponent $v = v(g) \geq 1$ such that $g^v = q$ is the least prime primitive root?

Trivially, $v = v(g) \leq p - 2$, but is it possible to have $v = v(g) = o(p)$ for any prime $p \geq 2$?

**Problem 12.8.** Let $K = \left\{ k_1, k_2, \ldots, k_{\varphi(p)} : \gcd(k, \varphi(p)) = 1 \right\} \subset G$. Show that $\sum_{n \in K} n \equiv \bmod \varphi(p)$. Hint: study the relation $g^{k_1 + k_2 + \cdots + k_{\varphi(p)}} \equiv ? \bmod p$.

**Problem 12.9.** Let $K = \{ k_1, k_2, \ldots, k_{\lambda(N)} : \gcd(k, \lambda(N)) = 1 \} \subset G$. Show that $\sum_{n \in K} n \equiv \bmod \lambda(N)$.

**Problem 12.10.** Given a large prime $p \geq 2$, and $1 < y \leq x \leq p$, determine the minimal short interval $[x, \, x + y]$ that contains primitive roots, and estimate the cardinality. In summatory function notation, this stated as

$(i) \quad \sum_{x \leq n \leq x+y} \Psi(n) \gg 1.$ $\qquad (ii) \quad \sum_{x \leq n \leq x+y} \Psi(n) = c\, y + O\!\left(x^{\beta}\right), \; c > 0, \; \beta > 0 \text{ constants.}$

| Chapter 13 |
|---|

# Subsets of Integers with Fixed Primitive Roots

Let $u \neq \pm 1, v^2$ be a fixed integer, and let $x \geq 1$ be a large number. Define the subset of integers $\mathcal{N}_u = \{ N \in \mathbb{N} : \gcd(u, N) = 1, \mathrm{ord}(u) = \lambda(N) \} \subset \mathbb{N}$ with a fixed primitive root $u$, and its counting function $N_u(x) = \# \{ N \leq x : \gcd(u, N) = 1, \mathrm{ord}(u) = \lambda(N) \}$. In this more general case, it is easy to demonstrate that an admissible fixed integer $u \neq \pm 1, v^2$ is a primitive root mod $N$ for infinitely many integers $N \geq 2$, $\gcd(u, N) = 1$. For example,

$$u = 2, \text{ and } N = 3^m, m \geq 1; \qquad u = 2, \text{ and } N = 5^m, m \geq 1;$$
$$u = 3, \text{ and } N = 7^m, m \geq 1; \qquad u = 2, \text{ and } N = 11^m, m \geq 1; \qquad (13.1)$$

and so on. Observe that the more restricted case over the primes, known as Artin conjecture for primitive root, is significantly more difficult to prove that the subset of primes $\mathcal{P}_u = \{ p \in \mathbb{N} : \mathrm{ord}(u) = \varphi(p) \}$ is infinite. As the subsets of primes and primes powers are included in the subset of integers $\mathcal{N}_u$ with a fixed primitive root $u$, it expected that $N_u(x) > x / \log x$. The heuristic for the asymptotic formula for the number of integers $N \leq x$ with a fixed primitive root $u$ states that

$$N_u(x) = \frac{\varphi(u)}{u} x \prod_q (1 - F(q, x)), \qquad (13.2)$$

where $F(q, x) = o(1)$ is a complicated function, see [LP03, p. 2]. A related result for the subset of integers $\mathcal{U} = \{ \lambda(n) : n \in \mathbb{N} \} \subset \mathbb{N}$, and its counting function $N'_u(x) = \# \{ \lambda(n) \leq x : n \in \mathbb{N} \}$, is proved here unconditionally. Two proofs are given. The first one in Theorem 13.1, has an approximate asymptotic formula, and the second one in Theorem 13.2 has the correct asymptotic formula, and the constant for the cardinality of the subset of integers

$$\{ \lambda(N) \leq x : \gcd(u, N) = 1, \mathrm{ord}(r) = \lambda(N) \} . \qquad (13.3)$$

## 13.1 First Approach to the Subsets of Integers Fixed Primitive Roots

The analysis takes in a proper subgroup $G \subset \mathbb{Z}_N^*$ of order $\# G = \lambda(N)$ of the multiplicative group $\mathbb{Z}_N^*$ of order $\varphi(N) \geq 1$. If the group index $\varphi(N) / \lambda(N) > 1$, the subgroup $G$ is much smaller than $\mathbb{Z}_N^*$. And if $\varphi(N) / \lambda(N) = 1$, the result seamlessly reduces to the prime case proved in Theorem 12.1 whenever $N \geq 2$ is a prime.

A new technique utilizing the moduli-free characteristic function, see Lemma 3.3, will be employed to illustrate this method.

***Theorem*** 13.1    A fixed integer $u \neq \pm 1, v^2$, and $\gcd(u, N) = 1$, is a primitive root mod $N$ for infinitely many integers

$N \geq 2$. In addition, the density of these integers satisfies

$$\# \{ N \leq x : \gcd(u, N) = 1, \operatorname{ord}(u) = \lambda(N) \} \geq C_2 \, \frac{\varphi(u)}{u} \, \frac{x}{\log\log\log x} + o\left( \frac{x}{\log\log\log x} \right), \tag{13.4}$$

where $C_2 > 0$ is a constant, for all large numbers $x \geq 1$.

**Proof**: Suppose that $u \neq \pm 1, v^2,$ and $\gcd(u, N) = 1$ is not a primitive root for all integers $N \geq x_0$, with $x_0 \geq 1$ constant. Let $p = \lambda(N) + o(\lambda(N))$ be a prime number, and let $x > x_0$ be a large number. Consider the nonexistence equation

$$0 = \sum_{\substack{x \leq N \leq 2\,x \\ \gcd(u,N)=1}} \left( \frac{1}{p} \sum_{\gcd(n,\lambda(N))=1} \frac{1}{n^s} \sum_{0 \leq k \leq \varphi(p-1)} \psi((\theta^n - u)\,k) \right), \tag{13.5}$$

where $\theta \in \mathbb{Z}_N$ is a primitive root mod $N$, and $\psi = e^{i\,2\pi/p} \neq 1$ is an additive character of order $\operatorname{ord} \psi = p$, see Lemma 3.3 for details on the characteristic function. The complex variable $s \in \mathbb{C}$ will be set to $s \geq 2$.

Regroup the triple finite sum as

$$0 = \sum_{\substack{x \leq N \leq 2\,x, \\ ,\gcd(u,N)=1}} \frac{1}{p} \sum_{\gcd(n,\lambda(N))=1} \frac{1}{n^s} + \sum_{\substack{x \leq N \leq 2\,x, \\ \gcd(u,N)=1}} \left( \frac{1}{p} \sum_{\gcd(n,\lambda(N))=1,} \sum_{0 \leq k \leq \varphi(p-1)} \frac{\psi((\theta^n - u)\,k)}{n^s} \right). \tag{13.6}$$

The main term is determined by a finite sum over the trivial additive character $\psi = 1$, and the error term is determined by a finite sum over the nontrivial additive characters $\psi \neq 1$.

Put $s = 2$. Replacing Lemma 8.7 into (13.6) return

$$0 = \sum_{\substack{x \leq N \leq 2\,x, \\ \gcd(u,N)=1}} \frac{1}{p} \sum_{\gcd(n,\lambda(N))=1} \frac{1}{n^s} + \sum_{\substack{x \leq N \leq 2\,x \\ \gcd(u,N)=1}} \frac{1}{p} \left( -1 + O\left( \frac{1}{p^{s-2+\epsilon}} \right) \right)$$

$$> \sum_{\substack{x \leq N \leq 2\,x, \\ \gcd(u,N)=1}} \frac{1}{p} \sum_{\substack{\gcd(n,\lambda(N))=1, \\ n \geq 3}} \frac{1}{n^s} + O\left( \sum_{x \leq N \leq 2\,x} \frac{1}{p^{1+\epsilon}} \right) \tag{13.7}$$

$$\geq \frac{1}{2} \frac{1}{3^2} \sum_{\substack{x \leq N \leq 2x, \\ \gcd(u, N)=1, \\ \lambda(N) \equiv 2,4 \bmod 6}} \frac{1}{p} + O\left( \sum_{x \leq N \leq 2x} \frac{1}{p^{1+\epsilon}} \right)$$

$$> 0.$$

Clearly, this is a contradiction for all sufficiently large numbers $x \geq x_0$. Id est, the short interval $[x, 2x]$ does contain integers $N \geq x$ such that $u \neq \pm 1$, $v^2$ is a primitive root mod $N$.

To compute the asymptotic integers counting formula, take

$$\sum_{\substack{N \leq x, \\ \gcd(u, N)=1}} \frac{1}{p} \sum_{\gcd(n, \lambda(N))=1,} \sum_{0 \leq k \leq \varphi(p-1))} \psi((\theta^n - r)\, k) = \sum_{\substack{N \leq x, \\ \gcd(u, N)=1}} \frac{\varphi(\lambda(N))}{p} + o\left( \sum_{N \leq x} \frac{\varphi(\lambda(N))}{p} \right).$$

(13.8)

This follows from (13.7). Then, by Lemma 7.6, and partial summation, the density of the corresponding subset of primes over the short interval $[x, 2x]$ is

$$\sum_{\substack{N \leq x, \\ \gcd(u, N)=1}} \frac{\varphi(\lambda(N))}{p} + o\left( \sum_{N \leq x} \frac{\varphi(\lambda(N))}{p} \right) = \sum_{\substack{N \leq x, \\ \gcd(u, N)=1}} \frac{\varphi(\lambda(N))}{\lambda(N)} + o\left( \sum_{N \leq x} \frac{\varphi(\lambda(N))}{\lambda(N)} \right) + O(\log x)$$

$$\geq C_2\, \frac{\varphi(u)}{u}\, \frac{x}{\log\log\log x} + o\left( \frac{x}{\log\log\log x} \right).$$

(13.9)

where $C_2 > 0$ is a constant. ∎

The determination of correct order of magnitude is a difficult problem. It involves an oscillating product as shown in equation (13.1). Likewise, the determination of the constant $\alpha \geq 0$ should be a very difficult problem in algebraic number theory. The simpler case for prime numbers took about half a century to settled, see Chapter 11, and [HY67, p. 218], [MO04, p. 10], [SN03].

## 13.2 Second Approach to the Subsets of Integers With Fixed Primitive Roots

Information on the subset of integers $\mathcal{U} = \{ \lambda(n) : n \in \mathbb{N} \}$ leads to a simple proof of the integers counting problem with fixed primitive roots mod $N$. Moreover, it has both the correct order of magnitude, and the constant. The analysis takes place in a maximal multiplicative subgroup $G \subset \mathbb{Z}_N^*$ of order $\# G = \lambda(N)$. For a large number $x \geq 1$, let

$$N_u'(x) = \# \{ \lambda(N) \leq x : \gcd(u, N) = 1,\ \mathrm{ord}(u) = \lambda(N) \} .$$

(13.10)

***Theorem* 13.2**   A fixed integer $u \neq \pm 1, v^2$, and $\gcd(u, N) = 1$, is a primitive root mod $N$ for infinitely many integers $N \geq 2$. In addition, the density of these integers satisfies

$$N'_u(x) = u_0 \, \frac{\varphi(u)}{u} \, \frac{x}{(\log x)^{\eta + o(1)}} + o\left( \frac{\varphi(u)}{u} \, \frac{x}{(\log x)^{\eta + o(1)}} \right), \tag{13.11}$$

where $\eta = 1 - (1 + \log \log 2) / \log 2 = .08607\ldots$, and $u_0 = \sum_{n \in \mathcal{U}} \mu(n) \, n^{-2} > 0$ are absolute constants, for all large numbers $x \geq 1$.

***Proof*:**  Suppose that $u \neq \pm 1, v^2$, with $\gcd(u, N) = 1$, is not a primitive root for all integers $N \geq x_0$, with $x_0 \geq 1$ constant. Let $x > x_0$ be a large number, and consider the nonexistence equation

$$0 = \sum_{\substack{x \leq N \leq 2x, \\ \gcd(u, N) = 1}} \left( \frac{1}{p} \sum_{\gcd(n, \lambda(N)) = 1} \frac{1}{n^s} \sum_{0 \leq k \leq \varphi(p)} \psi((\theta^n - u) \, k) \right), \tag{13.12}$$

where $\theta \in G \subset \mathbb{Z}_N$ is a primitive root mod $N$, $p = p(N) = \lambda(N) + o(\lambda(N))$ is a large prime, and $\psi = e^{i 2\pi / p} \neq 1$ is an additive character of order $\operatorname{ord} \psi = p$, see Lemma 3.3 for details on the characteristic function. The complex variable $s \in \mathbb{C}$ will be set to $s \geq 2$.

Regroup the triple finite sum as

$$0 = \sum_{\substack{x \leq N \leq 2x, \\ \gcd(r, N) = 1}} \frac{1}{p} \sum_{\gcd(n, \lambda(N)) = 1} \frac{1}{n^s} + \sum_{\substack{x \leq N \leq 2x, \\ \gcd(u, N) = 1}} \left( \frac{1}{p} \sum_{\gcd(n, \lambda(N)) = 1,} \frac{1}{n^s} \sum_{0 < k \leq \varphi(p)} \psi((\theta^n - u) \, k) \right). \tag{13.13}$$

The main term is determined by a finite sum over the trivial additive character $\psi = 1$, and the error term is determined by a finite sum over the nontrivial additive characters $\psi \neq 1$.

Put $s = 2$. Applying Lemma 8.5, the trivial upper bound $\varphi(\lambda(N)) / \lambda(N) \leq 1 / 2$, Lemma 5.2, and simplifying return

$$0 = \sum_{\substack{x \leq N \leq 2x, \\ \gcd(u, N) = 1}} \frac{1}{p} \sum_{\gcd(n, \lambda(N)) = 1} \frac{1}{n^s} + \sum_{\substack{x \leq N \leq 2x, \\ \gcd(u, N) = 1}} \frac{1}{p} \left( -\frac{s}{s-1} \, \frac{\varphi(\lambda(N))}{\lambda(N)} + O\left( \frac{1}{p^{s-1}} \right) \right)$$

$$\geq \sum_{\substack{x \leq N \leq 2x, \\ \gcd(u, N) = 1}} \frac{1}{p} \sum_{\gcd(n, \lambda(N)) = 1} \frac{1}{n^2} + \sum_{\substack{x \leq N \leq 2x, \\ \gcd(u, N) = 1}} \left( -\frac{1}{p} + O\left( \frac{1}{p^2} \right) \right) \tag{13.14}$$

$$\geq \sum_{\substack{x \leq N \leq 2\,x, \\ \gcd(u,\,N)=1}} \frac{1}{p} \sum_{\substack{\gcd(n,\,\lambda(N))=1, \\ n \geq 3}} \frac{1}{n^2} + O\!\left( \sum_{x \leq N \leq 2\,x} \frac{1}{p^2} \right).$$

To derive a concrete estimate, it is converted back (from the $p$-domain to the $\lambda(N)$-domain) using $\lambda(N)\,/\,2 < p = \lambda(N) + o(\lambda(N)) < 2\,\lambda(N)$. The precise steps are these:

$$0 \geq \sum_{\substack{x \leq N \leq 2\,x, \\ \gcd(u,\,N)=1}} \frac{1}{p} \sum_{\substack{\gcd(n,\,\lambda(N))=1, \\ n \geq 3}} \frac{1}{n^2} + O\!\left( \sum_{x \leq N \leq 2\,x} \frac{1}{p^2} \right)$$

$$> \frac{1}{2} \sum_{\substack{x \leq N \leq 2\,x, \\ \gcd(u,\,N)=1}} \frac{1}{\lambda(N)} \sum_{\substack{\gcd(n,\,\lambda(N))=1, \\ n \geq 3}} \frac{1}{n^2} + O\!\left( \sum_{x \leq N \leq 2\,x} \frac{1}{\lambda(N)^2} \right)$$

$$\geq \frac{1}{2} \frac{1}{3^2} \sum_{\substack{x \leq N \leq 2\,x, \\ \gcd(u,\,N)=1, \\ \lambda(N) \equiv 2,4 \bmod 6}} \frac{1}{\lambda(N)} + O\!\left( \sum_{x \leq N \leq 2\,x} \frac{1}{\lambda(N)^2} \right) \qquad (13.15)$$

$$> \frac{1}{2} \frac{1}{3} \frac{1}{3^2} \sum_{\substack{x \leq N \leq 2\,x, \\ \gcd(u,\,N)=1}} \frac{1}{N} + O\!\left( \sum_{x \leq N \leq 2\,x} \frac{1}{N^{2-\epsilon}} \right)$$

$$\geq \frac{1}{54} \log(2) + O\!\left( \frac{1}{x^{1-\epsilon}} \right)$$

$$> 0.$$

The third line follows from the fact that $\gcd(3,\,\lambda(N)) = 1$ for any $\lambda(N) \equiv 2,\,4 \bmod 6$. The fourth line follows from the density $1\,/\,3$ of the subset of integers in the arithmetic progression $\{\,\lambda(N) \equiv 2,\,4 \bmod 6 : N \in \mathbb{N}\,\}$, see Theorem 7.7. And the estimates $N^{1-\epsilon} < \lambda(N) < N,$ where $\epsilon > 0$ is a small number.

Clearly, this is a contradiction for all sufficiently large numbers $x \geq x_0$. Id est, the short interval $[x,\,2\,x]$ does contain integers $N \geq x$ such that $u \neq \pm 1$, $v^2$ is a primitive root mod $N$.

To derive the asymptotic integers counting formula, take the sum of the characteristic function over the interval $[1,\,x]$, that is,

$$\sum_{\substack{\lambda(N) \leq x, \\ \gcd(u,\,N)=1}} \frac{1}{p} \sum_{\gcd(n,\,\lambda(N))=1,} \sum_{0 \leq k \leq \varphi(p)} \psi((\theta^n - r)\,k) = \sum_{\substack{\lambda(N) \leq x, \\ \gcd(u,\,N)=1}} \frac{\varphi(\lambda(N))}{p} + o\!\left( \sum_{\lambda(N) \leq x} \frac{\varphi(\lambda(N))}{p} \right).$$

$$(13.16)$$

The right side of (13.16) follows from (13.15), which implies that the main term in (13.16) is larger than the error term. Now observe that since the prime satisfies the asymptotic relation $p \sim \lambda(N)$, the ratio $\varphi(\lambda(N)) \, / \, p$ can be converted back to

$$\frac{\varphi(\lambda(N))}{p} = \frac{\varphi(\lambda(N))}{\lambda(N)} \, \frac{\lambda(N)}{p} = \frac{\varphi(\lambda(N))}{\lambda(N)} + o\!\left( \frac{\varphi(\lambda(N))}{\lambda(N)} \right).$$

$$(13.17)$$

Then, by Theorem 7.10, the density of the corresponding subset of primes over the short interval $[x, \, 2\,x]$ is

$$\sum_{\substack{\lambda(N) \leqslant x, \\ \gcd(u,\,N)=1}} \frac{\varphi(\lambda(N))}{\lambda(N)} + o\!\left( \sum_{\lambda(N) \leqslant x} \frac{\varphi(\lambda(N))}{\lambda(N)} \right) = u_0 \, \frac{\varphi(u)}{u} \, \frac{x}{(\log x)^{\eta + o(1)}} + o\!\left( \frac{\varphi(u)}{u} \, \frac{x}{(\log x)^{\eta + o(1)}} \right),$$

$$(13.18)$$

where $\eta = 1 - (1 + \log\log 2) \, / \log 2 = .08607\ldots$ is a constant, and $u_0 = \sum_{n \,\in\, \mathcal{U}} \mu(n) \, n^{-2} > 0$ is an absolute constant, for all large numbers $x \geqslant 1$. ∎

# Chapter 14

## Other Densities Problems

The calculations of the density of primes with a fixed least primitive root is a refinement of the Artin conjecture. Some works by several authors have been done on this problem.

Let $m \neq \pm 1$, $b^2$ be a fixed integer. The density of primes with a fixed least primitive root $g(p) = m$ is defined by

$$D(m) = \lim_{x \to \infty} \frac{1}{\pi(x)} \sum_{p \leq x,\, g(p)=m} 1 \quad . \tag{14.1}$$

A few formulas and numerical data for $D(m)$, $m \geq 2$, are analyzed and provided in [PB02]. The first few cases are as follows.

$D(2) = \Delta_1$ .

$D(3) = \Delta_1 - \Delta_2$ .

$D(5) = \dfrac{20}{19}\Delta_1 - \dfrac{200}{91}\Delta_2 + \dfrac{500}{439}\Delta_3.$

$D(6) = \Delta_1 - \dfrac{282}{91}\Delta_2 + \dfrac{1000}{439}\Delta_3.$

$D(7) = \Delta_1 - \left(4 + \dfrac{9}{91} + \dfrac{5}{281}\right)\Delta_2 +$

$\qquad \left(6 + \dfrac{183}{439} + \dfrac{4836}{67\,585} + \dfrac{147\,193}{29\,669\,815}\right)\Delta_3 - \left(3 + \dfrac{1107}{2131} + \dfrac{71\,825}{12\,901\,197} + \dfrac{26\,503\,425}{2\,749\,409\,807}\right)\Delta_4.$

$$\tag{14.2}$$

The deltas are given by

$$\Delta_1 = \prod_{p \geq 2}\left(1 - \frac{1}{p(p-1)}\right).$$

$$\Delta_2 = \prod_{p \geq 2}\left(1 - \frac{2}{p(p-1)} + \frac{1}{p^2(p-1)}\right).$$

$$\Delta_3 = \prod_{p \geq 2}\left(1 - \frac{3}{p(p-1)} + \frac{3}{p^2(p-1)} - \frac{1}{p^3(p-1)}\right).$$

$$\tag{14.3}$$

$$\Delta_4 = \prod_{p \geq 2} \left( 1 - \frac{4}{p(p-1)} + \frac{6}{p^2(p-1)} - \frac{4}{p^3(p-1)} + \frac{1}{p^4(p-1)} \right).$$

This idea extends to other subsets of integers, with certain fixed primitive roots.

Let $q \geq 2$ be a fixed prime. The density of primes with a fixed least prime primitive root $g^*(p) = q$ is defined by

$$D_1(q) = \lim_{x \to \infty} \frac{1}{\pi(x)} \sum_{p \leq x, \, g^*(p) = q} 1 \quad .$$
(14.4)

Similarly, let $m \neq \pm 1, \; b^2$ be a fixed integer. The density of integers with a fixed least primitive root $g(N) = m$ is defined by

$$D_2(q) = \lim_{x \to \infty} \frac{1}{x} \sum_{p \leq x, \, g^*(p) = q} 1 \quad .$$
(14.5)

The last two integers counting problems appear be open problems.

**Chapter 15**

# Harmonic Sums

A few finite harmonic sums are stated in this section. These are useful in many applications. The literature on the harmonic sum is extensive, for starter, a wide range of problems are considered in [EP50], and sharp estimates in [VM05].

## 15.1. Harmonic Sums Over The Integers

*Lemma* **15.1.**   Let $x \geq 1$ be a large number, then

$$\sum_{n \leq x} \frac{1}{n} = \log x + \gamma + \frac{1}{2x} + O\left(\frac{1}{x^2}\right). \tag{15.1}$$

The number $\gamma = 0.5772156649 ...$ is Euler constant, see (18.2) for the definition.

## 15.2 Harmonic Sums Over Large Subsets of Integers

The error term for a harmonic sum over a subset of integers of nonzero density, but not 1, is considerably larger than the harmonic sum over the entire set of integers. In some cases there are very sharp estimates but in other cases there are not.

*Lemma* **15.2.**   Let $x \geq 1$ be a large number, and let $\mathcal{A} \subset \mathbb{N}$ be a subset of integers of density $\alpha = \delta(\mathcal{A}) > 0$. Then

$$\sum_{n \leq x, n \in \mathcal{A}} \frac{1}{n} = \alpha \log x + \gamma_\alpha + o(\log x). \tag{15.2}$$

The number $\gamma_\alpha > 0$ is a generalized Euler constant, see (18.6) for the definition.

## 15.3. Harmonic Sums Over Squarefree Integers

*Lemma* **15.3.**   Let $x \geq 1$ be a large number, and let $Q = \{ n \in \mathbb{N} : \mu(n) = \pm 1 \} \subset \mathbb{N}$ be a subset of squarefree integers of density $\alpha = \delta(Q) = 6/\pi^2$. Then

(i) $\displaystyle \sum_{n \leq x, n \in Q} \frac{1}{n} = \frac{6}{\pi^2} \log x + \gamma_\alpha + O\left(\frac{1}{x^{1/2}}\right).$

$(15.3)$

(ii) $\displaystyle\sum_{n \le x, \, n \in Q} \frac{1}{n} = \frac{6}{\pi^2} \log x + \gamma_\alpha + \Omega\left(\frac{1}{x^{3/4}}\right).$

The constant has a simple expression as $\gamma_\alpha = 6 \, \pi^{-2} \, \gamma$.

**Proof**: For (ii), rewrite the finite sum as an integral with respect to the discrete counting measure $Q(x) = 6 \, \pi^{-2} \, x + \Omega\big(x^{1/4}\big)$.
∎

## 15.4 Harmonic Sums For Fixed Primitive Roots

Let $\mathcal{N}_2 = \{\, n \in \mathbb{N} : \mathrm{ord}(2) = \lambda(n) \,\}$ be the subset of integers such that 2 is a primitive root modulo $n \ge 1$, and let $N_2(x) = \#\{\, n \le x : \mathrm{ord}(2) = \lambda(n) \,\}$ be the discrete counting measure, see Theorem 23.1. An asymptotic formula for the harmonic sum over the subset of integers $\mathcal{N}_2$ is determined here.

**Lemma 15.4.** Let $x \ge 1$ be a large number, and let $\mathcal{N}_2 \subset \mathbb{N}$ be a subset of integers of density $\alpha_2 = \delta(\mathcal{N}_2) > 0$. Then

$$\sum_{n \le x, \, n \in \mathcal{N}_2} \frac{1}{n} = \kappa_2 \, (\log x)^{\alpha_2} + \gamma_2 + O\left(\frac{1}{(\log x)^{1-\alpha_2}}\right). \tag{15.4}$$

The number $\alpha_2 = 0.373955813619202\ldots$ is Artin constant, and $\gamma_2 > 0$ is the Artin-Euler constant, see (18.6) for the definition. The other constant is

$$\kappa_2 = \frac{1}{2 \, \alpha_2 \, \Gamma(\alpha_2)} \prod_{\mathrm{ord}(2) = p-1, \, \mathrm{ord}(2) \ne p(p-1)} \left(1 - \frac{1}{p^2}\right) = 0.562432224942\ldots. \tag{15.5}$$

The index of the product ranges over the Wieferich primes, which are also characterized by the congruence $2^{p-1} \equiv 1 \bmod p^2$.

**Proof**: Use the discrete counting measure $N_2(x) = (\alpha_2 \, \kappa_2 + o(1)) \, x(\log x)^{\alpha_2 - 1}$, Theorem 23.1, to write the finite sum as an integral, and evaluate it:

$$\sum_{n \le x, \, n \in \mathcal{N}_2} \frac{1}{n} = \int_{x_0}^{x} \frac{1}{t} \, d \, N_2(t) = \left.\frac{N_2(t)}{t}\right|_{x_0}^{x} + \int_{x_0}^{x} \frac{N_2(t)}{t^2} \, d \, t, \tag{15.6}$$

where $x_0 > 0$ is a constant. Continuing the evaluation yields

$$\int_{x_0}^{x} \frac{1}{t} \, d \, N_2(t) = \frac{(\alpha_2 \, \kappa_2 + o(1))}{(\log x)^{1-\alpha_2}} + c_0(x_0) + \int_{x_0}^{x} \frac{(\alpha_2 \, \kappa_2 + o(1))}{t(\log t)^{1-\alpha_2}} \, d \, t$$

$$= \kappa_2 (\log x)^{\alpha_2} + \gamma_2 + O\left(\frac{1}{(\log x)^{1-\alpha_2}}\right),$$

where $c_0(x_0)$ is a constant. Moreover,

$$\gamma_2 = \lim_{x \to \infty} \left( \sum_{n \leq x, \, n \in \mathcal{N}_2} \frac{1}{n} - \kappa_2 \log^{\alpha_2} x \right) = c_0(x_0) + \int_{x_0}^{\infty} \frac{(\alpha_2 \, \kappa_2 + o(1))}{t(\log t)^{1-\alpha_2}} \, d \, t \qquad (15.8)$$

is a second definition of this constant. ∎

The integration lower limit $x_0 = 2$ appears to be correct one since the subset of integers is $\mathcal{N}_2 = \left\{ 3, \, 5, \, 3^2, \, 11, \, \ldots \right\}$.

## 15.5 Harmonic Sums For The Subsets Of Quadratic Residues And Nonresidues

Let $(r \mid p) = \chi(r)$ denotes the quadratic residue symbol. Let $\mathcal{A}_r = \{ n \in \mathbb{N} : p \mid n \Rightarrow (r \mid p) = 1 \}$ be the subset of integers with a fixed quadratic nonresidue $|r| > 1$.

**Lemma 15.5.** If $x \geq 1$ is a large number, then

$$\sum_{\phantom{x}} \frac{1}{\phantom{x}} = \frac{2}{\phantom{xxx}} \prod (1 - p^{-1})^{-1/2} \prod (1 - p^{-2})^{-1/2} \frac{1}{\phantom{xxxx}} \log(x)^{1/2} \left( 1 + O\left( \frac{1}{\phantom{xxx}} \right) \right). \qquad (15.9)$$

Let $\mathcal{B}_r = \{ n \in \mathbb{N} : p \mid n \Rightarrow (r \mid p) = -1 \}$ be the prime divisors constrained subset of integers.

**Lemma 15.6.** If $x \geq 1$ is a large number, then

$$\sum_{\phantom{x}} \frac{1}{\phantom{x}} = \frac{2}{\phantom{xxx}} \prod (1 - p^{-1})^{1/2} \prod (1 - p^{-2})^{1/2} L(1, \chi)^{1/2} \log(x)^{1/2} \left( 1 + O\left( \frac{1}{\phantom{xxx}} \right) \right). \qquad (15.10)$$

The proof of Lemmas 15.5 and 19.6 are derived from Riegler formula, see Lemma 20.2. Additional details are developed in [WS75].

## 15.6. Fractional Finite Sums

The finite sums of fractional parts provide information about the error terms of various finite sums. The paper [PF10] has a recent survey of many fractional finite sums.

**Lemma 15.7.** (delaVallee Poussin) Let $x \geq 1$ be a large number, and let $\{x\} = x - [x]$ be the fractional part function. Then

(i) $\displaystyle\sum_{n \leq x} \left\{ \frac{x}{n^k} \right\} = (1 - \gamma_{1/k})\, x + O\big(x^{1/(k+1)}\big),$             fixed $k \geq 1$,

(ii) $\displaystyle\sum_{n \leq x,\ p \equiv a \bmod q} \left\{ \frac{x}{n} \right\} = \frac{(1 - \gamma)}{\varphi(q)}\, x + O\big(\sqrt{x}\big),$      small parameters $a < q,\ \gcd(a, q) = 1.$

$\hspace{12cm}(15.11)$

Modulo the Riemann hypothesis, IT expected that the error term in the second finite sum would be $O\big(x^{1/4+\epsilon}\big)$, $\epsilon > 0$. This is essentially equivalent to the circle problem, and the divisor problem.

## 15.7 Problems And Exercises

**Problem 15.1.** Let $x \geq 1$ be a large number, and let $\mathcal{N}_u \subset \mathbb{N}$ be a subset of integers with a fixed primitive root $u \bmod n$ of density $\alpha_u = \delta(\mathcal{N}_u) > 0$, see Theorem 23.2. Prove that

$$\sum_{n \leq x,\, n \in \mathcal{N}_u} \frac{1}{n} = \kappa_u \log^{\alpha_u} x + \gamma_u + O\left( \frac{1}{(\log x)^{1-\alpha_u}} \right).$$
$\hspace{12cm}(15.12)$

The number $\alpha_u > 0$ is Artin constant, $\gamma_u > 0$ is the Artin-Euler constant, see Definition 18.2, and the other constant is

$$\kappa_u = \frac{\varphi(u)}{u\, \alpha_u\, \Gamma(\alpha_u)} \prod_{\text{ord}(u)=p-1,\, \text{ord}(u) \neq p(p-1)} \left( 1 - \frac{1}{p^2} \right).$$
$\hspace{12cm}(15.13)$

Here the product ranges over the of the Abel primes, which are also characterized by the congruence $u^{p-1} \equiv 1 \bmod p^2$.

**Problem 15.2.** Let $x \geq 1$ be a large number, and let $\mathcal{N}_u \subset \mathbb{N}$ be a subset of integers with a fixed primitive root $u \bmod n$ of density $\alpha_u = \delta(\mathcal{N}_u) > 0$, see Theorem 23.2. For a large $z \leq x$, calculate the average constant

$$\frac{1}{z \log x} \sum_{u \leq z,\, u \neq 1,\, v^2} \left( \sum_{n \leq x,\, n \in \mathcal{N}_u} \frac{1}{n} \right) \overset{?}{=} \frac{1}{\alpha_u\, \Gamma(\alpha_u)} \prod_{\text{ord}(u)=p-1,\, \text{ord}(u) \neq p(p-1)} \left( 1 - \frac{1}{p^2} \right).$$
$\hspace{12cm}(15.14)$

**Problem 15.3.** Use Delange theorem, see Lemma 20.3, to reprove Lemma 15.4.

**Problem 15.4.** Let $x \geq 1$ be a large number, and let $\mathcal{N}_\varphi = \{ n \in \mathbb{N} : p \mid n \Rightarrow \mu(p-1) = \pm 1 \}$ be a subset of integers divisible by squarefree totients $p \geq 2$ such that $\mu(p-1) = \pm 1$. Find an asymptotic formula for $N_\varphi = \#\{ n \leq x : p \mid n \Rightarrow \mu(p-1) = \pm 1 \}$. Hint: Apply Wirsing formula, Lemma 20.1 or something similar.

**Problem 15.5.** Find a collection of functions $f : \mathbb{Z} \longrightarrow \mathbb{C}$, and constants $c_f > 0$ and $d_f > 0$ such that

(i) $\displaystyle\sum_{n \leq x} f(n) \left\{ \frac{x}{n} \right\} \sim c_f \sum_{n \leq x} f(n),$     (ii) $\displaystyle\sum_{n \leq x} f(n)\, \Lambda(n) \left\{ \frac{x}{n} \right\} \sim d_f \sum_{n \leq x} f(n)\, \Lambda(n),$     see [PL10, p. 82].    (15.15)

**Problem 15.6.** Let $x \geq 1$ be a large number, and let $\mathcal{N}_u \subset \mathbb{N}$ be a subset of integers with a fixed primitive root $u \bmod n$ of density $\alpha_u = \delta(\mathcal{N}_u) > 0$, see Theorem 23.2. For a large $z \leq x$, calculate the average constant.

**Problem 15.7.** Let $a, q \in \mathbb{N}$ be integers, $1 \leq a < q$, $\gcd(a, q) = 1$. Generalize the digamma function to arithmetic progressions, that is,

$$\psi(n+1) = -\gamma + \sum_{n \leq x} \frac{1}{n} \quad \text{to} \quad \psi(n+1, q, a) = -\gamma(q, a) + \sum_{n \leq x,\, n \equiv a \bmod q} \frac{1}{n}. \tag{15.16}$$

## Chapter 16

# Prime Numbers Theorems

The primes counting function is defined by the usual symbol $\pi(x) = \#\{\, p \leq x : p \text{ is prime}\,\}$, and for an arithmetic progression $\{\, n \equiv a \bmod q : \gcd(q, a) = 1 \,\}$, it is $\pi(x, q, a) = \#\{\, p \leq x : p \equiv a \bmod q \text{ is prime}\,\}$. Extensive details on the prime number theorems are provided in [EL85], [IK04], [MV07], [NW00], [TM95], et alii.

## 16.1 Asymptotic Formulas

The logarithm integral $\operatorname{li}(x)$ arises as the main term in the prime counting function $\pi(x)$. This improper integral can be estimated using a recursive evaluation.

**Lemma 16.1.** Let $x \geq 1$ be a large number. Then

$$\text{(i)} \quad \operatorname{li}(x) = \int_c^x \frac{1}{\log t}\, d\,t = \frac{x}{\log x} + \int_c^x \frac{1}{\log^2 t}\, d\,t,$$

$$\text{(ii)} \quad \operatorname{li}(x) = x \sum_{1 \leq n \leq m} \frac{(n-1)!}{\log(x)^n} + O\!\left(\frac{x}{\log(x)^m}\right),$$

$$\text{(iii)} \quad \operatorname{li}_n(x) = \int_c^x \frac{1}{\log^n t}\, d\,t = \frac{x}{\log^n x} + \int_c^x \frac{1}{\log^{n+1} t}\, d\,t,$$

(16.1)

where $c > 0$ is a constant. More compactly, statements (i) and (iii) can be written as $\operatorname{li}_n(x) = x \log^{-n}(x) + \operatorname{li}_{n+1}(x)$, $n \geq 1$.

Often, the constant is set to $c > 1$, away from the singularity of the integrand at $t = 1$. Other details are given in [CO07, p. 273], [MV07, p. 188], [DF12, p. 7].

**Theorem 16.2.** (delaVallee Poussin, Hadamard)   Let $x \geq 1$ be a large number. Then, there is an absolute constant $c > 0$ such that

$$\text{(i)} \quad \pi(x) = \operatorname{li}(x) + O\!\left(x\, e^{-c\sqrt{\log x}}\right), \qquad\qquad \text{Unconditionally,}$$

$$\text{(ii)} \quad \pi(x) = \operatorname{li}(x) + O\!\left(x^{1/2} \log^2 x\right), \qquad\qquad \text{Conditional on the RH.}$$

(16.2)

**Theorem 16.3.** (Siegel-Walfisz) Let $x \geq 1$ be a large number. If $a < q$ are relatively prime, and $q = O(\log^B x)$, $B \geq 0$ constant, then, there is an absolute constant $c > 0$ such that

(i) $\pi(x, q, a) = \dfrac{1}{\varphi(q)} \operatorname{li}(x) + O\left(x \, e^{-c \sqrt{\log x}}\right),$          Unconditionally,

                                                                           (16.3)

(ii) $\pi(x, q, a) = \dfrac{1}{\varphi(q)} \operatorname{li}(x) + O\left(x^{1/2} \log^2 x\right),$         Conditional on the GRH.

**Theorem 16.4.** (Brun–Titchmarsh) Let $x \geq 1$ be a large number. If $a < q$ are relatively prime, then,

(i) $\pi(x, a, q) \leq \dfrac{2 \, x}{\varphi(q) \log(x \, / \, q)},$

                                                                           (16.4)

(ii) $\pi(x + h, a, q) - \pi(x, a, q) \leq \dfrac{2 \, h}{\varphi(q) \log(h \, / \, q)} \, (1 + O(1 \, / \log(h \, / \, q)),$

for all $1 \leq q < x$, and $2 \leq 2 \, q \leq h$.

## 16.2 Densities Of Primes In irrational Sequences

There are many subsets of primes specified by various irrationality constraints. Two of these are presented here. Given a pair of irrational numbers $\alpha, \beta \in \mathbb{R} - \mathbb{Q}$, these are written as $\pi_\alpha = \{p \in \mathbb{P} : p = [\alpha \, n], \, n \in \mathbb{N}\}$, and $\pi_\beta = \{p \in \mathbb{P} : p = [n^\beta], \, n \in \mathbb{N}\}$, and the corresponding counting functions as $\pi_\alpha(x) = \#\{p \leq x : p = [\alpha \, n], \, n \in \mathbb{N}\}$, and $\pi_\beta(x) = \#\{p \leq x : \mathbb{P} : p = [n^\beta], \, n \in \mathbb{N}\}$ respectively.

**Theorem 16.5.** (Vinogradov) Let $x \geq 1$ be a large number, and let $\alpha \in \mathbb{R} - \mathbb{Q}$ be an irrational number. Then

$$\pi_\alpha(x) = \frac{1}{\alpha} \operatorname{li}(x) + O\left(x \, e^{-c \sqrt{\log x}}\right). \tag{16.5}$$

**Theorem 16.6.** (Piateski-Shapiro) Let $x \geq 1$ be a large number, and let $\alpha \in \mathbb{R} - \mathbb{Q}$ be an irrational number, $1 < \alpha < 12 \, / \, 11$. Then

$$\pi_\alpha(x) = \frac{1}{\alpha} \operatorname{li}(x) + O\left(x \, e^{-c \sqrt{\log x}}\right). \tag{16.6}$$

**Chapter 17**

# Prime Harmonic Sums Over Primes

The basic techniques of harmonic sums over the prime numbers are illustrated in the next few sections. These results are useful in the study of subset of primes with various properties or characteristics. An algorithm for computing large value is given in [BK09].

## 17.1 Prime Harmonic Sums Over The Primes

The prime harmonic sums associated with the set of primes $\mathbb{P} = \{2, 3, 5, \ldots\}$ have well known estimates.

***Lemma* 17.1.** (Mertens)  Let $x \geq 1$ be a large number. Then, there exists a pair of constants $\beta_1$ and $\gamma$ such that

(i) $\displaystyle\sum_{p \leq x} \frac{1}{p} = \log\log x + \beta_1 + O\!\left(e^{-c\sqrt{\log x}}\right).$

$\hspace{10cm}$ (17.1)

(ii) $\displaystyle\sum_{p \leq x} \frac{\log p}{p - 1} = \log x - \gamma + O\!\left(e^{-c\sqrt{\log x}}\right).$

For large $x \geq 1$, the second finite sum can be rewritten in term of the vonMangoldt function as, see [MV07, p. 182]:

$$\sum_{p \leq x} \frac{\log p}{p - 1} = \sum_{n \leq x} \frac{\Lambda(n)}{n} = \log x - \gamma + O\!\left(e^{-c\sqrt{\log x}}\right).$$

$\hspace{10cm}$ (17.2)

These classical results generalize in many different ways. One of the possible generalizations to arbitrary subsets $\mathcal{A} \subset \mathbb{P}$ of primes of nonzero densities $\alpha = \delta(\mathcal{A}) > 0$ has the following structures.

## 17.2 Prime Harmonic Sums Over Large Subsets Of Primes

***Lemma* 17.2.**  Let $\alpha = \delta(\mathcal{A}) > 0$ be the density of a subset of primes $\mathcal{A} \subset \mathbb{P}$, and let $x \geq 1$ be a large number. Then, there exists a pair of constants $\beta_\alpha = \alpha\beta_1$ and $\gamma_\alpha = \alpha\gamma$ such that

(i) $\displaystyle\sum_{p \leq x,\ p \in \mathcal{A}} \frac{1}{p} = \alpha\log\log x + \beta_\alpha + O\!\left(\frac{1}{\log^2 x}\right).$

$\hspace{10cm}$ (17.3)

(ii) $\displaystyle\sum_{p \leq x, \, p \in \mathcal{A}} \frac{\log p}{p-1} = \alpha \log x - \gamma_\alpha + O\left(\frac{1}{(\log x)^2}\right).$

**Proof:** Since the density of the subset of primes $\mathcal{A} \subset \mathbb{P}$ is $\lim_{x \to \infty} \pi_{\mathcal{A}}(x) / \pi(x) = \alpha > 0$, the corresponding primes counting function has the form

$$\pi_{\mathcal{A}}(x) = \#\{\, p \leq x : p \in \mathcal{A} \,\} = \alpha \operatorname{li}(x) + o(\operatorname{li}(x)) \,, \tag{17.4}$$

where $\operatorname{li}(x)$ is the logarithm integral, and $\alpha > 0$ is the density, for all large numbers $x \geq 1$. The corresponding Stieltjes integral representation is

$$\sum_{p \leq x, \, p \in \mathcal{A}} \frac{1}{p} = \int_a^x \frac{1}{t} \, d\, \pi_{\mathcal{A}}(t) = \alpha \int_a^x \frac{1}{t} \, d\, \pi(t) \,, \tag{17.5}$$

where $a > 0$ is a constant. Use the recursive expression for the logarithm integral, see Lemma 16.1, to complete the evaluation as

$$
\begin{aligned}
\alpha \int_a^x \frac{1}{t} \, d\, \pi(t) &= \frac{\alpha \pi(x)}{x} + \alpha \int_a^x \frac{\pi(t)}{t^2} \, d\, t \\[2mm]
&= \alpha \left( \frac{1}{\log x} + \frac{1}{x}\left(\operatorname{li}_2(x) + O\left(x \, e^{-c\sqrt{\log x}}\right)\right) \right) \\[2mm]
&\qquad + \alpha \int_a^x \left( \frac{1}{t \log t} + \frac{1}{t^2}\left(\operatorname{li}_2(t) + O\left(x \, e^{-c\sqrt{\log x}}\right)\right) \right) d\, t \\[2mm]
&= \alpha \log\log x + \beta_\alpha + O\left(\frac{1}{\log^2 x}\right),
\end{aligned}
$$

$$\tag{17.6}$$

where the constant $\beta_\alpha = \alpha \beta_1$, and $\beta_1 = 0.261497212847\ldots$ is Mertens constant, and $c > 0$ is an absolute constant, see Theorem 16.2. This proves the claim (i). The partial summation, id est, evaluate the integral

$$\int_c^x (\log t) \, d\, R(t) \,, \tag{17.7}$$

where $R(t) = \alpha \log\log t + \beta_\alpha + O\left(\log^{-2} t\right)$, and the constant $\gamma_\alpha = \alpha \gamma$. $\blacksquare$

## 17.3 Prime Harmonic Sums Over Subsets Of Arithmetic Progressions

Let $\mathcal{A} \subset \mathbb{P}$ be a subset of primes of nonzero density $\alpha = \delta(\mathcal{A}) > 0$. The subsets of primes numbers $\mathcal{A}(q, a) = \{\, p \in \mathcal{A} : p \equiv a \bmod q \,\}$ contained in the arithmetic progressions $\{\, p \in \mathbb{P} : p \equiv a \bmod q \,\}$, are investigated in this section.

**Lemma** 17.3. Let $x \geq 1$ be a large number, and let $a < q$ are relatively prime numbers, with $q = O(\log^B x)$, $B \geq 0$ constant. If $\alpha / \varphi(q) = \delta(\mathcal{A}(q, a)) > 0$ is the density of a subset of primes $\mathcal{A}(q, a) \subset \mathbb{P}$, and $x \geq 1$ is a large number, then, there exists a pair of constants $\beta_\alpha(q, a) > 0$ and $\gamma_\alpha(q, a) > 0$ such that

(i) $\displaystyle\sum_{p \leq x,\ p \in \mathcal{A}(q, a)} \frac{1}{p} = \frac{\alpha}{\varphi(q)} \log\log x + \beta_\alpha(q, a) + O\left(\frac{1}{\log^2 x}\right).$

(ii) $\displaystyle\sum_{p \leq x,\ p \in \mathcal{A}(q, a)} \frac{\log p}{p - 1} = \frac{\alpha}{\varphi(q)} \log x - \gamma_\alpha(q, a) + O\left(\frac{1}{\log^2 x}\right).$

$(17.8)$

A table of the values of the constants for small parameters $a, q$ was compiled in [LZ09].

## 17.3 Prime Harmonic Sums Over Primes With Fixed Quadratic Nonresidues

The subset of primes $\mathcal{P}_{2,r} = \{\, p \in \mathbb{P} : \operatorname{ord}(r) = 2 \,\} \subset \mathbb{P}$ consists of all the primes with a fixed quadratic nonresidue $r \in \mathbb{Z}$, and its density is $\kappa_{2,r} = \delta(\mathcal{P}_{2,r})$.

**Lemma** 17.4. Let $\kappa_{2,r} = \delta(\mathcal{P}_{2,r})$ be the density of a subset of primes a fixed quadratic nonresidue $r \in \mathbb{Z}$, and let $x \geq 1$ be a large number. Then, there exists a pair of constants $\beta_{2,r} > 0$ and $\gamma_{2,r} > 0$ such that

(i) $\displaystyle\sum_{p \leq x,\ p \in \mathcal{P}_{2,r}} \frac{1}{p} = \kappa_{2,r} \log\log x + \beta_{2,r} + O\left(\frac{1}{\log^2 x}\right).$

(ii) $\displaystyle\sum_{p \leq x,\ p \in \mathcal{P}_{2,r}} \frac{\log p}{p - 1} = \kappa_{2,r} \log x - \gamma_{2,r} + O\left(\frac{1}{\log^2 x}\right).$

$(17.9)$

**Example** 17.3.1. The subset of primes $\mathcal{P}_{2,2} = \{\, p \in \mathbb{P} : \operatorname{ord}(2) = 2 \,\} = \{\, p = 8\,m \pm 3 \in \mathbb{P} : m \geq 1 \}$ consists of all the primes with a fixed quadratic nonresidue $r = 2$. The density is $1 / 2 = \delta(\mathcal{P}_{2,2})$. This follows from the quadratic reciprocity law

$$\left(\frac{2}{p}\right) = (-1)^{(p^2 - 1)/8}.$$

$(17.10)$

## 17.4 Prime Harmonic Sums Over Primes With Fixed Cubic Nonresidues

The subset of primes $\mathcal{P}_{3,r} = \{\, p \in \mathbb{P} : \mathrm{ord}(r) = 3 \,\} \subset \mathbb{P}$ consists of all the primes with a fixed cubic nonresidue $r \in \mathbb{Z}$, and its density is $\kappa_{3,r} = \delta(\mathcal{P}_{2,r})$.

**Lemma 17.5.** Let $\kappa_{3,r} = \delta(\mathcal{P}_{3,r})$ be the density of a subset of primes a fixed cubic nonresidue $r \in \mathbb{Z}$, and let $x \geq 1$ be a large number. Then, there exists a pair of constants $\beta_{3,r} > 0$ and $\gamma_{3,r} > 0$ such that

(i) $\displaystyle\sum_{p \leq x,\ p \in \mathcal{P}_{3,r}} \frac{1}{p} = \kappa_{3,r} \log\log x + \beta_{3,r} + O\left(\frac{1}{\log^2 x}\right).$

(ii) $\displaystyle\sum_{p \leq x,\ p \in \mathcal{P}_{3,r}} \frac{\log p}{p-1} = \kappa_{3,r} \log x - \gamma_{3,r} + O\left(\frac{1}{\log^2 x}\right).$

$$(17.11)$$

**Example 17.4.1.** The subset of primes $\mathcal{P}_{3,2} = \{\, p \in \mathbb{P} : \mathrm{ord}(2) = 2 \,\} = \{\, p = 6m + 1 \in \mathbb{P} : m \geq 1 \,\}$ consists of all the primes with a fixed cubic nonresidue $r = 2$. The density is $1 \,/\, 2 = \delta(\mathcal{P}_{3,2})$.

## 17.5 Prime Harmonic Sums Over Primes With Fixed Quartic Nonresidues

The subset of primes $\mathcal{P}_{4,r} = \{\, p \in \mathbb{P} : \mathrm{ord}(r) = 4 \,\} \subset \mathbb{P}$ consists of all the primes with a fixed quadratic nonresidue $r \in \mathbb{Z}$, and its density is $\kappa_{4,r} = \delta(\mathcal{P}_{4,r})$.

**Lemma 17.6.** Let $\kappa_{4,r} = \delta(\mathcal{P}_{4,r})$ be the density of a subset of primes a fixed quadratic nonresidue $r \in \mathbb{Z}$, and let $x \geq 1$ be a large number. Then, there exists a pair of constants $\beta_{4,r} > 0$ and $\gamma_{4,r} > 0$ such that

(i) $\displaystyle\sum_{p \leq x,\ p \in \mathcal{P}_{4,r}} \frac{1}{p} = \kappa_{4,r} \log\log x + \beta_{4,r} + O\left(\frac{1}{\log^2 x}\right).$

(ii) $\displaystyle\sum_{p \leq x,\ p \in \mathcal{P}_{4,r}} \frac{\log p}{p-1} = \kappa_{4,r} \log x - \gamma_{4,r} + O\left(\frac{1}{\log^2 x}\right).$

$$(17.12)$$

## 17.6 Prime Harmonic Sums Over Primes With Fixed Primitive Roots

The subset of primes $\mathcal{P}_u = \{\, p \in \mathbb{P} : \mathrm{ord}(u) = p - 1 \,\} \subset \mathbb{P}$ consists of all the primes with a fixed primitive root $u \in \mathbb{Z}$, and its density is $\alpha_u = \delta(\mathcal{P}_u) > 0$.

**Lemma 17.7.** Let $\alpha_u = \delta(\mathcal{P}_u) > 0$ be the density of a subset of primes, and let $x \geq 1$ be a large number. Then, there exists a pair of constants $\alpha_u > 0$, and $\gamma_u > 0$ such that

(i) $\displaystyle\sum_{p \leq x,\ p \in \mathcal{P}_u} \frac{1}{p} = \alpha_u \log\log x + \beta_u + O\left(\frac{1}{\log^2 x}\right).$

(17.13)

(ii) $\displaystyle\sum_{p \leq x,\ p \in \mathcal{P}_u} \frac{\log p}{p-1} = \alpha_u \log x - \gamma_u + O\left(\frac{1}{\log^2 x}\right).$

The constants are $\beta_u = \beta_1\, \alpha_u$ and $\gamma_u = \gamma \alpha_u$.

## 17.7 Prime Harmonic Sums Over Squarefree Totients

The subset of squarefree prime totients is denoted by $\mathcal{T} = \{\, p \in \mathbb{P} : \mu(p-1) = \pm 1 \,\}$. It has nonzero density $\tau = \delta(\mathcal{T}) = \alpha_2$ in the set of primes $\mathbb{P}$.

**Lemma 17.8.** Let $\tau = \delta(\mathcal{T}) > 0$ be the density of a subset of squarefree prime totients $\mathcal{T}$, and let $x \geq 1$ be a large number. Then, there exists a pair of constants $\beta_\tau$ and $\gamma_\tau$ such that

(i) $\displaystyle\sum_{p \leq x,\ p \in \mathcal{T}} \frac{1}{p} = \tau \log\log x + \beta_\tau + O\left(\frac{1}{\log^2 x}\right).$

(17.14)

(ii) $\displaystyle\sum_{p \leq x,\ p \in \mathcal{T}} \frac{\log p}{p} = \tau \log x - \gamma_\tau + O\left(\frac{1}{\log^2 x}\right).$

The constants for this case are $\beta_\tau = \alpha_2\, \beta_1$ and $\gamma_u = \alpha_2\, \gamma$.

## 17.8. Fractional Finite Sums Over The Primes

The finite sums of fractional parts provide information about the error terms of various finite sums over the primes. The paper [PF10] has a recent survey of many fractional finite sums.

**Lemma 17.9.** (delaVallee Possin) Let $x \geq 1$ be a large number, and let $\{\, x \,\} = x - [\, x \,]$ be the fractional part function. Then

(i) $\displaystyle\sum_{p \leq x} \left\{ \frac{x}{p^k} \right\} = (1 - \gamma_{1/k})\, \frac{x}{\log x} + O\left( x\, e^{-c\sqrt{\log x}} \right),$ \qquad fixed $k \geq 1$,

(17.15)

(ii) $\displaystyle\sum_{p \leq x,\ p \equiv a \bmod q} \left\{ \frac{x}{p} \right\} = \frac{(1-\gamma)}{\varphi(q)}\, \frac{x}{\log x} + O\left( x\, e^{-c\sqrt{\log x}} \right),$ \qquad small parameters $a < q$, $\gcd(a, q) = 1$, where $\gamma > 0$ is Euler constant, and $c > 0$ is an absolute constant.

## 17.9 Problems And Exercises

**Problem 17.1.** Determine the finite sums over the primes, and the associated constants:

(i) $\displaystyle\sum_{p \leq x} \frac{1}{p \log p}$. (ii) $\displaystyle\sum_{p \leq x} \frac{1}{p \log\log p}$ .

$$(17.16)$$

**Problem 17.2.** Let $(r \mid p)$ denotes the quadratic residue symbol, and let $r \neq \pm 1$, $s^2$ be a fixed integer. Verify the finite sums over the primes, and the associated constant:

$$\sum_{p \leq x,\, (r|p)=1} \frac{1}{p} = \frac{1}{2} \log\log x + \tau_+ + O\!\left(\frac{1}{\log x}\right).$$

$$(17.17)$$

**Problem 17.3.** Let $(r \mid p)$ denotes the quadratic residue symbol, and let $r \neq \pm 1$, $s^2$ be a fixed integer. Verify the finite sums over the primes, and the associated constant:

$$\sum_{p \leq x,\, (r|p)=-1} \frac{1}{p} = \frac{1}{2} \log\log x + \tau_- + O\!\left(\frac{1}{\log x}\right).$$

$$(17.18)$$

**Problem 17.4.** Compute approximate quantities of the constants:

$$\tau_- = \lim_{x \to \infty} \left( \sum_{p \leq x,\, (r|p)=-1} \frac{1}{p} - \frac{1}{2} \log\log x \right) \text{ and } \tau_+ = \lim_{x \to \infty} \left( \sum_{p \leq x,\, (r|p)=1} \frac{1}{p} - \frac{1}{2} \log\log x \right).$$

$$(17.19)$$

**Chapter 18**

# Connecting Constants

Let $\mathcal{A} \subset \mathbb{P}$ be a subset of prime numbers, and let $\mathcal{B}(\mathcal{A}) = \{\, n \in \mathbb{N} : p \mid n \Rightarrow p \in \mathcal{A} \,\}$ be the multiplicative subset of integers generated buy the primes in $\mathcal{A}$. The subset arises naturally in the $L$-series

$$L(s, \mathcal{A}) = \prod_{p \in \mathcal{A}} \left(1 - \frac{1}{p^s}\right)^{-1} = \sum_{n \in \mathcal{B}} \frac{1}{n^s} , \qquad (18.1)$$

where $s \in \mathbb{C}$ is a complex number.

## 18.1 Definitions

The *connecting constants* link the additive and multiplicative structures of an arbitrary subset of prime numbers $\mathcal{A}$ and the multiplicative subset of integers $\mathcal{B}(\mathcal{A})$ generated by the primes in $\mathcal{A}$. This idea generalizes the classical connection between the Euler constant of the harmonic sum, which is defined by the limit

$$\gamma = \lim_{x \to \infty} \left( \sum_{n \leq x} \frac{1}{n} - \log x \right), \qquad (18.2)$$

and the Mertens constant of the harmonic sum, which is defined by the limit

$$B_1 = \lim_{x \to \infty} \left( \sum_{p \leq x} \frac{1}{p} - \log\log x \right). \qquad (18.3)$$

**Lemma 18.1.** ([HW08, Theorem 427]) Let $\mathcal{A} = \mathbb{P}$ be the set of primes of density $\alpha = \delta(\mathcal{A}) = 1$, and let $\mathcal{B}(\mathcal{A}) = \{\, n \in \mathbb{N} : p \mid n \Rightarrow p \in \mathcal{A} \,\} = \mathbb{N}$ be the set of integers generated by $\mathcal{A}$. Then

(i) $B_1 = \gamma - \sum_{p \geq 2} \left( \log\left(1 - \frac{1}{p}\right) + \frac{1}{p} \right) = \gamma - \sum_{p \geq 2} \sum_{k \geq 2} \frac{1}{k \, p^k} ,$

$(18.4)$

(ii) $B_1 = \gamma - \sum_{p \leq x} \sum_{k \geq 2} \frac{1}{k\, p^k} + O\left(\frac{1}{x}\right),$

for any large number $x \geq 1$.

Some of the associated constants are

$$B_2 = -\sum_{p \geq 2} \left( \log\left(1 - \frac{1}{p}\right) + \frac{1}{p} \right) = \sum_{p \geq 2} \sum_{k \geq 2} \frac{1}{k\, p^k} = 0.31559103872859671283668 9 \ldots \, .$$

$$B_3 = \sum_{p \geq 2} \left( \log\left(1 + \frac{1}{p}\right) - \frac{1}{p} \right) = \sum_{p \geq 2} \sum_{k \geq 2} \frac{(-1)^{k+1}}{k\, p^k} = 0.18185460401679502838197 8 \ldots \, .$$

(18.5)

There is a vast literature on the arithmetic properties of the constants $\gamma$ and $\beta_1$. The statement of Lemma 18.1 is essentially a linear independence result, see [UW12, p. 34]. The generalizations of the Euler constant, in several directions, are discussed in [LR75], [DF08], [DR92], and many other references. The Euler constants has a labyrinths of equivalent definitions, refer to [LJ13]. In addition to (18.2) a few other definitions are

(ii) $\gamma = \lim_{s \to 1} \left( \frac{\zeta'(s)}{\zeta(s)} + \frac{1}{s-1} \right),$ \qquad this arises from \quad $\frac{\zeta'(s)}{\zeta(s)} = -\frac{1}{s-1} + \gamma + \sum_{n \geq 1} \frac{\gamma_n}{n!}(s-1)^n.$

(iii) $\gamma = -\lim_{s \to \infty} \left( \sum_{n \leq x} \frac{\Lambda(n)}{n} - \log x \right),$ \qquad this arises from \quad $\sum_{n \leq x} \frac{\Lambda(n)}{n} = \log x - \gamma + O\left(\frac{1}{\log x}\right).$

(18.6)

A uniform generalization, suitable for the theory of primitive root, and in terms of the subsets of primes of interest will be used here.

**Definition** 18.2. Let $\mathcal{A} \subset \mathbb{P}$ be an arbitrary subset of primes of density $\alpha = \delta(\mathcal{A}) > 0$. The *connecting constants* are defined by the Kronecker-Euler constant

$$\gamma_\alpha = \lim_{x \to \infty} \left( \sum_{p \leq x,\; p \in \mathcal{A}} \frac{\log p}{p-1} - \alpha \log x \right),$$

(18.7)

and the Kronecker-Mertens constant

$$\beta_\alpha = \lim_{x \to \infty} \left( \sum_{p \leq x,\, p \in \mathcal{A}} \frac{1}{p} - \alpha \log\log x \right).$$

(18.8)

Some advance theory of the Kronecker-Euler constants appears in [LM14], and the references within. The special case $\alpha = \alpha_u > 0$ is Artin constant, the corresponding constant $\gamma_u > 0$ is referred to as the Artin-Euler constant, which is a more appropriate nomenclature.

***Lemma* 18.3.** Let $\mathcal{A} \subset \mathbb{P}$ be an arbitrary subset of primes of density $\alpha = \delta(\mathcal{A}) > 0$. If $x \geq 1$ is a large number, then the connecting constants satisfy the equations

(i) $\quad \beta_\alpha = \gamma_\alpha - \sum\limits_{p \in \mathcal{A},\, k \geq 2} \dfrac{1}{k\, p^k}$ .

(18.9)

(ii) $\quad \beta_\alpha = \gamma_\alpha - \sum\limits_{p \leq x,\, p \in \mathcal{A},\, k \geq 2} \dfrac{1}{k\, p^k} + O\!\left(\dfrac{1}{x}\right).$

**Proof**: (ii) From the definition of the constant

$$\beta_\alpha = \lim_{x \to \infty} \left( \sum_{p \leq x,\, p \in \mathcal{A}} \frac{1}{p} - \alpha \log\log x \right) = \gamma_\alpha - \lim_{x \to \infty} \left( \sum_{p \leq x,\, p \in \mathcal{A},\, k \geq 2} \frac{1}{k\, p^k} \right). \tag{18.10}$$

Express the finite sum as a difference of infinite sums, and estimate the error:

$$\sum_{p \leq x,\, p \in \mathcal{A},\, k \geq 1} \frac{1}{k\, p^k} = \sum_{p \in \mathcal{A},\, k \geq 2} \frac{1}{k\, p^k} - \sum_{p > x,\, p \in \mathcal{A},\, k \geq 2} \frac{1}{k\, p^k} = \sum_{p \in \mathcal{A},\, k \geq 2} \frac{1}{k\, p^k} + O\!\left(\frac{1}{x}\right). \tag{18.11}$$

This proves the assertion. ∎

## 18.2 Problems And Exercises

**Problem 18.1.** Explain why the "definition" of Euler constant

$$\gamma_1 = -\lim_{x \to \infty} \left( \sum_{p \leq x} \frac{\log p}{p} - \log x \right) \tag{18.12}$$

is incorrect. For example, at $x = 500\,000$, it has the value

$$\log x - \sum_{p \leq x} \frac{\log p}{p} = 1.3308718305277886045 8134\ldots. \tag{18.13}$$

which is inaccurate. The actual value should be $\gamma = 0.5772156649015328606065 12\ldots$. Hint: for $x = 500\,000$, consider

$$\sum_{n \leq x} \frac{\Lambda(n)}{n} = \sum_{p \leq x} \frac{\log p}{p-1} = 0.5755072182908552112823 26\ldots. \tag{18.14}$$

**Problem 18.2.** Explain the convergence property of the definition of Mertens constant

$$\beta_1 \ = \ \lim_{x \to \infty} \left( \sum_{p \leq x} \frac{1}{p} - \log\log x \right) . \tag{18.15}$$

For example, at $x = 500\,000$, it has the value

$$\sum_{p \leq x} \frac{1}{p} - \log\log x = 0.261613989133765329377640 \ldots . \tag{18.16}$$

**Problem 18.3.** Use the product to compute Artin constant:

$$\alpha_1 \ = \ \lim_{x \to \infty} \prod_{p \leq x} \left( 1 - \frac{1}{p(p-1)} \right) . \tag{18.17}$$

and explain its accuracy. For example, at $x = 500\,000$, it has the value

$$\prod_{p \leq x} \left( 1 - \frac{1}{p(p-1)} \right) = 0.3739558667768911078453 79 \ldots . \tag{18.18}$$

**Problem 18.4.** Find a way for computing the Artin-Mertens constant for the fixed primitive root 2. The definition gives inaccurate approximation: For $x = 797$,

$$\beta_2 \ = \ \left( \sum_{p \leq x, \ p \in \mathcal{P}_2} \frac{1}{p} - \alpha_2 \log\log x \right) = -0.292545942620583575595681 \ldots , \tag{18.19}$$

which is inaccurate. Explain if it is multiplicative, that is, $\beta_2 = \alpha_2 \beta_1$, and the actual value should be

$$\beta_2 = \alpha_2 \beta_1 = (0.373955813619202288054 \ldots)\,(0.261497212847642783755 \ldots) = 0.0977884029895939748480384 \ldots .$$

**Problem 18.5.** Find a way for computing the Artin-Euler constant for the fixed primitive root 2. The definition gives inaccurate approximation: For $x = 797$,

$$\gamma_2 \ = \ \left( \sum_{p \leq x, \ p \in \mathcal{P}_2} \frac{\log p}{p} - \alpha_2 \log x \right) = 0.431765030391543036913277 \ldots , \tag{18.20}$$

which is inaccurate. Explain if it is multiplicative, that is, $\gamma_2 = \alpha_2 \gamma$, and the actual value should be

$$\gamma_2 = \alpha_2 \gamma = (0.373955813619202288054 \ldots)\,(0.577215664901532860606512 \ldots) = 0.215853184285452242986624 \ldots .$$

**Problem 18.6.** Use the rapidly convergent series, [BC00, p. 251],

$$\gamma \ = \ 1 - \log 2 + \sum_{n \geq 2} \frac{(-1)^n \, (\zeta(n) - 1)}{n} \tag{18.21}$$

to compute Euler constant.

**Problem 18.7.** Use the rapidly convergent series, [BC00, p. 251],

$$B_1 \ = \ 1 - \log 2 + \sum_{n \geq 2} \frac{1}{n} \left( \mu(n) \log \zeta(n) + (-1)^n \, (\zeta(n) - 1) \right) \tag{18.22}$$

to compute Mertens constant.

# Chapter 19

## Prime Harmonic Products

Various prime harmonic products associated with the set of primes $\mathbb{P} = \{\, 2,\ 3,\ 5,\ \ldots \,\}$ have well known estimates.

### 19.1 Prime Harmonic Products Over The Primes

**Lemma 19.1.**    Let $x \geq 1$ be a large number. Then, there exists a pair of constants $\gamma > 0$ and $\kappa_1 > 0$ such that

(i)    $\displaystyle \prod_{p \leq x} \left( 1 - \frac{1}{p} \right)^{-1} = e^{\gamma} \log x + O\left( \frac{1}{\log^2 x} \right).$

(ii)    $\displaystyle \prod_{p \leq x} \left( 1 + \frac{1}{p} \right) = \frac{6\, e^{\gamma}}{\pi^2} \log x + O\left( \frac{1}{\log^2 x} \right).$    (19.1)

(iii)    $\displaystyle \prod_{p \leq x} \left( 1 - \frac{\log p}{p - 1} \right)^{-1} = e^{\nu_1 - \gamma}\, x + O\left( \frac{x}{\log x} \right).$

**Proof**: (iii) Since $0 < \log p / p < 1$, the logarithm of the product can be expressed as

$$\sum_{p \leq x} \log \left( 1 - \frac{\log p}{p - 1} \right)^{-1} = \sum_{p \leq x,\ k \geq 1} \frac{1}{k} \left( \frac{\log p}{p - 1} \right)^{k} = \sum_{p \leq x} \frac{\log p}{p - 1} + \sum_{p \leq x,\ k \geq 2} \frac{1}{k} \left( \frac{\log p}{p - 1} \right)^{k}$$

$$= \log x + \nu_1 - \gamma + O\left( \frac{1}{\log x} \right),$$    (19.2)

where the constant $\nu_1 > 0$ is defined by the double power series

$$\nu_1 = \sum_{p \geq 2,\ k \geq 2} \frac{1}{k} \left( \frac{\log p}{p - 1} \right)^{k} = 0.0152535216156812960614733 \ldots .$$    (19.3)

This proves the assertion in (iii).    ∎

### 19.2 Prime Harmonic Products Over Large Subsets Of Primes

The classical results can be extended in various ways. One of the possible extensions to large subsets of prime numbers of nonzero densities have the following structures.

***Lemma* 19.2.** Let $\alpha = \delta(\mathcal{A}) > 0$ be the density of a subset of primes $\mathcal{A} \subset \mathbb{P}$, and let $x \geq 1$ be a large number. Then, there exist a pair of constants $\gamma_\alpha > 0$ and $\nu_\alpha > 0$ such that

(i) $\displaystyle \prod_{p \leq x, \, p \in \mathcal{A}} \left(1 - \frac{1}{p}\right)^{-1} = e^{\gamma_\alpha} \log^\alpha x + O\left(\frac{1}{\log^2 x}\right).$

(ii) $\displaystyle \prod_{p \leq x, \, p \in \mathcal{A}} \left(1 + \frac{1}{p}\right) = e^{\gamma_\alpha} \prod_{p \in \mathcal{A}} \left(1 - p^{-2}\right) \log^\alpha x + O\left(\frac{1}{\log^2 x}\right).$  (19.4)

(iii) $\displaystyle \prod_{p \leq x, \, p \in \mathcal{A}} \left(1 - \frac{\log p}{p-1}\right)^{-1} = e^{\nu_\alpha - \gamma_\alpha} x^\alpha + O\left(\frac{x^\alpha}{\log x}\right).$

***Proof*:** (i) Express the logarithm of the product as

$$\sum_{p \leq x, \, p \in \mathcal{A}} \log\left(1 - \frac{1}{p}\right)^{-1} = \sum_{p \leq x, \, p \in \mathcal{A}, \, k \geq 1} \frac{1}{k \, p^k} = \sum_{p \leq x, \, p \in \mathcal{A}} \frac{1}{p} + \sum_{p \leq x, \, p \in \mathcal{A}, \, k \geq 2} \frac{1}{k \, p^k}.$$  (19.5)

Applying Lemma 17.2, and Lemma 18.3 yield

$$\begin{aligned} \sum_{p \leq x, \, p \in \mathcal{A}} \log\left(1 - \frac{1}{p}\right)^{-1} &= \alpha \log\log x + \beta_\alpha + O\left(\frac{1}{\log^2 x}\right) + \sum_{p \leq x, \, p \in \mathcal{A}, \, k \geq 1} \frac{1}{k \, p^k} \\ &= \alpha \log\log x + \gamma_\alpha + O\left(\frac{1}{\log^2 x}\right), \end{aligned}$$  (19.6)

compare this to Lemma 19.1. Reversing the logarithm proves the assertion. For (ii) use the expansion

$$\sum_{p \leq x, \, p \in \mathcal{A}} \log\left(1 + \frac{1}{p}\right) = \sum_{p \leq x, \, p \in \mathcal{A}} \log\left(1 - \frac{1}{p^2}\right) + \sum_{p \leq x, \, p \in \mathcal{A}} \log\left(1 - \frac{1}{p}\right)^{-1}.$$  (19.7)

This identity leads to the claim.  ∎

## 19.3 Prime Harmonic Products For Fixed Primitive Roots

The subset of primes $\mathcal{P}_u = \{ \, p \in \mathbb{P} : \mathrm{ord}(u) = p - 1 \, \} \subset \mathbb{P}$ consists of all the primes with a fixed primitive root $u \in \mathbb{Z}$. By Hooley theorem, it has nonzero density $\alpha_u = \delta(\mathcal{P}_u) > 0$. The real number $\alpha_u > 0$ coincides with the corresponding Artin constant, see Chapter 11 or [HY67, p. 220] for the formula.

**Lemma 19.3.** Let $\alpha_u = \delta(\mathcal{P}_u) > 0$ be the density of a subset of primes $\mathcal{P}_u \subset \mathbb{P}$, and let $x \geq 1$ be a large number. Then, there exists a pair of constants $\gamma_u > 0$ and $\kappa_u > 0$ such that

(i) $\displaystyle\prod_{p \leq x,\ p \in \mathcal{P}_u} \left(1 - \frac{1}{p}\right)^{-1} = e^{\gamma_u} \log(x)^{\alpha_u} + O\left(\frac{1}{\log x}\right).$

(ii) $\displaystyle\prod_{p \leq x,\ p \in \mathcal{P}_u} \left(1 + \frac{1}{p}\right) = e^{\gamma_u} \prod_{p \in \mathcal{P}_u} \left(1 - p^{-2}\right) \log(x)^{\alpha_u} + O\left(\frac{1}{\log x}\right).$ (19.8)

(iii) $\displaystyle\prod_{p \leq x,\ p \in \mathcal{P}_u} \left(1 - \frac{\log p}{p-1}\right)^{-1} = e^{\nu_u - \gamma_u} x^{\alpha_u} + O\left(\frac{x^{\alpha_u}}{\log x}\right).$

**Proof:** Express the logarithm of the product as

$$\sum_{p \leq x,\ p \in \mathcal{P}_u} \log\left(1 - \frac{1}{p}\right)^{-1} = \sum_{p \leq x,\ p \in \mathcal{P}_u} \sum_{k \geq 1} \frac{1}{k\,p^k} = \sum_{p \leq x,\ p \in \mathcal{P}_u} \frac{1}{p} + \sum_{p \leq x,\ p \in \mathcal{P}_u} \sum_{k \geq 2} \frac{1}{k\,p^k}.$$ (19.9)

The connecting constants are defined by the Artin-Mertens and Artin-Euler constants respectively:

$$\beta_u = \lim_{x \to \infty} \left(\sum_{p \leq x,\ p \in \mathcal{P}_u} \frac{1}{p} - \alpha_u \log\log x\right) \qquad \text{and} \qquad \beta_u = \gamma_u - \sum_{p \in \mathcal{P}_u,\ k \geq 2} \frac{1}{k\,p^k},$$ (19.10)

where $\alpha_u = \delta(\mathcal{P}_u) > 0$ is the density of the subset of primes $\mathcal{P}_u \subset \mathbb{P}$, see Lemma 17.2. ∎

The *Artin-Mertens* $\beta_u$ and *Artin-Euler* $\gamma_u$ constants are defined by respectively:

$$\beta_u = \lim_{x \to \infty} \left(\sum_{p \leq x,\ p \in \mathcal{P}_u} \frac{1}{p} - \alpha_u \log\log x\right) \qquad \text{and} \qquad \beta_u = \gamma_u - \sum_{p \in \mathcal{P}_u,\ k \geq 2} \frac{1}{k\,p^k},$$ (19.11)

where $\alpha_u = \delta(\mathcal{P}_u) > 0$ is the density of the subset of primes $\mathcal{P}_u \subset \mathbb{P}$. These constants satisfy $\beta_u = \beta_1 \alpha_u$ and $\gamma_u = \gamma \alpha_u$. If the density $\alpha_u = 1$, these definitions reduce to the usual Euler constant and the Mertens constant, which are defined by the limits $\gamma = \lim_{x \to \infty} \left(\sum_{p \leq x} \log p / (p-1) - \log x\right)$, and $\beta_1 = \lim_{x \to \infty} \left(\sum_{p \leq x} 1/p - \log\log x\right)$, or some other equivalent definitions, respectively. Moreover, the linear independence relation becomes $\beta = \gamma - \sum_{p \in \mathcal{P}_u,\ k \geq 2} \left(k\,p^k\right)^{-1}$, see in Lemma 18.1 or [HW08, Theorem 427].

A numerical experiment for the primitive root 2 gives the values

$$\alpha_2 = \prod_{p \geq 2}\left(1 - \frac{1}{p(p-1)}\right) = 0.373955813619202\ldots$$

$$\beta_2 = \sum_{p \leq 1000,\, p \in \mathcal{P}_2} \frac{1}{p} - \alpha_2 \log\log x = 0.167302820886\ldots, \quad \text{and}$$

$$\gamma_2 = \sum_{p \leq 1000,\, p \in \mathcal{P}_2} \frac{\log p}{p-1} - \alpha_2 \log x = 0.424902273366\ldots.$$

The constant $\nu_u > 0$ is defined by the double power series (an approximate the numerical value for $\nu_1$, which is for the set of all primes $\{2,\ 3,\ 5,\ 7,\ 11,\ \ldots,\}$ is shown):

$$\nu_1 = \sum_{p \in \mathcal{P}_u,\, k \geq 2} \frac{1}{k}\left(\frac{\log p}{p-1}\right)^k = 1.141000520953940297743767\ldots. \tag{19.13}$$

## 19.4 Prime Harmonic Products Over Arithmetic Progressions

The materials provided here is concerned with the basic asymptotic structures of the products over the primes in arithmetic progressions. Deeper analysis of the error terms are provided in [LZ09], and similar references. Let $\mathcal{A} \subset \mathbb{P}$ be a subset of primes of nonzero density $\alpha = \delta(\mathcal{A}) > 0$. The products over the subsets of primes $\mathcal{A}(q, a) = \{\, p \in \mathcal{A} : p \equiv a \bmod q \,\} \subset \mathbb{P}$, are estimated in this subsection.

**Lemma 19.4.** Let $x \geq 1$ be a large number. If $a < q$ are relatively prime, and $q = O\big(\log^B x\big)$, $B \geq 0$ constant, then, there exists a pair of constants $\gamma(q, a) > 0$, and $\kappa(q, a) > 0$ such that

(i) $\displaystyle \prod_{p \leq x,\, p \equiv a \bmod q}\left(1 - \frac{1}{p}\right)^{-1} = e^{\gamma(q,a)} \log(x)^{\frac{1}{\varphi(q)}} + O\left(\frac{1}{\log x}\right).$

(ii) $\displaystyle \prod_{p \leq x,\, p \equiv a \bmod q}\left(1 + \frac{1}{p}\right) = e^{\gamma(q,a)} \prod_{p \leq x,\, p \equiv a \bmod q}\left(1 - p^{-2}\right) \log(x)^{\frac{1}{\varphi(q)}} + O\left(\frac{1}{\log x}\right).$ $\tag{19.14}$

(iii) $\displaystyle \prod_{p \leq x,\, p \equiv a \bmod q}\left(1 - \frac{\log p}{p-1}\right)^{-1} = e^{\nu(q,a)-\gamma(q,a)} x^{\frac{1}{\varphi(q)}} + O\left(\frac{x^{\frac{1}{\varphi(q)}}}{\log x}\right).$

The Lehmer-Euler constant is defined by the real number

$$\gamma(q, a) = \lim_{x \to \infty} \left( \sum_{n \leq x,\, n \equiv a \bmod q} \frac{1}{n} - \log^{1/\varphi(q)} x \right) \quad \text{and} \quad \nu(q, a) = \sum_{p \equiv a \bmod q} \sum_{k \geq 2} \frac{1}{k} \left( \frac{\log p}{p - 1} \right)^k \tag{19.15}$$

This result can be extended in various ways. The extension in term of subsets of the prime numbers of nonzero densities have the following shapes.

**Lemma** 19.5.      Let $x \geq 1$ be a large number, and let $a < q$ be relatively prime, with $q = O(\log^B x)$, $B \geq 0$ constant. If $\alpha / \varphi(q) = \delta(\mathcal{A}(q, a)) > 0$ is the density of a subset of primes $\mathcal{A}(q, a) \subset \mathbb{P}$, then, there exists a pair of constants $\gamma_a(q, a) > 0$ and $\kappa_a(q, a) > 0$ such that

(i) $\displaystyle \prod_{p \leq x,\, p \in \mathcal{A}(q,a)} \left( 1 - \frac{1}{p} \right)^{-1} = e^{\gamma_a(q,a)} \log(x)^{\frac{a}{\varphi(q)}} + O\left( \frac{1}{\log x} \right).$

(ii) $\displaystyle \prod_{p \leq x,\, p \in \mathcal{A}(q,a)} \left( 1 + \frac{1}{p} \right) = e^{\gamma_a(q,a)} \prod_{p \in \mathcal{A}(q,a)} \left( 1 - p^{-2} \right) \log(x)^{\frac{a}{\varphi(q)}} + O\left( \frac{1}{\log x} \right).$      (19.16)

(iii) $\displaystyle \prod_{p \leq x,\, p \in \mathcal{A}(q,a)} \left( 1 - \frac{\log p}{p} \right)^{-1} = e^{\gamma_a(q,a) - \gamma_a(q,a)} x^{\frac{a}{\varphi(q)}} + O\left( \frac{x^{\frac{a}{\varphi(q)}}}{\log x} \right).$

**Proof**: Express the logarithm of the product as

$$\sum_{p \leq x,\, p \in \mathcal{A}(q,a)} \log\left( 1 - \frac{1}{p} \right)^{-1} = \sum_{p \leq x,\, p \in \mathcal{A}(q,a),\, k \geq 1} \sum \frac{1}{k\, p^k} = \sum_{p \leq x,\, p \in \mathcal{A}(q,a)} \frac{1}{p} + \sum_{p \leq x,\, p \in \mathcal{A}(q,a),\, k \geq 2} \sum \frac{1}{k\, p^k}. \tag{19.17}$$

The constants are defined by the limit and infinite sum

$$\beta_a(q, a) = \lim_{x \to \infty} \left( \sum_{p \leq x,\, p \in \mathcal{A}(q,a)} \log\left( 1 - \frac{1}{p} \right)^{-1} - \frac{1}{p} \right), \qquad \beta_a(q, a) = \gamma_a(q, a) - \sum_{p \in \mathcal{A}(q,a),\, k \geq 2} \sum \frac{1}{k\, p^k}, \tag{19.18}$$

where $\alpha / \varphi(q) = \delta(\mathcal{A}(q, a)) > 0$ is the density of the subset of primes $\mathcal{A}(q, a) \subset \mathbb{P}$, see the previous examples of the classical definitions.  For (ii), use the identity

$$\sum_{p \leq x,\, p \in \mathcal{A}(q,a)} \log\left( 1 + \frac{1}{p} \right) = \sum_{p \leq x,\, p \in \mathcal{A}(q,a)} \log\left( 1 - \frac{1}{p^2} \right) + \sum_{p \leq x,\, p \in \mathcal{A}(q,a)} \log\left( 1 - \frac{1}{p} \right)^{-1}. \tag{19.19}$$

This identity easily yields the claim.  ∎

## 19.5 Prime Harmonic Products For Fixed Primitive Roots In Arithmetic Progressions

The subset of primes $\mathcal{P}_u(q, a) = \{\, p \in \mathbb{P} : \operatorname{ord}(u) = p - 1 \text{ and } p \equiv a \bmod q \,\} \subset \mathbb{P}$ consists of all the primes with a fixed primitive root $u \in \mathbb{Z}$, and the primes in the arithmetic progression defined by the relatively prime integers $a < q$.

**Lemma** 19.6. Let $x \geq 1$ be a large number, and let $a < q$ be relatively prime integers, with $q = O\big(\log^B x\big)$, $B \geq 0$ constant. If $\alpha_u / \varphi(q) = \delta(\mathcal{P}_u(q, a)) > 0$ is the density of a subset of primes $\mathcal{P}_u(q, a) \subset \mathbb{P}$, then, there exists a pair of constants $\gamma_u(q, a) > 0$ and $\kappa_u(q, a) > 0$ such that

(i) $\displaystyle \prod_{p \leq x,\; p \in \mathcal{P}_u(q,a)} \left(1 - \frac{1}{p}\right)^{-1} = e^{\gamma_u(q,a)} \log(x)^{\frac{\alpha_u}{\varphi(q)}} + O\left(\frac{1}{\log x}\right).$

(ii) $\displaystyle \prod_{p \leq x,\; p \in \mathcal{P}_u(q,a)} \left(1 + \frac{1}{p}\right) = e^{\gamma_u(q,a)} \prod_{p \in \mathcal{P}_u(q,a)} \left(1 - p^{-2}\right) \log(x)^{\frac{\alpha_u}{\varphi(q)}} + O\left(\frac{1}{\log x}\right).$  (19.20)

(iii) $\displaystyle \prod_{p \leq x,\; p \in \mathcal{P}_u(q,a)} \left(1 - \frac{\log p}{p - 1}\right)^{-1} = e^{\kappa_u(q,a) - e^{\gamma_u(q,a)}} x^{\frac{\alpha_u}{\varphi(q)}} + O\left(\frac{x^{\frac{\alpha_u}{\varphi(q)}}}{\log x}\right).$

**Proof**: Express the logarithm of the product as

$$\sum_{p \leq x,\; p \in \mathcal{P}_u(q,a)} \log\left(1 - \frac{1}{p}\right)^{-1} = \sum_{p \leq x,\; p \in \mathcal{P}_u(q,a),\; k \geq 1} \sum \frac{1}{k} \frac{1}{p^k} = \sum_{p \leq x,\; p \in \mathcal{P}_{u(q,a)}} \frac{1}{p} + \sum_{p \leq x,\; p \in \mathcal{P}_u(q,a),\; k \geq 2} \sum \frac{1}{k} \frac{1}{p^k}.$$  (19.21)

The associated constants $\beta_u = \beta_u(q, a)$ and $\gamma_u = \gamma_u(q, a)$ are defined by the limit and infinite sum

$$\beta_u = \lim_{x \to \infty} \sum_{p \leq x,\; p \in \mathcal{P}_u(q,a)} \left(\log\left(1 - \frac{1}{p}\right)^{-1} - \frac{1}{p}\right) = \lim_{x \to \infty} \left(\sum_{p \leq x,\; p \in \mathcal{P}_u(q,a)} \frac{1}{p} - \frac{\alpha_u}{\varphi(q)} \log\log x\right)$$

and

$$\gamma_u(q, a) = \beta_u(q, a) + \sum_{p \in \mathcal{P}_u(q,a),\; k \geq 2} \sum \frac{1}{k} \frac{1}{p^k},$$  (19.22)

where $\alpha_u / \varphi(q) = \delta(\mathcal{P}_u(q, a)) > 0$ is the density of the subset of primes $\mathcal{P}_u(q, a) \subset \mathbb{P}$. ∎

## 19.6 Prime Harmonic Products For The Subset Of Integers With Fixed Quadratic Residues

The subset of primes $\mathcal{P}_{2,r} = \{\, p \in \mathbb{P} : \operatorname{ord}(r) = 2 \,\} \subset \mathbb{P}$ consists of all the primes with a fixed quadratic residue $r \in \mathbb{Z}$. The analysis of the product for this subset of primes is known, but is included here to demonstrate it importance in the analysis of similar subsets of primes.

**Lemma** 19.7. Let $x \geq 1$ be a large number, and let $(r \mid p) = \chi(r)$ be the quadratic residue symbol. If $\rho_{2,r} = \delta(\mathcal{P}_{2,r}) > 0$ is the density of a subset of primes $\mathcal{P}_{2,r}$, then, there exists a pair of constants $\gamma / 2 > 0$ and $\kappa_+ > 0$ such that

(i) $\displaystyle \prod_{p \leq x, \, (e \mid p) = 1} \left(1 - \frac{1}{p}\right)^{-1} = e^{\gamma/2} \prod_{p \mid r} \left(1 - p^{-1}\right)^{1/2} \prod_{(r \mid p) = -1} \left(1 - p^{-2}\right)^{1/2} L(1, \chi) \log(x)^{1/2} \left(1 + O\left(\frac{1}{\log\log x}\right)\right).$

(ii) $\displaystyle \prod_{p \leq x, \, (e \mid p) = 1} \left(1 + \frac{1}{p}\right) = e^{\gamma/2} \prod_{p \mid r} \left(1 - p^{-1}\right)^{1/2} \prod_{p \geq 2} \left(1 - p^{-2}\right)^{1/2} \prod_{(r \mid p) = 1} \left(1 - p^{-2}\right)^{1/2} L(1, \chi) \log(x)^{1/2} \left(1 + O\left(\frac{1}{\log\log x}\right)\right).$

(iii) $\displaystyle \prod_{p \leq x, \, (e \mid p) = 1} \left(1 - \frac{\log p}{p - 1}\right)^{-1} = e^{\nu_+ - e^{\gamma/2}} x^{\frac{1}{2}} + O\left(\frac{x^{\frac{1}{2}}}{\log x}\right).$

$$(19.23)$$

## 19.7 Problems And Exercises

**Problem 19.1.** Let $\alpha = \delta(\mathcal{A}) > 0$ be the density of a subset of primes $\mathcal{A} \subset \mathbb{P}$, and let $x \geq 1$ be a large number. Show that

$$\prod_{p \leq x, \, p \in \mathcal{A}} \left(1 - \frac{1}{p}\right) = \prod_{p \leq x} \left(1 - \frac{\alpha}{p}\right). \tag{19.24}$$

**Problem 19.2.** Compute approximations of the constants

$$\nu_- = \sum_{(2 \mid p) = -1, \, n \geq 2} \frac{1}{n} \left(\frac{\log p}{p - 1}\right)^n \quad \text{and} \quad \nu_+ = \sum_{(2 \mid p) = 1, \, n \geq 2} \frac{1}{n} \left(\frac{\log p}{p - 1}\right)^n. \tag{19.25}$$

## Chapter 20

# Primes Decompositions Formulae

The finite product

$$\sum_{n \leq x, \ p \mid n \Rightarrow p \leq x} \frac{1}{n} = \prod_{p \leq x} \left(1 - \frac{1}{p}\right)^{-1} \tag{20.1}$$

is a sort of a basic prototype for primes decompositions formulas. The primes decompositions of certain summatory multiplicative functions as finite products over the primes supports of the functions are achieved by several important results in the literature.

## 20.1. Wirsing Formula

This formula provides a decompositions of some summatory multiplicative functions as products over the primes supports of the functions. This technique works well with certain multiplicative functions, which have supports on subsets of primes numbers of nonzero densities.

**Lemma** 20.1.   (Wirsing)    Suppose that $f : \mathbb{N} \longrightarrow \mathbb{C}$ is a multiplicative function with the following properties.

$(i)$ $f(n) \geq 0$          for all integers $n \in \mathbb{N}$.

$(ii)$ $f(n) \geq c^k$,          for all integers $k \in \mathbb{N}$, and $c > 0$ constant.

$(iii)$ $\sum_{p \leq x} f(p) = (\tau + o(1)) \, x / \log x,$      where $\tau > 0$ is a constant as $x \longrightarrow \infty.$      (20.2)

Then

$$\sum_{n \leq x} f(n) = \left(\frac{1}{e^{\gamma \tau} \, \Gamma(\tau)} + o(1)\right) \frac{x}{\log x} \prod_{p \leq x} \left(1 + \frac{f(p)}{p} + \frac{f(p^2)}{p^2} + \cdots\right). \tag{20.3}$$

The gamma function appearing in the above formula is defined by $\Gamma(s) = \int_0^\infty t^{s-1} \, e^{-st} \, d \, t, \, s \in \mathbb{C}$. The intricate proof of Wirsing formula is assembled in various papers, such as [HD87], [PV88, p. 195], and discussed in [MV07, p. 70], [TM95, p. 308]. Various applications are provided in [MP11], [PP03], et alii.

## 20.2. Rieger Formula

A related but less restrictive prime decomposition formula was developed around the same time.

**Lemma 20.2.** (Rieger) Let $\mathcal{T}$ be a subset of primes of density $\tau = \delta(\mathcal{T}) > 0$, and let

$$\sum_{p \leq x, \, p \in \mathcal{T}} \frac{\log p}{p} = \tau \log x - \gamma_\tau + O\left(e^{-c\sqrt{\log x}}\right), \tag{20.4}$$

where $c > 0$ is an absolute constant. Then

$$\sum_{n \leq x, \, p \mid n \Rightarrow p \in \mathcal{T}} \frac{1}{n} \sim \frac{1}{e^{\gamma \tau} \, \Gamma(\tau + 1)} \prod_{p \leq x, \, p \in \mathcal{T}} \left(1 - \frac{1}{p}\right)^{-1}. \tag{20.5}$$

An application of this theorem is given in [WS74, p. 356].

## 20.3. Delange Formula

Another related prime decomposition formula usable for other types of subsets of primes and multiplicative functions is given below. This result is concerned with the mean value

$$\lim_{x \longrightarrow \infty} \frac{1}{x} \sum_{n \leq x} f(n) \tag{20.6}$$

of multiplicative functions $f : \mathbb{N} \longrightarrow \mathbb{C}$ on the complex unit disk $D = \{\, z \in \mathbb{C} : |z| \leq 1 \,\}$.

**Lemma 20.3.** (Delange) Let $f : \mathbb{N} \longrightarrow \mathbb{C}$ be a multiplicative function. Assume that

(i) $|f(n)| \leq 1$                         for all integers $n \in \mathbb{N}$.

(ii) $\displaystyle\sum_{p \leq x} \frac{1 - f(p)}{p} < \infty.$ \hfill (20.7)

Then

$$\sum_{n \leq x} f(n) = x \prod_{p \leq x} \left(1 - \frac{1}{p}\right)\left(1 + \frac{f(p)}{p} + \frac{f(p^2)}{p^2} + \cdots\right) + o(x). \tag{20.8}$$

The product in the above expression converges to a complex number, which depends on the function $f$.

## 20.4. Halasz Formula

The generalization of Delange theorem to complex valued multiplicative functions widen the sphere of applicability of these results, a proof appears in [TM95, p. 335].

**Lemma 20.4.** (Halasz) Let $f : \mathbb{N} \longrightarrow \mathbb{C}$ be a complex valued multiplicative function such that

(i) $| f(n) | \leq 1$                          for all integers $n \in \mathbb{N}$.

(ii) $\displaystyle\sum_{p \leq x} \frac{1 - \mathcal{R}e(f(p)) \, p^{-i\tau}}{p} < \infty$        converges for some $\tau \in \mathbb{R}$.             (20.9)

Then

(iii) $\displaystyle\frac{1}{x} \sum_{n \leq x} f(n) = \frac{x^{i\tau}}{1 + i\,\tau} \prod_{p \leq x} \left(1 - \frac{1}{p}\right)\left(1 + \frac{f(p)}{p} + \frac{f\left(p^2\right)}{p^2} + \cdots\right) + o(x)\,,$     as $x \longrightarrow \infty$.

(iv) $\displaystyle\frac{1}{x} \sum_{n \leq x} f(n) = o(1)\,,$         if there is no some $\tau \in \mathbb{R}$ such that (ii) converges as $x \longrightarrow \infty$.        (20.10)

Chapter 21

# Characteristic Functions In The Finite Rings $\mathbb{Z} / p^k \, \mathbb{Z}$

Let $\mathbb{Z} = \{ 0, \ \pm 1, \ \pm 2, \ \dots \}$ denotes the set of integers, and let $\mathbb{P} = \{ 2, \ 3, \ 5, \ \dots \}$ denotes the set of prime numbers. Some materials on the characteristic functions of certain elements in the finite rings $\mathbb{Z} / p^k \, \mathbb{Z} = \{ \, n \bmod p^k : n \in \mathbb{Z} \, \}$ are determined in this section.

## 21.1 Characteristic Function of Quadratic Nonresidues in the Finite Rings $\mathbb{Z} / p^k \, \mathbb{Z}$

The subset of primes $\mathcal{K}_{2,r} = \{ \, p \in \mathbb{P} : \operatorname{ord}(r) = 2 \, \} \subset \mathbb{P}$ with a fixed quadratic nonresidue $r \in \mathbb{Z}$ has a nonzero density $\kappa_{2,r} = \delta(\mathcal{K}_{2,r}) = 1/2$. This is implied by the quadratic reciprocity law. The characteristic function of a fixed quadratic nonresidue $r \in \mathbb{Z}$ over the finite rings $\mathbb{Z} / p^k \, \mathbb{Z}$, $k \geq 1$, the integers modulo $p^k$, is investigated here. Given a fixed integer $r \in \mathbb{Z}$, and a prime power $p^k$ such that $\gcd(r, \ p^k) = 1$, $k \geq 1$. Let

$$\left( \frac{r}{p^k} \right) = \begin{cases} 1 & \text{if } r \text{ is a quadratic residue,} \\ -1 & \text{if } r \text{ is a quadratic nonresidue,} \end{cases} \tag{21.1}$$

be the quadratic symbol mod $p^k$, and let $\operatorname{ord} r$ denotes the order of the element $r \in \mathbb{Z}_{p^k}^*$ in the multiplicative group of the integers modulo $p^k$. The quadratic symbol is a different way of identifying the elements of order $\operatorname{ord} r = 2$.

**Lemma 21.1.** Let $p^k$, $k \geq 1$ be a prime power, and let $r \in \mathbb{Z}$ be an integer such that $\gcd(r, \ p^k) = 1$. Then,

(i) The characteristic function of the quadratic nonresidue $r \bmod p^k$ is given by

$$f_2(p^k) = \begin{cases} 1 & \text{if } \left( \frac{r}{p^k} \right) = -1, \ k \geq 1, \\ 0 & \text{if } \left( \frac{r}{p^k} \right) = 1, \ k \geq 1. \end{cases} \tag{21.2}$$

(ii) The function $f_2 : \mathbb{N} \longrightarrow \{ 0, \ 1 \}$ is completely multiplicative:

$$\begin{aligned} &\text{(i)} \, f(p \, q) \ = f(p) \, f(q), & \gcd(p, \, q) = 1, \\ &\text{(ii)} \, f(p^2) = f(p) \, f(p), & \text{if } \operatorname{ord} u \ = 2 \bmod p^k. \end{aligned} \tag{21.3}$$

**Proof**: From the information $\left( \frac{r}{p} \right) \equiv r^{(p-1)/2} \equiv -1 \bmod p$, and $r^{p^{k-1}(p-1)/2} = (-1 + a \, p)^{p^{k-1}} \equiv -1 \bmod p^k$, it follows that $\left(\frac{r}{p^k}\right) \equiv r^{p^{k-1}(p-1)/2} \equiv -1 \bmod p^k$ for all $k \geq 1$. Thus, the function has the value $f\left(p^k\right) = 1$ if and only if the element $r \in \mathbb{Z}_{p^k}^*$ is a quadratic nonresidue modulo $p^k$. Otherwise, it vanishes: $f\left(p^k\right) = 0$. ∎

## 21.2 Characteristic Function of Cubic Nonresidues in the Finite Rings $\mathbb{Z}/p^k\mathbb{Z}$

The subset of primes $\mathcal{K}_{3,r} = \{\, p \in \mathbb{P} : \mathrm{ord}(r) = 3 \,\} \subset \mathbb{P}$ with a fixed quartic nonresidue $r \in \mathbb{Z}$ has a nonzero density $\kappa_{3,r} = \delta(\mathcal{K}_{3,r})$. This is implied by the cubic reciprocity law.

The characteristic function of a fixed cubic nonresidue $r \in \mathbb{Z}$ in the finite rings $\mathbb{Z}/p^k\mathbb{Z}$, $k \geq 1$, the integers modulo $p^k$, is determined here. Given a fixed integer $r \in \mathbb{Z}$, and a prime power such that $\gcd(r, p^k) = 1$, $k \geq 1$. Let

$$\left(\frac{r}{p^k}\right)_3 = \begin{cases} \pm 1 & \text{if } r \text{ is a cubic residue,} \\ \omega^{\pm 1} & \text{if } r \text{ is a cubic nonresidue,} \end{cases} \tag{21.4}$$

where $\omega \neq 1$ is a cubic root of unity, be the cubic symbol $\bmod\, p^k$, and let $\mathrm{ord}\, r$ denotes the order of the element $r \in \mathbb{Z}_{p^k}^*$ in the multiplicative group of the integers modulo $p^k$. The cubic symbol is a different way of identifying the elements of order $\mathrm{ord}\, r = 3$.

***Lemma*** **21.2.** Let $p^k$, $k \geq 1$ be a prime power, and let $r \in \mathbb{Z}$ be an integer such that $\gcd(r, p^k) = 1$. Then

(i) The characteristic function of the cubic nonresidue $r \bmod p^k$ is given by

$$f_3\left(p^k\right) = \begin{cases} 1 & \text{if } \left(\frac{r}{p^k}\right)_3 = \omega, \ k \geq 1, \\ 0 & \text{if } \left(\frac{r}{p^k}\right)_3 = \pm 1, \ k \geq 1. \end{cases} \tag{21.5}$$

(ii) The function $f_3 : \mathbb{N} \longrightarrow \{\, 0, \ 1 \,\}$ is completely multiplicative:

$$\begin{aligned} &\text{(i) } f(p\,q) = f(p)\,f(q), && \gcd(p, q) = 1, \\ &\text{(ii) } f\left(p^2\right) = f(p)\,f(p), && \text{if } \mathrm{ord}\, r = 3 \bmod p^k. \end{aligned} \tag{21.6}$$

***Proof***: From the information $\left(\frac{r}{p}\right)_3 \equiv r^{(p-1)/3} \equiv \omega^{\pm 1} \bmod p$, and $r^{p^{k-1}(p-1)/3} = \left(\omega^{\pm 1} + a\,p\right)^{p^{k-1}} \equiv \omega^{\pm 1} \bmod p^k$, it follows that $\left(\frac{r}{p^k}\right)_3 \equiv r^{p^{k-1}(p-1)/3} \equiv \omega^{\pm 1} \bmod p^k$ for all odd prime powers $p^k$, $k \geq 1$. Thus, the function has the value $f\left(p^k\right) = 1$ if and only if the element $r \in \mathbb{Z}_{p^k}^*$ is a cubic nonresidue modulo $p^k$. Otherwise, it vanishes: $f\left(p^k\right) = 0$. ∎

## 21.3 Characteristic Function of Quartic Nonresidues in the Finite Rings $\mathbb{Z}/p^k\mathbb{Z}$

The subset of primes $\mathcal{K}_{4,r} = \{\, p \in \mathbb{P} : \mathrm{ord}(r) = 4 \,\} \subset \mathbb{P}$ with a fixed quartic nonresidue $r \in \mathbb{Z}$ has a nonzero density $\kappa_{4,r} = \delta(\mathcal{K}_{4,r})$. This is implied by the biquadratic reciprocity law.

The characteristic function of a fixed quartic nonresidue $r \in \mathbb{Z}$ in the finite rings $\mathbb{Z}/p^k\mathbb{Z}$, $k \geq 1$, the integers modulo $p^k$, is determined here. Given a fixed integer $r \in \mathbb{Z}$, and a prime power such that $\gcd(r, p^k) = 1$, $k \geq 1$. Let

$$\left(\frac{r}{p^k}\right)_4 = \begin{cases} \pm 1 & \text{if } r \text{ is a quartic residue,} \\ \pm i & \text{if } r \text{ is a quartic nonresidue,} \end{cases} \tag{21.7}$$

be the quartic symbol $\bmod\, p^k$, and let $\operatorname{ord} r$ denotes the order of the element $r \in \mathbb{Z}_{p^k}^*$ in the multiplicative group of the integers modulo $p^k$. The quartic symbol is a different way of identifying the elements of order $\operatorname{ord} r = 4$.

**Lemma 21.3.** Let $p^k$, $k \geq 1$, be a prime power, and let $r \in \mathbb{Z}$ be an integer such that $\gcd(r, p^k) = 1$. Then

(i) The characteristic function of the quartic nonresidue $r \bmod p^k$ is given by

$$f_4(p^k) = \begin{cases} 1 & \text{if } \left(\frac{r}{p^k}\right) = -1, \ k \geq 1, \\ 0 & \text{if } \left(\frac{r}{p^k}\right) = 1, \ k \geq 1. \end{cases} \tag{21.8}$$

(ii) The function $f_4 : \mathbb{N} \longrightarrow \{0, 1\}$ is completely multiplicative:

$$\begin{array}{lll} \text{(i) } f(p\,q) = f(p)\,f(q), & \gcd(p, q) = 1, \\ \text{(ii) } f(p^2) = f(p)\,f(p), & \text{if } \operatorname{ord} r = 4 \bmod p^k. \end{array} \tag{21.9}$$

**Proof**: From the information $\left(\frac{r}{p}\right)_4 \equiv r^{(p-1)/4} \equiv \pm i \bmod p$, and $r^{p^{k-1}(p-1)/4} \equiv (\pm i + a\,p)^{p^{k-1}} \equiv \pm i \bmod p^k$, it follows that $\left(\frac{r}{p^k}\right)_4 \equiv r^{p^{k-1}(p-1)/4} \equiv \pm i \bmod p^k$ for all odd prime powers $p^k$, $k \geq 1$. Thus, the function has the value $f(p^k) = 1$ if and only if the element $r \in \mathbb{Z}_{p^k}^*$ is a quartic nonresidue modulo $p^k$. Otherwise, it vanishes: $f(p^k) = 0$. ∎

## 21.4 Characteristic Function of Primitive Roots in the Finite Rings $\mathbb{Z}/p^k\mathbb{Z}$

The symbol $\operatorname{ord} u$ denotes the order of the element $u \in \mathbb{Z}_{p^k}^*$ in the multiplicative group of the integers modulo $p^k$. The subset of primes $\mathcal{P}_u = \{p \in \mathbb{P} : \operatorname{ord}(u) = p - 1\} \subset \mathbb{P}$ with a fixed primitive root $u \in \mathbb{Z}$ such that $u \neq \pm 1$, $v^2$, is expected to have a nonzero density $\alpha_u = \delta(\mathcal{P}_u) > 0$. The characteristic function of primitive root for the finite ring $\mathbb{Z}_{p^k}$, the integers modulo $p^k$, is determined here.

**Lemma 21.4.** Let $p^k$, $k \geq 1$, be a prime power, and let $u \in \mathbb{Z}$ be an integer such that $\gcd(u, p^k) = 1$. Then

(i) The characteristic function of the primitive root $u \bmod p^k$ is given by

$$f\left(p^k\right) = \begin{cases} 1 & \text{if } p^k = 2^k, \ k \leq 2, \\ 0 & \text{if } p^k = 2^k, \ k > 2, \\ 1 & \text{if ord } u = p^{k-1}(p-1), \ p > 2, \ \text{for any } k \geq 1, \\ 0 & \text{if ord } u \neq p^{k-1}(p-1), \ \text{and } p > 2, \ k \geq 2. \end{cases}$$

(ii) The function $f : \mathbb{N} \longrightarrow \{0, \ 1\}$ is multiplicative, but not completely multiplicative since

$$\begin{aligned} &\text{(i) } f(p\,q) = f(p)\,f(q), &&\gcd(p, q) = 1, \\ &\text{(ii) } f(p^2) \neq f(p)\,f(p), &&\text{if ord } u \neq p(p-1) \bmod p^2. \end{aligned} \qquad (21.11)$$

***Proof***: The function has the value $f\left(p^k\right) = 1$ if and only if the element $u \in \mathbb{Z}_{p^k}^*$ is a primitive root modulo $p^k$. Otherwise, it vanishes: $f\left(p^k\right) = 0$. ∎

## 21.5 Problems And Exercises

**Problem 21.1.** Determine the characteristic function of primitive roots in the finite ring $\mathbb{Z} \big/ p^k \mathbb{Z}$, $k \geq 1$, for primes in the arithmetic progression $\mathcal{P}_u = \{\, p \in \mathbb{P} : p \equiv a \bmod q \text{ and ord}(u) = p - 1 \,\}$ with a fixed primitive root $u \in \mathbb{Z}$, and $a < q$, $\gcd(a, q) = 1$.

**Problem 21.2.** Determine the characteristic function of squarefree totient $p$ prime and $\mu(p - 1) \neq 0$.

**Problem 21.3.** Find an expression for the density of fixed cubic nonresidue modulo primes, and integers.

**Problem 21.34.** Find an expression for the density of fixed quartic nonresidue modulo primes, and integers.

Chapter 22

# Primitive Roots In The Finite Rings $\mathbb{Z} / p^k \mathbb{Z}$

The primitive roots of the finite rings of the integers modulo $p^k$, $k \geq 1$, have some exceptional structures of interest in the analysis of the fixed primitive roots of the integers. Some of these structures are described here.

## 22.1 Primitive Roots Test and Primality Tests

***Lemma* 22.1.**    (Gauss)   The integers $N = 2, 4, p^m, 2\,p^m$,  $p$ prime and $m \geq 1$, are the only integers which have primitive roots mod $N$.

This is a standard result in elementary number theory. A proof is available in [RE95, p. 91], et alii.

***Lemma* 22.2.**    (Gauss Primitive Root Test)   Let $p \geq 2$ be an odd prime, and let $u \in \mathbb{Z}$ be an integer. Then, $u$ is primitive element in the finite ring $\mathbb{Z} / p^k \mathbb{Z}$, $k \geq 1$, if and only if

$$u^{p^{k-1}(p-1)/q} \not\equiv 1 \bmod p^k \tag{22.1}$$

for all prime divisors $q \mid p^{k-1}(p-1)$.

A few versions of this result are explicated in [GF86, Articles 55 to 58, p. 35]. A few primality tests are simple extension of the Gauss primitive root test. One of these is known as Lucas primality test. The precise claim is follows.

***Lemma* 22.3.**    (Lucas)   Let $n \geq 2$ be an integer, and let $u \in \mathbb{Z}$ be a fixed integer. Then, $n$ is a prime if and only if

(i) $u^{n-1} \equiv 1 \bmod n$ ,                    (ii) $u^{(n-1)/q} \not\equiv 1 \bmod n$ ,          (22.2)

for all prime divisors $q \mid (n-1)$.

Other specialized primality tests are derived from these cluster of ideas, these are explicated in [CP05, p. 174], et alii.

The Gauss test is concerned with the primitivity of a fixed integer $u \neq \pm 1$, $v^2$ with respect to a single prime $p \geq 2$, while Artin conjecture for primitive roots is concerned with the primitivity of a fixed integer $u \neq \pm 1$, $v^2$ with respect to an infinite subset of primes $p \in \{p : \operatorname{ord}(u) = p - 1\}$ .

**Lemma** 22.4.    (i)  Let $p \geq 2$ be an odd prime, and let $r_1$, $r_2$, ..., $r_{\varphi(p-1)}$ be the primitive roots modulo $p$. Then, the correspondence

$$r_i \longrightarrow s_{i,j} = r_i + p\, x_j, \tag{22.3}$$

where $x_j \in \mathbb{Z}\,/\,p\,\mathbb{Z}$, $1 \leq i \leq \varphi(p-1)$, and $0 \leq j < p-1$, generates the set of primitive roots modulo $p^2$.

(ii)  Let  $p^k$, $k \geq 2$,  be an  odd prime  power,  and let $s_1$, $s_2$, ..., $s_m$  be the  primitive  roots modulo  $p^k$, $m = p^{k-2}(p-1)\,\varphi(p-1)$. Then, the correspondence

$$s_i \longrightarrow t_{i,j} = s_i + p^k x_j, \tag{22.4}$$

where $x_j \in \mathbb{Z}\,/\,p\,\mathbb{Z}$, $1 \leq i \leq \varphi(p-1)$,  and $0 \leq j \leq p-1$, generates the set of primitive roots modulo $p^{k+1}$.

This result is proved in [AP86, p. 209], [NZ66, Theorem 2.39], [RE95].

Observe that the finite ring $\mathbb{Z}\,\big/\,p^2\,\mathbb{Z}$ of integers modulo $p^2$ has $\varphi\big(\varphi(p^2)\big) = (p-1)\,\varphi(p-1)$ primitive roots, so there is an extra value among the integers 0,  1,  2,  3,  ...,  $p-1$. This exceptional value has a precise formulation as

$$x \equiv \frac{1-g^{p-1}}{p}\,\big((p-1)\,g^{p-2}\big)^{-1} \bmod p\,, \tag{22.5}$$

where $g = r_i$ is a primitive root mod $p$. If exceptional value $x = 0$, then the primitive root $g \bmod p$ cannot be lifted to a primitive root mod $p^2$. This condition is often characterized by the congruence equation $u^{p-1} \equiv 1 \bmod p^2$. Given a fixed primitive root $u$, very few such *sporadic primes* are known.

## 22.2 Abel, and Wieferich Primes

For a fixed integer $u \in \mathbb{Z}$, the prime solutions $p \in \mathbb{P}$ of the congruence $u^{p-1} \equiv 1 \bmod p^2$ was investigated by Abel, [RN, p. 334]. This topic has a two centuries long history, and important applications in mathematics, confer the literature for more details. However, the nature of the subsets of primes $\big\{\, p \in \mathbb{P} : u^{p-1} \equiv 1 \bmod p^2 \,\big\}$ remains unknown. The primes solutions of the special case $u = 2$, inter alia, are known as Wieferich primes.

Some of the best known applications are

(1) The criterion $a^{p-1} \equiv 1 \bmod p^2$ for the resolution of the Fermat equation $x^n + y^n = z^n$.

(2) The dual criteria $q^{p-1} \equiv 1 \bmod p^2$   and   $p^{q-1} \equiv 1 \bmod q^2$ used in the resolution of the Catalan equation $x^p + y^q = 1$.

(3) The ratio $q_p(u) = \big(u^{p-1} - 1\big)\big/\,p$, called Fermat quotient, is used in the theory of the discrete logarithm. It has logarithm

type properties:

(i) $q_p(p-1) \equiv 1 \bmod p$,

(ii) $q_p(p+1) \equiv -1 \bmod p$,

(iii) $q_p(u\,v) \equiv q_p(u) + q_p(v) \bmod p$. $\qquad$ (22.6)

The associated Diophantine equation

$$u^{X-1} = 1 + X^2\,Y\,, \qquad (22.7)$$

$|u| > 1$ fixed, should have a finite number of integers solutions $X = p$ prime, $Y \in \mathbb{Z}$. A few of the known primes solutions are listed below.

$u = 2, \ \ p = 1093,\ 3511,\ \dots$.

$u = 3, \ \ p = 11,\ 1\,006\,003,\ \dots$

$u = 5, \ \ p = 2,\ 20\,771,\ 40\,487,\ 53\,471\,161,\ 1\,645\,333\,507,\ 6\,692\,367\,337,\ 188\,748\,146\,801,\ \dots$.

$u = 6, \ \ p = 66\,161,\ 534\,851,\ 3\,152\,573,\ \dots$.

***Theorem* 22.5.** ([SN88]) Fix $u \in \mathbb{Q}^*$, $u \neq \pm 1$. If the abc conjecture is true, then

$$\#\left\{\, p \leq x : u^{p-1} \not\equiv 1 \bmod p^2 \,\right\} \gg \log x \ \text{as} \ x \longrightarrow \infty. \qquad (22.8)$$

An extension of the Abel-Wieferich condition to the more general groups of rational points of algebraic curves of genus $g = 1$ is given in the same paper.

## 22.3 Correction Factor

The sporadic subsets of Abel-Wieferich primes, see [RN96] for other details, have roles in the determination of the densities of the subsets of integers $\mathcal{N}_u = \{\, n \in \mathbb{N} \ \text{and} \ \mathrm{ord}(u) = \lambda(n) \,\}$ with fixed primitive roots $u \in \mathbb{Z}$. The finite prime product arising from the sporadic existence of the Abel-Wieferich primes $p \geq 3$ is implemented by the equivalent expression

$$\prod_{\substack{p \leq x,\ \mathrm{ord}(u)=p-1, \\ \mathrm{ord}(u)\neq p(p-1)}} \left(1 + \frac{1}{p}\right) \prod_{\substack{p \leq x, \\ \mathrm{ord}(u)=p(p-1)}} \left(1 + \frac{1}{p} + \frac{1}{p^2} + \cdots\right) = \prod_{\substack{p \leq x,\ \mathrm{ord}(u)=p-1, \\ \mathrm{ord}(u)\neq p(p-1)}} \left(1 - \frac{1}{p^2}\right) \prod_{\substack{p \leq x, \\ \mathrm{ord}(u)=p-1}} \left(1 - \frac{1}{p}\right)^{-1}. \qquad (22.9)$$

Note that ord $u = p - 1 \bmod p$, but ord $u \neq p(p-1) \bmod p^2$ is equivalent to the usual characterization of Abel-Wieferich primes $p \geq 3$ such that $u^{p-1} \equiv 1 \bmod p^2$.

This finite product seems to be a correction factor similar to case for primitve roots over the prime numbers. The correction required for certain densities of primes with respect to fixed primitve roots was discovered by the Lehmers, see [ST03].

## 22.4 Problems And Exercises

**Problem 22.1.** A heuristic argument shows that an Abel-Wieferich prime $p \geq 3$ such that $u^{p-1} \equiv 1 \bmod p^2$ exists with probability $1/p$. Explain why this argument contradicts the conjecture number of solutions of the Diophantine equation $u^{X-1} = 1 + X^2\,Y$, $|u| > 1$ fixed, and $X = p$ prime, $Y \in \mathbb{Z}$. Hint: $\sum_{p \leq x} 1/p > .25 \log\log x$.

## Chapter 23

# Generalized Artin Conjecture

A generalized Artin conjecture for subsets of composite integers with fixed primitive roots is presented in this chapter. This analysis spawn new questions about the structure of the $L$-series associated with the multiplicative subset of integers with a fixed primitive root, and related ideas.

## 23.1 Subset of Integers With Fixed Primitive Root 2

The underlining structure of an asymptotic counting formula $N_2(x) = \#\{n \le x : \mathrm{ord}(2) = \lambda(n)\}$ for the subset of integers $\mathcal{N}_2 = \{n \in \mathbb{N} : \mathrm{ord}(2) = \lambda(n)\}$ will be demonstrated here. This analysis, based on Theorem 12.1 or Theorem 11.1, see [HY67]. This theorem states that the integer 2 is a primitive root mod $p$ for infinitely many primes. Id est,

$$\mathcal{P}_2 = \{3,\ 5,\ 11,\ 13,\ 19,\ 29,\ 37,\ 53,\ 59,\ 61,$$
$$67,\ 83,\ 101,\ 107,\ 107,\ 131,\ 139,\ 149,\ 163,\ 173,\ 179,\ 181,\ 197,\ ...\,\}. \tag{23.1}$$

The subset of primes $\mathcal{P}_2 = \{p \in \mathbb{P} : \mathrm{ord}(2) = p - 1\} \subset \mathbb{P}$ of density $\alpha_2 = \delta(\mathcal{P}_2) > 0$, which consists of all the primes with a fixed primitive root $2 \in \mathbb{Z}$, is utilized to generate the subset of composite integers

$$\mathcal{N}_2 = \{3,\ 5,\ 3^2,\ 11,\ 3,\ 3 \cdot 5,\ 19,\ 5^2,\ 3^3,\ 29,\ 3 \cdot 11,\ 37,\ 3 \cdot 13,\ 45,\ 53,\ 55,\ 61,\ 65,\ 67,\ 3 \cdot 5^2,\ ...\,\}, \tag{23.2}$$

which have 2 as a primitive root. A few associated constants arise in this analysis. These are

$$\frac{1}{2\,\Gamma(\alpha_2)} \prod_{\mathrm{ord}(2)=p-1,\,\mathrm{ord}(2)\ne p(p-1)} \left(1 - \frac{1}{p^2}\right) = 0.210324810565075\,..., \tag{23.3}$$

where $\Gamma(s) = \int_0^\infty t^{s-1}\,e^{-s\,t}\,d\,t$, $s \in \mathbb{C}$, is the gamma function, and

$$\alpha_2 = \prod_{p \ge 2}\left(1 - \frac{1}{p(p-1)}\right) = 0.373955813619202\,... \tag{23.4}$$

is Artin constant, and the *Artin-Euler* constant, which is defined by

$$\gamma_2 \ = \ \lim_{x \longrightarrow \infty} \left( \sum_{p \leq x, \, p \in \mathcal{P}_2} \frac{\log p}{p-1} - \alpha_2 \log x \right), \tag{23.5}$$

confer Lemma 15.4, Definition 18.2, and Lemma 19.3 for details. This result is consistent with the heuristic explained in [LP03, p. 10].

**Theorem 23.1.** The integer 2 is primitive root mod $n$ for infinitely many composite integers $n \geq 1$. Moreover, the number of integers $n \leq x$ such that 2 is a primitive root mod $n$ has the asymptotic formulae

(i) $N_2(x) = \left( \dfrac{1}{2 \, \Gamma(\alpha_2)} + o(1) \right) \dfrac{x}{(\log x)^{1-\alpha_2}} \displaystyle\prod_{\substack{\mathrm{ord}(2)=p-1, \\ \mathrm{ord}(2) \neq p(p-1)}} \left( 1 - \dfrac{1}{p^2} \right).$

(ii) $N_2(x) \geq \left( \dfrac{3}{\pi^2 \, \Gamma(\alpha_2)} + o(1) \right) \dfrac{x}{(\log x)^{1-\alpha_2}},$

$$\tag{23.6}$$

where $\alpha_2 = .373955813619202 \ldots$ is Artin constant, and the index of the product ranges over the Wieferich primes $2^{p-1} \equiv 1 \bmod p^2$, for all large numbers $x \geq 1$.

**Proof**: By Theorem 12.1, the density $\alpha_2 > 0$ is nonzero. Put $\tau = \alpha_2$, in Wirsing formula with $\tau = \alpha_2$, Lemma 20.1, and replace the characteristic function $f(n)$ of primitive roots in the finite ring $\mathbb{Z}/n\mathbb{Z}$, see Lemma 21.4, to produce

$$\begin{aligned}
\sum_{n \leq x} f(n) &= \left( \frac{1}{e^{\gamma \tau} \, \Gamma(\tau)} + o(1) \right) \frac{x}{\log x} \prod_{p \leq x} \left( 1 + \frac{f(p)}{p} + \frac{f(p^2)}{p^2} + \cdots \right) \\
&= \frac{\varphi(2)}{2} \left( \frac{1}{e^{\gamma \alpha_2} \, \Gamma(\alpha_2)} + o(1) \right) \frac{x}{\log x} \prod_{\substack{p \leq x, \, \mathrm{ord}(2)=p-1, \\ \mathrm{ord}(2) \neq p(p-1)}} \left( 1 + \frac{1}{p} \right) \prod_{\substack{p \leq x, \\ \mathrm{ord}(2)=p(p-1)}} \left( 1 + \frac{1}{p} + \frac{1}{p^2} + \cdots \right),
\end{aligned} \tag{23.7}$$

where the index $p \in \mathcal{P}_2$ runs over the primes such that 2 is a primitive root modulo $p$ or modulo $p^2$. The ratio $\varphi(2)/2$ accounts for the condition $\gcd(2, n) = 1$. Replacing the equivalent product, see Section 22.3, and using Lemma 19.3, yield

$$\sum_{n \leq x} f(n) = \frac{1}{2} \left( \frac{1}{e^{\gamma \alpha_2} \, \Gamma(\alpha_2)} + o(1) \right) \frac{x}{\log x} \prod_{\substack{p \leq x, \, \mathrm{ord}(2)=p-1, \\ \mathrm{ord}(2) \neq p(p-1)}} \left( 1 - \frac{1}{p^2} \right) \prod_{\substack{p \leq x, \\ \mathrm{ord}(2)=p-1}} \left( 1 - \frac{1}{p} \right)^{-1}$$

$$= \frac{1}{2} \left( \frac{e^{\gamma_2 - \gamma\alpha_2}}{\Gamma(\alpha_2)} + o(1) \right) \frac{x}{(\log x)^{1-\alpha_2}} \prod_{\substack{p \leq x, \, \mathrm{ord}(2) = p-1, \\ \mathrm{ord}(2) \neq p(p-1)}} \left( 1 - \frac{1}{p^2} \right)$$

$$= \frac{1}{2} \left( \frac{1}{\Gamma(\alpha_2)} + o(1) \right) \frac{x}{(\log x)^{1-\alpha_2}} \prod_{\substack{\mathrm{ord}(2) = p-1, \\ \mathrm{ord}(2) \neq p(p-1)}} \left( 1 - \frac{1}{p^2} \right),$$

where $\alpha_2 > 0$ is the Artin constant for the primitive root 2, refer to Chapter 11 or [HY67, p. 220], and $\gamma_2 = \gamma\alpha_2 > 0$ is the Artin-Euler constant, see Lemma 19.3 for detail. Lastly, the convergent product was replaced with the approximation

$$\prod_{\substack{p \leq x, \, \mathrm{ord}(2) = p-1, \\ \mathrm{ord}(2) \neq p(p-1)}} \left( 1 - \frac{1}{p^2} \right) = \prod_{\substack{\mathrm{ord}(2) = p-1, \\ \mathrm{ord}(2) \neq p(p-1)}} \left( 1 - \frac{1}{p^2} \right) + O\left( \frac{1}{x} \right). \tag{23.9}$$

This yields line 3 above.  ∎

The subset of primes $\{ p \in \mathbb{P} : \mathrm{ord}(2) = p-1 \text{ and } \mathrm{ord}(2) \neq p(p-1) \}$ is the same as the subset of Wieferich primes. These primes are usually characterized by the congruence $\{ p \in \mathbb{P} : 2^{p-1} \equiv 1 \bmod p^2 \} = \{ 1093, \ 3511, \ \dots \}.$

The graphs display the index $\mathrm{ord}(2)/\lambda(n)$ of the integer 2 in the finite ring $\mathbb{Z}/n\mathbb{Z} = \{ 0, \ 1, \ 2, \ 3, \ \dots, \ n-1 \}$, two ranges of integers are demonstrated. Any even integers seems to have somewhat similar patterns.

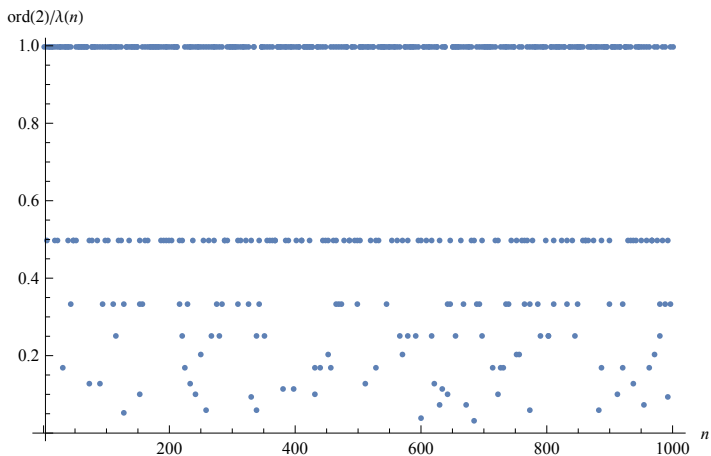

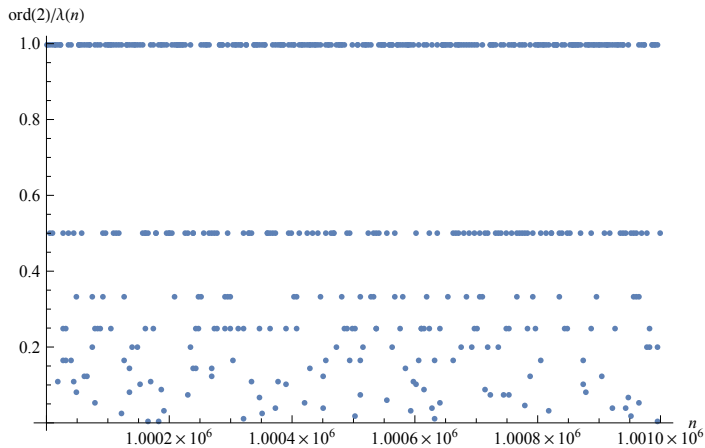

## 23.2. Subset of Integers With Fixed Primitive Root *u*

An asymptotic counting formula $N_u(x) = \#\{ n \leq x : \operatorname{ord}(u) = \lambda(n) \}$ for the subset of integers $\mathcal{N}_u = \{ n \in \mathbb{N} : \operatorname{ord}(u) = \lambda(n) \}$ is proved here. By Theorem 12.1, this analysis is unconditional. This theorem states that $u \neq \pm 1$, $v^2$ is a primitive root for infinitely many primes. The subset of primes $\mathcal{P}_u = \{ p \in \mathbb{P} : \operatorname{ord}(u) = p - 1 \} \subset \mathbb{P}$ of density $\alpha_u = \delta(\mathcal{P}_u) > 0$, which consists of all the primes with a fixed primitive root $u \in \mathbb{Z}$, is utilized to generate the subset of composite integers $\mathcal{N}_u$, which have $u$ as a primitive root. The associated constant

$$C_u = \frac{\varphi(u)}{u} \, \frac{1}{\Gamma(\alpha_u)} \prod_{\operatorname{ord}(u) = p-1, \, \operatorname{ord}(u) \neq p(p-1)} \left( 1 - \frac{1}{p^2} \right), \tag{23.1}$$

where $\alpha_u$ is Artin constant, see Chapter 11, and [HY67, p. 220] for the formula, is not completely determined since it has a dependence on the Abel primes $u^{p-1} \equiv 1 \bmod p^2$, but can be approximated. Very likely, this constant is a rational multiple of the real number $1 / \Gamma(\alpha_u)$.

**Theorem** 23.2. Let $u \neq \pm 1$, $v^2$ be a fixed integer. Then, the counting function for number of integers $n \leq x$, $\gcd(u, n) = 1$, with a primitive root mod $n$, has the asymptotic formulae

(i) $\displaystyle N_u(x) = \frac{\varphi(u)}{u} \left( \frac{1}{\Gamma(\alpha_u)} + o(1) \right) \frac{x}{(\log x)^{1-\alpha_u}} \prod_{\operatorname{ord}(u) = p-1, \, \operatorname{ord}(u) \neq p(p-1)} \left( 1 - \frac{1}{p^2} \right),$

(ii) $\displaystyle N_u(x) \gg \frac{x}{(\log x)^{1-\alpha_u}},$

$$\tag{23.2}$$

where the index of the product ranges over the Abel primes $u^{p-1} \equiv 1 \bmod p^2$, $\gcd(u, p) = 1$, and $\alpha_u > 0$, $\gamma_u > 0$ are constants, for all large numbers $x \geq 1$.

**Proof**: Start with Wirsing formula with $\tau = \alpha_u$, Lemma 20.1, and replace the characteristic function $f(n)$ of primitive roots in the finite ring $\mathbb{Z}/n\mathbb{Z}$, see Lemma 7.4, to produce

$$\sum_{n \leq x} f(n) = \left( \frac{1}{e^{\gamma \tau}\,\Gamma(\tau)} + o(1) \right) \frac{x}{\log x} \prod_{p \leq x} \left( 1 + \frac{f(p)}{p} + \frac{f(p^2)}{p^2} + \cdots \right)$$

$$= \frac{\varphi(u)}{u} \left( \frac{1}{e^{\gamma \alpha_u}\,\Gamma(\alpha_u)} + o(1) \right) \frac{x}{\log x} \prod_{\substack{p \leq x,\ \mathrm{ord}(u) = p-1, \\ \mathrm{ord}(u) \neq p(p-1)}} \left( 1 + \frac{1}{p} \right) \prod_{\substack{p \leq x, \\ \mathrm{ord}(u) = p(p-1)}} \left( 1 + \frac{1}{p} + \frac{1}{p^2} + \cdots \right),$$

(23.3)

where the index $p \in \mathcal{P}_u$ runs over the primes such that $\gcd(u,\,p) = 1$, and $u$ is a primitive root modulo $p$ or modulo $p^2$. The ratio $\varphi(u)/u$ accounts for the condition $\gcd(u,\,n) = 1$. The corrected product, see Section 22.3, yields

$$\sum_{n \leq x} f(n) = \frac{\varphi(u)}{u} \left( \frac{1}{e^{\gamma \alpha_u}\,\Gamma(\alpha_u)} + o(1) \right) \frac{x}{\log x} \prod_{\substack{p \leq x,\ \mathrm{ord}(u) = p-1, \\ \mathrm{ord}(u) \neq p(p-1)}} \left( 1 - \frac{1}{p^2} \right) \prod_{\substack{p \leq x, \\ \mathrm{ord}(u) = p-1}} \left( 1 - \frac{1}{p} \right)^{-1}$$

$$= \frac{\varphi(u)}{u} \left( \frac{e^{\gamma_u - \gamma \alpha_u}}{\Gamma(\alpha_u)} + o(1) \right) \frac{x}{(\log x)^{1 - \alpha_u}} \prod_{\substack{p \leq x,\ \mathrm{ord}(u) = p-1, \\ \mathrm{ord}(u) \neq p(p-1)}} \left( 1 - \frac{1}{p^2} \right)$$

(23.4)

$$= \frac{\varphi(u)}{u} \left( \frac{1}{\Gamma(\alpha_u)} + o(1) \right) \frac{x}{(\log x)^{1 - \alpha_u}} \prod_{\substack{\mathrm{ord}(u) = p-1, \\ \mathrm{ord}(u) \neq p(p-1)}} \left( 1 - \frac{1}{p^2} \right),$$

where $\alpha_u > 0$ is the Artin constant for the primitive root $u$, Chapter 11 or [HY67, p. 220], and $\gamma_u = \gamma \alpha_u > 0$ is the *Artin-Euler* constant, see Lemma 19.3 for detail. Lastly, the convergent product was replaced with the approximation

$$\prod_{\substack{p \leq x,\ \mathrm{ord}(u) = p-1, \\ \mathrm{ord}(u) \neq p(p-1)}} \left( 1 - \frac{1}{p^2} \right) = \prod_{\substack{\mathrm{ord}(u) = p-1, \\ \mathrm{ord}(u) \neq p(p-1)}} \left( 1 - \frac{1}{p^2} \right) + O\left( \frac{1}{x} \right).$$

(23.5)

This yields line 3 above. ∎

The subset of primes $\{\, p \in \mathbb{P} : \mathrm{ord}(u) \neq p(p-1) \,\}$ is the same as the subset of Abel primes $\{\, p \in \mathbb{P} : u^{p-1} \equiv 1 \bmod p^2 \,\}$.

The graphs display the index of the integer 3 in the finite ring $\mathbb{Z}/n\mathbb{Z} = \{\, 0,\ 1,\ 2,\ 3,\ \ldots,\ n-1 \,\}$, two ranges of integers are demonstrated.

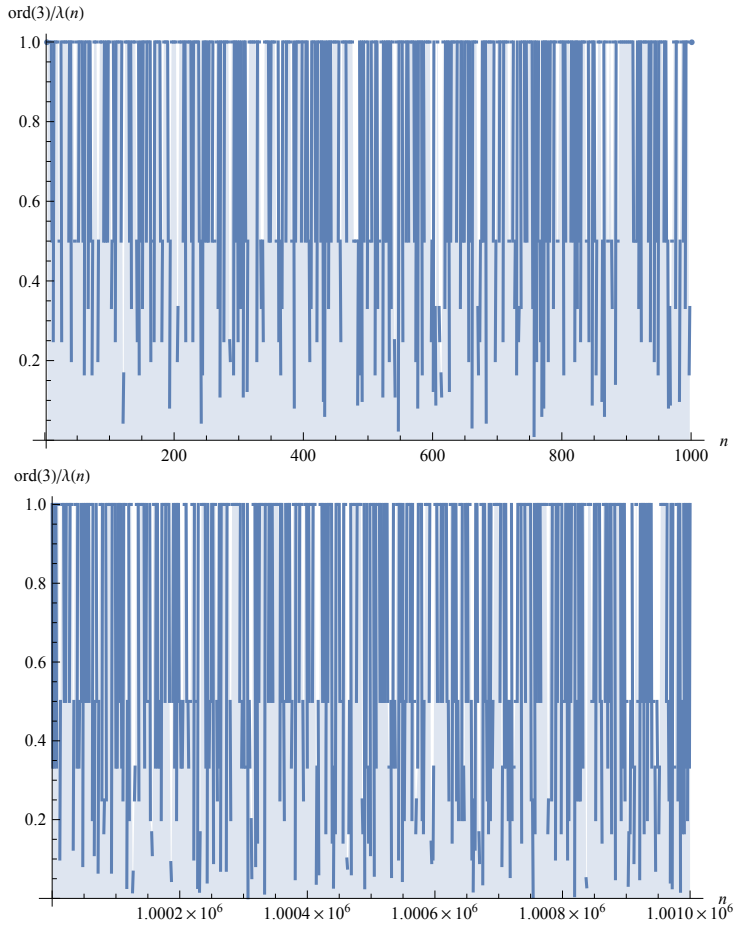

## 23.3 Problems And Exercises

**Problem 23.1.** Compute 24 digits approximation of the constant

$$\frac{1}{2\,\Gamma(\alpha_2)} \prod_{\mathrm{ord}(2)=p-1,\,\mathrm{ord}(2)\neq p(p-1)} \left(1 - \frac{1}{p^2}\right) = \frac{1}{2\,\Gamma(\alpha_2)}\left(1 - \frac{1}{1093^2}\right)\left(1 - \frac{1}{3511^2}\right)\cdots, \tag{23.6}$$

where $\alpha_2 = 0.373955813619202\ldots$ is Artin constant, and the product is over the Wieferich primes $2^{p-1} \equiv 1 \bmod p^2$.

**Problem 23.2.** Compute 24 digits approximation of the constant

$$\frac{1}{2\,\Gamma(\alpha_3)} \prod_{\mathrm{ord}(3)=p-1,\,\mathrm{ord}(3)\neq p(p-1)} \left(1 - \frac{1}{p^2}\right) = \frac{2}{3\,\Gamma(\alpha_3)}\left(1 - \frac{1}{11^2}\right)\left(1 - \frac{1}{1\,006\,003^2}\right)\cdots, \tag{23.7}$$

where $\alpha_3 = 0.373955813619202\ldots$ is Artin constant, and the product is over the Wieferich primes $3^{p-1} \equiv 1 \bmod p^2$.

**Problem 23.3.** Compute 24 digits approximation of the constant

$$\frac{4}{5\,\Gamma(\alpha_5)}\prod_{\mathrm{ord}(5)=p-1,\,\mathrm{ord}(5)\neq p(p-1)}\left(1-\frac{1}{p^2}\right)=\frac{4}{5\,\Gamma(\alpha_2)}\left(1-\frac{1}{2^2}\right)\left(1-\frac{1}{20\,771^2}\right)\cdots,\tag{23.8}$$

where $\alpha_5=(20\,/\,19)\cdot 0.373955813619202\ldots$ is Artin constant, and the product is over the Wieferich primes $5^{p-1}\equiv 1\bmod p^2$.

**Problem 23.4.** Let $\alpha_u>0$ be the Artin constant attached to the primitive root $u\in\mathbb{Z}$, and let $\Gamma(s)=\int_0^\infty t^{s-1}\,e^{-s\,t}\,d\,t,\ s\in\mathbb{C}$, be the gamma function. Are the numbers $\Gamma(\alpha_2),\ \Gamma(\alpha_3),\ \Gamma(\alpha_5),\ \Gamma(\alpha_6),\ \ldots$ irrationals?

**Problem 23.5.** The $L$-series for the primitive root $2\bmod n$ has the form

$$L_2(s,f)=\sum_{n\geq 1}\frac{f(n)}{n^s}=\prod_{\mathrm{ord}(2)=p-1,\,\mathrm{ord}(2)\neq p(p-1)}\left(1-\frac{1}{p^{2\,s}}\right)\prod_{\mathrm{ord}(2)=p-1}\left(1-\frac{1}{p^s}\right)^{-1},\tag{23.9}$$

where $s\in\mathbb{C}$ is a complex number, and $f(n)$ is the characteristic function of primitive root $2\bmod n$. Show that the product

$$\prod_{p>2,\,\mathrm{ord}(2)\neq p(p-1)}\left(1-\frac{1}{p^{2\,s}}\right)=\left(1-\frac{1}{1091^{2\,s}}\right)\left(1-\frac{1}{3511^{2\,s}}\right)\prod_{p>3511,\,\mathrm{ord}(2)=p-1,\,\mathrm{ord}(2)\neq p(p-1)}\left(1-\frac{1}{p^{2\,s}}\right)\tag{23.10}$$

has finitely many factors, this is implied by the unknown distribution of Wieferich primes.

**Problem 23.6.** Let $x\geq 1$ be a large number. Show that

$$\prod_{\substack{p\leq x,\ \mathrm{ord}(u)=p-1,\\ \mathrm{ord}(u)\neq p(p-1)}}\left(1-\frac{1}{p^2}\right)=\prod_{\substack{\mathrm{ord}(u)=p-1,\\ \mathrm{ord}(u)\neq p(p-1)}}\left(1-\frac{1}{p^2}\right)+O\!\left(\frac{1}{x}\right).\tag{23.11}$$

## Chapter 24

# Asymptotic Formulae For Quadratic, Cubic, And Quartic Nonresidues

## 24.1 Subset of Integers With Fixed Quadratic Nonresidues

The quadratic nonresidues has the asymptotic counting formula $\mathcal{K}_{2,r}(x) = \# \{ n \leq x : \text{ord}(r) = 2 \}$. This function coincides with the discrete measure of the subset of integers $\mathcal{K}_{2,r} = \{ n \in \mathbb{N} : \text{ord}(r) = 2 \}$ such that $1 \neq r \in \mathbb{Z}$ is a nontrivial second power residue.

**Theorem 24.1.** A fixed integer $r \neq \pm 1$, $s^2$ is quadratic nonresidue mod $n$ for infinitely many composite integers. More-over, number of integers $n \leq x$ such that $r$ is a quadratic nonresidue mod $n$ has the asymptotic formula

(i) $\mathcal{K}_{2,r}(x) = \dfrac{\varphi(r)}{r} \dfrac{1}{\sqrt{\pi}} \dfrac{x}{(\log x)^{1/2}} + O\left(\dfrac{x}{(\log x)^{3/2}}\right),$

(ii) $\mathcal{K}_{2,r}(x) \gg \dfrac{\varphi(r)}{r} \dfrac{x}{(\log x)^{1/2}} .$

$$(24.1)$$

**Proof**: Start with Wirsing formula, Lemma 20.1, and replace the characteristic function $f_2(n)$ of quadratic nonresidues in the finite rings $\mathbb{Z} / n \, \mathbb{Z}$, Lemma 21.2, to produce

$$
\begin{aligned}
\sum_{n \leq x} f(n) &= \left( \frac{1}{e^{\gamma \tau} \, \Gamma(\tau)} + o(1) \right) \frac{x}{\log x} \prod_{p \leq x} \left( 1 + \frac{f(p)}{p} + \frac{f(p^2)}{p^2} + \cdots \right) \\
&= \frac{\varphi(u)}{u} \left( \frac{1}{e^{\gamma/2} \, \Gamma(1/2)} + o(1) \right) \frac{x}{\log x} \prod_{\text{ord}(u)=2, \, p \leq x} \left( 1 + \frac{1}{p} + \frac{1}{p^2} + \cdots \right),
\end{aligned}
$$

$$(24.2)$$

where the product index $p$ runs over the primes such that $\gcd(r, p) = 1$, and $u$ is a quadratic nonresidue modulo $p^k$, $k \geq 1$. Applying Lemma 19.4, yields

$$
\log \prod_{p \leq x, \, \text{ord}(r)=2} \left( 1 + \frac{1}{p} + \frac{1}{p^2} + \cdots \right) = \log \prod_{p \leq x, \, \text{ord}(r)=2} \left( 1 - \frac{1}{p} \right)^{-1}
$$

$$(24.3)$$

$$= \frac{1}{2} \log\log x + \gamma_r + O\left(\frac{1}{\log x}\right),$$

where $\gamma_2 = \gamma / 2$ is a constant, since the density of primes $\{\, p \in \mathbb{P} : \mathrm{ord}(r) = 2 \,\}$ with a fixed quadratic nonresidue is $1 / 2$. This follows from the quadratic reciprocity law:

$$\left(\frac{r}{p}\right)\left(\frac{p}{r}\right) = (-1)^{\frac{p-1}{2}\cdot\frac{r-1}{2}}, \qquad (24.4)$$

see [AP76] for details the quadratic reciprocity law. Lastly, the gamma value $\Gamma(1/2) = \sqrt{\pi}$. $\blacksquare$

## 24.2. Subset of Integers With Fixed Cubic Nonresidues

Let $\mathcal{P}_{3,r} = \{\, p \in \mathbb{P} : \mathrm{ord}(r) = 3 \,\}$ be the subset of primes with a fixed cubic nonresidue, and let its density be defined by $\kappa_{3,r} = \delta(\mathcal{P}_{3,r}) > 0$. The information on this subset is used to obtain information on the subset of integers $\mathcal{K}_{3,r} = \{\, n \in \mathbb{N} : \mathrm{ord}(r) = 3 \,\}$ with a fixed cubic nonresidue $r \in \mathbb{Z}$. The corresponding counting function is denoted by $K_{3,r}(x) = \{\, n \leq x : \mathrm{ord}(r) = 3 \,\}$.

**Theorem 24.2.** A fixed integer $r \neq \pm 1$, $s^3$ is cubic nonresidue mod $n$ for infinitely many composite integers. Moreover, number of integers $n \leq x$ such that $r$ is a cubic nonresidue mod $n$ has the asymptotic formula

$$K_{3,r}(x) = C_3 \, \frac{\varphi(r)}{r} \, \frac{x}{(\log x)^{1/2}} + O\left(\frac{x}{(\log x)^{3/2}}\right), \qquad (24.5)$$

where $C_3 > 0$ is a constant.

## 24.3. Subset of Integers With Fixed Quartic Nonresidues

Let $\mathcal{P}_{4,r} = \{\, p \in \mathbb{P} : \mathrm{ord}(r) = 4 \,\}$ be the subset of primes with a fixed cubic nonresidue, and let its density be defined by $\kappa_{4,r} = \delta(\mathcal{P}_{4,r}) > 0$. The information on this subset is used to obtain information on the subset of integers $\mathcal{K}_{4,r} = \{\, n \in \mathbb{N} : \mathrm{ord}(r) = 4 \,\}$ with a fixed cubic nonresidue $r \in \mathbb{Z}$. The corresponding counting function is denoted by $K_{4,r}(x) = \{\, n \leq x : \mathrm{ord}(r) = 4 \,\}$.

**Theorem 24.3.** A fixed integer $r \neq \pm 1$, $s^2$ is quartic nonresidue mod $n$ for infinitely many composite integers. Moreover, number of integers $n \leq x$ such that $r$ is a quartic nonresidue mod $n$ has the asymptotic formula

$$\mathcal{K}_{4,r}(x) = C_4 \, \frac{\varphi(r)}{r} \, \frac{x}{(\log x)^{1/2}} + O\left(\frac{x}{(\log x)^{3/2}}\right). \qquad (24.6)$$

where $C_4 > 0$ is a constant.

## 24.4 Problems And Exercises

**Problem 24.1.** Let $\mathcal{P}_{m,r}(x) = \#\{\, p \leqslant x : \operatorname{ord}(r) = m \,\}$ be the counting function for the subset of primes $\mathcal{P}_{m,r} = \{\, p \in \mathbb{P} : \operatorname{ord}(r) = m \,\}$ such that $r \in \mathbb{Z}$ is an mth power residue. Compute an asymptotic formula for $\mathcal{P}_{m,r}(x)$.

**Problem 24.2.** Let $\mathcal{K}_{m,r}(x) = \#\{\, n \leqslant x : \operatorname{ord}(r) = m \,\}$ be the counting function for the subset of integers $\mathcal{K}_{m,r} = \{\, n \in \mathbb{N} : \operatorname{ord}(r) = m \,\}$ such that $r \in \mathbb{Z}$ is an mth power residue. Compute an asymptotic formula for $\mathcal{K}_{m,r}(x)$.

**Problem 24.3.** Let $A_\mu(x) = \#\{\, n \leqslant x : p \mid n \Rightarrow \mu(p-1) = \pm 1 \,\}$ be the counting function for the subset of integers $\mathcal{A}_\mu = \{\, n \in \mathbb{N} : p \mid n \Rightarrow \mu(p-1) = \pm 1 \,\}$. Compute an asymptotic formula for $A_\mu(x)$. Hint: Apply Wirsing formula, and Lemma 5.12.

# Chapter 25

## The $L$-Functions Of Fixed Primitive Roots

Let $\mathcal{A} \subset \mathbb{P}$ be a subset of primes of density $\alpha = \delta(\mathcal{A}) > 0$, and let $\mathcal{B} = \{ n \in \mathbb{N} : p \mid n \Rightarrow p \in \mathcal{A} \}$. Then, the associated $L$-function

$$L(s, \mathcal{A}) = \sum_{n \in \mathcal{B}} \frac{1}{n^s} = \prod_{p \in \mathcal{A}} \left( 1 - \frac{1}{p^s} \right)^{-1}, \tag{25.1}$$

where $s \in \mathbb{C}$ is a complex number, is absolutely convergent on the half plane $\{ s \in \mathbb{C} : \mathcal{R}e(s) > 1 \}$, and it has a pole at $s = 1$. The logarithm derivative of the associated $L$-function is

$$-\frac{L'(s, \mathcal{A})}{L(s, \mathcal{A})} = \sum_{n \in \mathcal{B}} \frac{\Lambda_{\mathcal{A}}(n)}{n^s}, \tag{25.2}$$

where the vonMangold function is defined by

$$\Lambda_{\mathcal{A}}(n) = \begin{cases} \log p & \text{if } n = p^k, \, k \geq 1, \text{ and } p \in \mathcal{A}, \\ 0 & \text{if } n \neq p^k, \, k \geq 1, \text{ or } p \notin \mathcal{A}. \end{cases} \tag{25.3}$$

The constraints forced on some subsets of primes induce certain structures on the associated $L$-functions. The existence of Abel-Wieferich primes, id est, primes $p \geq 2$ for which $\text{ord}(u) = p - 1$, but $\text{ord}(u) \neq p(p - 1)$, induces a complex structure on the $L$-series of the fixed primitive root $u \neq \pm 1, v^2$. The function in question is

$$L(s, f) = \sum_{n \geq 1} \frac{f(n)}{n^s} = \prod_{\text{ord}(u) = p-1, \, \text{ord}(u) \neq p(p-1)} \left( 1 - \frac{1}{p^{2s}} \right) \prod_{\text{ord}(u) = p-1} \left( 1 - \frac{1}{p^s} \right)^{-1}, \tag{25.4}$$

where $f(n)$ is the characteristic function of the primitive root $u \bmod n$, see Lemma 21.3, and $s \in \mathbb{C}$ is a complex number.

This is akin to the structure of the $L$-series attached to some modular forms. Recall that the $L$-series of an algebraic curve

$C : f(x, y) = 0$ over the integers $\mathbb{Z}$, and genus $g = 1$, is defined by

$$L(s, C) = \sum_{n \geq 1} \frac{a(n)}{n^s} = \prod_{p \mid N} \left(1 - \frac{1}{p^s}\right)^{-1} \prod_{p \nmid N} \left(1 - \frac{a(p)}{p^s} + \frac{1}{p^{2s-1}}\right)^{-1}, \tag{25.5}$$

where $a(n)$ counts the number of rational points over the finite ring $\mathbb{Z} / n \mathbb{Z}$, and $N \geq 1$ is the conductor, refer to the literature for the precise information.

**Theorem 25.1.** Let $\mathcal{P}_u \subset \mathbb{P}$ be a subset of primes with a fixed primitive root $u$, and let $\alpha_u = \delta(\mathcal{P}_u) > 0$ be its density. Then, the associated $L$-function

$$L(s, f) = \sum_{n \geq 1} \frac{f(n)}{n^s} \tag{25.6}$$

where $s \in \mathbb{C}$ is a complex number, is absolutely convergent on the half plane $\{ s \in \mathbb{C} : \mathcal{R}e(s) > 1 \}$, and it has a pole at $s = 1$.

**Proof**: Apply Delange theorem, confer Lemma 20. 3. ∎

The logarithm derivative of the associated $L$-function is

$$-\frac{L'(s, f)}{L(s, f)} = \sum_{n \geq 1} \frac{\Lambda_f(n)}{n^s}, \tag{25.7}$$

where the vonMangold function is defined by

$$\Lambda_f(n) = \begin{cases} \log p & \text{if } n = p^k, \, k \geq 1, \text{ and ord } u = p - 1, \\ 0 & \text{if } n \neq p^k, \, k \geq 1. \end{cases} \tag{25.8}$$

## 25.1 Problems And Exercises

**Problem 25.1.** Prove or disprove that the subset of integers $\mathcal{N}_5 = \{ n \in \mathbb{N} : \text{ord}(5) = \lambda(n) \}$ such that 5 is a primitive root mod $n$ has the smallest density

$$\delta_5 = \lim_{x \to \infty} \frac{\{ n \leq x : \text{ord}(5) = \lambda(n) \}}{x} > 0 \tag{25.9}$$

Hint: The product

$$\prod_{\text{ord}(5) = p-1, \, \text{ord}(5) \neq p(p-1)} \left(1 - \frac{1}{p^{2s}}\right) = \left(1 - \frac{1}{2^{2s}}\right) \prod_{\text{ord}(5) = p-1, \, \text{ord}(5) \neq p(p-1)} \left(1 - \frac{1}{p^{2s}}\right) \tag{25.10}$$

at $s = 1$, confer Theorem 23.2, and show that it has finitely many factors. This is implied by the unknown distribution of Abel-Wieferich primes $p \geq 2$ such that $5^{p-1} \equiv 1 \bmod p^2$.

**Problem 25.2.** Let $u \in \mathbb{Z}$ be a fixed nonsquare integer. Prove or disprove that the product

$$\prod_{\text{ord}(u)=p-1,\,\text{ord}(u)\neq p(p-1)} \left(1 - \frac{1}{p^{2s}}\right) \tag{25.11}$$

has finitely many factors, this is implied by the unknown distribution of Abel-Wieferich primes $u^{p-1} \equiv 1 \bmod p^2$.

## Chapter 26

# Class Numbers And Primitive Roots

Let $p \geq 2$ be a prime. The length of the $\ell$-adic expansion of the rational number $a / p = 0.\overline{x_1 \, x_2 \, ... \, x_m}$, where $a < p$, was the original motivation for study of primitive roots.

***Lemma 26.1.*** (Gauss) The decimal expansion $1 / p = 0.\overline{x_1 \, x_2 \, ... \, x_m}$ has maximal length $m = p - 1$ if and only if $p \geq 2$ is a prime and 10 is a primitive root mod $p$.

## 26.1 Recent Formulae

***Theorem 26.2.*** (Girstmair formula) Let $p = 4\,m + 3 \geq 7$ be a prime, and let $\ell \geq 2$ be a fixed primitive root mod $p$. Let

$$\frac{1}{p} = \sum_{n \geq 1} \frac{x_n}{\ell^n}, \tag{26.1}$$

where $x_i \in \{ 0, \, 1, \, 2, \, ..., \, \ell - 1 \}$ be the unique expansion. Let $h(-p) \geq 1$ be the class number of the quadratic field $\mathbb{Q}\left( \sqrt{-p} \right)$. Then

$$\sum_{1 \leq n \leq p-1} (-1)^n \, x_n = (\ell + 1) \, h(-p) \, . \tag{26.2}$$

The analysis of the Girstmair formula appears in [GR94], [GRS94]. A generalization to elements of any orders mod $p$ was later proved in [RT11].

***Theorem 26.3.*** ([RT11]) Let $p \geq 3$ be a prime, and $r \mid p - 1$. Suppose that a character $\chi$ of order $r$ is odd and that $\ell \geq 2$ has order $(p - 1) / r$ mod $p$. Let

$$\frac{1}{p} = \sum_{n \geq 1} \frac{x_n}{\ell^n}, \tag{26.3}$$

where $x_i \in \{ 0, \, 1, \, 2, \, ..., \, \ell - 1 \}$ be the unique $\ell$- expansion. Then

$$\sum_{1 \leq n \leq (p-1)/r} x_n = \frac{\ell - 1}{2} \, \frac{p - 1}{r} + \frac{\ell - 1}{r} \sum_{\mathrm{ord}\chi = r} B_{1, \chi} \tag{26.4}$$

where $B_{1,\chi}$ is the generalized Bernoulli number.

A similar relationship between the class number and the continued fraction of the number $\sqrt{p}$ was discovered earlier, see [HF73, p. 241].

**Theorem** 26.4. (Hirzebruch formula) Let $p = 4m + 3 \geq 7$ be a prime. Suppose that the class number of the quadratic field $\mathbb{Q}\left(\sqrt{p}\right)$ is $h(p) = 1$, and $h(-p) \geq 1$ is the class number of the quadratic field $\mathbb{Q}\left(\sqrt{-p}\right)$. Then

$$\sum_{1 \leq n \leq 2t} (-1)^n a_n = 3 h(-p),$$ (26.5)

where $\sqrt{p} = [a_0, \overline{a_1, a_2, \ldots, a_{2t}}], a_i \geq 1$, is the continued fraction of the real number $\sqrt{p}$.

Simple and straight Forward applications of these results are presented here. One of these applications is related to and reminiscent of the square root cancellation rule.

**Corollary** 26.5. Let $p = 4m + 3 \geq 7$ be a prime. Let $\ell \geq 2$ be a fixed primitive root mod $p$, and let $1/p = 0. \overline{x_1 x_2 \ldots x_m}$ be the $\ell$-adic expansion of the rational number $1/p$. Then, for any arbitrarily small number $\epsilon > 0$, the alternating sum satisfies

(i) $\displaystyle\sum_{1 \leq n \leq p-1} (-1)^n x_n \ll p^{1/2+\epsilon}$.

(ii) $\displaystyle\sum_{1 \leq n \leq p-1} (-1)^n x_n \gg p^{1/2-\epsilon}$. (26.6)

**Proof**: (i) By Theorem 26.3, and Theorem 9.1, the alternating sum has the upper bound

$$\sum_{1 \leq n \leq p-1} (-1)^n x_n = (\ell + 1) h(-p)$$
$$\ll p^{\epsilon/2} h(-p)$$
$$\ll p^{1/2+\epsilon},$$ (26.7)

where $h(-p) \ll p^{1/2+\epsilon/2}$ is an upper bound of the class number of the quadratic field $\mathbb{Q}\left(\sqrt{-p}\right)$, with $\epsilon > 0$.

For (ii), apply Siegel theorem $L(1, \chi) \gg p^\epsilon$, see [MV07, p. 372]. ∎

The class numbers of quadratic fields have various formulations. Among these is the expression

$$h(-p) = \begin{cases} \dfrac{w\sqrt{p}}{2\pi}\, L(1,\chi), & d < 0, \\[2ex] \dfrac{\sqrt{p}}{\log \epsilon}\, L(1,\chi), & d > 0, \end{cases}$$

where $\chi \ne 1$ is the quadratic symbol mod $p$, $\epsilon = a + b\sqrt{d}$ is a fundamental unit, and $w = 6, 4,$ or $2$ according to $d = -3, -4,$ or $d < -4$ respectively, see [UW00, p. 28], [MV07, p. 391]. This is utilized below to exhibit two ways of computing the special values of the $L$-functions $L(s,\chi) = \sum_{n\ge 1}\chi(n)\, n^{-s}$ at $s = 1$.

**Corollary 26.6.** Let $p = 4m + 3 \ge 7$ be a prime. Let $\ell \ge 2$ be a fixed primitive root mod $p$. Let $1/p = 0.\overline{x_1 x_2 \ldots x_m}$ be the $\ell$-adic expansion of the rational number $1/p$, and let $\sqrt{p} = [a_0, \overline{a_1, a_2, \ldots, a_{2t}}]$, $a_i \ge 1$, be the continued fraction of the real number $\sqrt{p}$. Then,

(i) $L(1, \chi) = \dfrac{(\ell + 1)\pi}{\sqrt{p}} \displaystyle\sum_{1 \le n \le p-1} (-1)^n x_n$, for any $h(-p) \ge 1$.

(ii) $L(1, \chi) = \dfrac{\pi}{3\sqrt{p}} \displaystyle\sum_{1 \le n \le 2t} (-1)^n a_n$, if $h(p) = 1$, and $h(-p) \ge 1$. $\qquad$ (26.9)

(iii) $\displaystyle\sum_{1 \le n \le 2t} (-1)^n a_n = 3(\ell + 1) \sum_{1 \le n \le p-1} (-1)^n x_n$, if $h(p) = 1$, and $h(-p) \ge 1$.

## 26.2 Problems And Exercises

**Problem 26.1.** Let $\pi \in \mathcal{K} = \mathbb{Q}(\sqrt{d})$ be a prime, and let $\beta \ne \pm\alpha$, $\gamma^2$, $N(\alpha) = \pm 1$, be a fixed primitive root mod $\pi$. Let

$$\frac{1}{\pi} = \sum_{n \ge 1} \frac{x_n}{\beta^n}, \qquad (26.10)$$

where $x_i \in \mathbb{Z}[\sqrt{d}]/(\beta)$ be the unique expansion. Note: the norm of an element $\alpha = a + b\sqrt{d} \in \mathbb{Q}(\sqrt{d})$, in the quadratic field is given by $N(\alpha) = a^2 - d\, b^2$. Generalize the Girstmair formula to quadratic numbers fields: Prove some form of relation such as this: If $h \ge 1$ is the class number of the numbers field $\mathcal{L} = \mathbb{Q}(\sqrt{\pi})$. Then

$$\sum_{1 \le n \le p-1} (-1)^n x_n \overset{?}{=} (\beta + 1)\, h\,. \qquad (26.11)$$

**Problem 26.2.** Show that Corollary 26.65 implies that there exists infinitely many primes $p = 4m + 3 \ge 7$ for which the class number of the quadratic field $\mathbb{Q}(\sqrt{p})$ is $h(p) = 1$. If not state why not. For example, in Corollary 26.6-iii, it is show

that

$$\sum_{1 \le n \le 2\,t} (-1)^n\, a_n = 3\,(\ell + 1) \sum_{1 \le n \le p-1} (-1)^n\, x_n\,, \qquad\qquad \text{if } h(p) = 1, \text{ and } h(-p) \ge 1. \tag{26.12}$$

**Problem 26.3.** The Cohen-Lenstra heuristic claims that the probability that the class number of the quadratic field $\mathbb{Q}\!\left(\sqrt{p}\,\right)$ is $h(p) = 1$ is $\Pr(h(p) = 1) \approx .754458173$, the derivation is given in [JW09, p. 167], and similar sources. Show that it has some direct correlation to the Artin constant:

$$\prod_{p \le x}\left(1 - \frac{1}{p(p-1)}\right) + \epsilon = \tag{26.13}$$

$$0.373955\ldots + \epsilon \; \overset{?}{=} \; \frac{1}{2}\Pr(h(p) = 1) \approx \frac{.754458173}{2} = 0.377229 \quad \text{for some constant } \epsilon \in \mathbb{R}.$$

## Appendices

# Appendix A

***Theorem A*.1** (Abel summation formula) Let $x \geq 1$ be a large number. Let $\{\, a_n : n \in \mathbb{N} \,\} \subset \mathbb{C}$ be a sequence of complex numbers, and let $A(x) = \sum_{n \leq x} a_n$. If a function $f : [1, \infty] \longrightarrow \mathbb{C}$ has continuous derivative on the interval $(1, \infty)$, then

$$\sum_{n \leq x} a_n \, f(n) = A(x) \, f(x) - \int_1^x A(t) \, f'(t) \, d\,t \,. \tag{27.1}$$